\documentclass[brochure,12pt,french]{smfbourbaki}

\usepackage[T1]{fontenc}
\usepackage{lmodern,amssymb,bm}
\usepackage{graphicx}
\usepackage{tikz}
\usepackage{braids}
\usetikzlibrary{matrix,arrows,calc,fit,cd,positioning,intersections,arrows.meta,braids}

\usepackage{babel}

\usepackage[utf8]{inputenc}

\usepackage{url}

\usepackage{microtype} 

\usepackage[shortlabels]{enumitem}

\usepackage[
backend=bibtex,
style=authoryear, 
citestyle=authoryear-comp,
maxnames=7,
sortcites=false
]{biblatex}
\usepackage{csquotes}

\DeclareDelimFormat{nameyeardelim}{\addcomma\space}

\DeclareNameAlias{sortname}{given-family}

\renewcommand{\bibnamedash}{\leavevmode\raise3pt\hbox to3em{\hrulefill}\space}

\AtEveryBibitem{%
  \clearfield{issn} 
  \clearfield{isbn} 
  \clearfield{doi} 

  \ifentrytype{online}{}{
    \clearfield{url}
  }
}

\renewbibmacro{in:}{%
    \ifentrytype{article}{}{\printtext{\bibstring{in}\intitlepunct}}}

\DeclareFieldFormat[article,periodical,inreference]{number}{\mkbibparens{#1}}
\DeclareFieldFormat[article,periodical,inreference]{volume}{\mkbibbold{#1}}
\renewbibmacro*{volume+number+eid}{%
    \printfield{volume}%
    \setunit*{\addthinspace}
    \printfield{number}%
    \setunit{\addcomma\space}%
    \printfield{eid}}

\DeclareFieldFormat[article,inbook,incollection]{title}{\enquote{#1}\addcomma} 

\date{Juin 2022}
\bbkannee{74\textsuperscript{e} année, 2021--2022}  
\bbknumero{1195}                                      

\title{La conjecture du $K(\pi,1)$ pour les groupes d'Artin affines}
\subtitle{d'après Giovanni Paolini et Mario Salvetti}

\author{Thomas Haettel}
\address{IMAG, Univ Montpellier, CNRS, France\\ Place Eugène Bataillon\\ 34090 Montpellier, France}
\email{thomas.haettel@umontpellier.fr}

\newcommand{\bdf}{\begin{defi}}
\newcommand{\edf}{\end{defi}}

\newcommand{\bpro}{\begin{prop}}
\newcommand{\epro}{\end{prop}}

\newcommand{\bthm}{\begin{theo}}
\newcommand{\ethm}{\end{theo}}

\newcommand{\bex}{\begin{exem}}
\newcommand{\eex}{\end{exem}}

\newcommand{\bp}{\begin{proof}}
\newcommand{\ep}{\end{proof}}

\newcommand{\blem}{\begin{lemm}}
\newcommand{\elem}{\end{lemm}}

\newcommand{\beq}{\begin{eqnarray*}}
\newcommand{\eeq}{\end{eqnarray*}}

\newcommand{\beqn}{\begin{equation}}
\newcommand{\eeqn}{\end{equation}}

\newcommand{\ben}{\begin{enumerate}}
\newcommand{\een}{\end{enumerate}}
\newcommand{\bit}{\begin{itemize} \renewcommand{\labelitemi}{\textemdash} \renewcommand{\labelitemii}{$\star$}}
\newcommand{\eit}{\end{itemize}}

\newcommand{\R}{\mathbb{R}}

\newcommand{\N}{\mathbb{N}}
\newcommand{\Z}{\mathbb{Z}}
\newcommand{\C}{\mathbb{C}}

\renewcommand{\O}{\mathbb{O}}

\newcommand{\D}{\mathbb{D}}
\renewcommand{\SS}{\mathbb{S}}

\renewcommand{\L}{\mathbb{L}}

\newcommand{\ov}{\overline}

\renewcommand{\tilde}{\widetilde}

\newcommand{\Cay}{\operatorname{Cay}}

\newcommand{\codim}{\operatorname{codim}}

\renewcommand{\t}{ ^t\!}

\renewcommand{\dim}{\operatorname{dim}}

\newcommand{\Isom}{\operatorname{Isom}}

\newcommand{\Min}{\operatorname{Min}}

\newcommand{\Fix}{\operatorname{Fix}}
\newcommand{\Dep}{\operatorname{Dep}}

\newcommand{\eps}{\varepsilon}
\newcommand{\st}{\, | \,}

\newcommand{\ra}{\rightarrow}

\newcommand{\f}{\frac}

\renewcommand{\geq}{\geqslant}
\renewcommand{\leq}{\leqslant}

\newcommand{\gothique}{\mathfrak}

\renewcommand{\sl}{\gothique{sl}}

\newcommand{\<}{\langle}
\renewcommand{\>}{\rangle}

\newcommand{\mk}{\medskip}

\newtheorem*{thms}{\theoname}
\newtheorem*{conjs}{\conjname}

\begin{document}

\maketitle

\tableofcontents

\section*{Introduction}

Les groupes de tresses ont été définis formellement par 
\textcite{artin_original}, et sont depuis un objet d'étude fascinant, grâce à leurs multiples définitions et connexions avec d'autres domaines (nous renvoyons le lecteur à \textcite{kassel_turaev} et \textcite{farb_margalit} pour des introductions aux groupes de tresses). L'une des manières de définir le groupe de tresses à $n$ brins est de le voir comme le groupe fondamental de l'espace de configuration de $n$ points du plan
$$Y = \{F \subset \C, |F|=n\}.$$
Ainsi, si on considère l'arrangement d'hyperplans $\{z_i=z_j\}_{i \neq j}$ dans $\C^n$, alors l'espace de configuration $Y$ peut aussi s'interpréter comme le complémentaire de cet arrangement d'hyperplans
$$\{z \in \C^n \st z_i \neq z_j \mbox{ si } i \neq j\},$$
quotienté par l'action naturelle du groupe symétrique $\frak{S}_n$ sur $\C^n$.

\begin{figure}
\begin{center}
\begin{tikzpicture}
\definecolor{col1}{rgb}{0,0,0.6}
\definecolor{col2}{rgb}{0,0,0.8}
\definecolor{col3}{rgb}{0.1,0.1,1}
\definecolor{col4}{rgb}{0.3,0.3,1}
\definecolor{col5}{rgb}{0.5,0.5,1}
\pic[
braid/.cd,
every strand/.style={ultra thick},
strand 1/.style={col1},
strand 2/.style={col2},
strand 3/.style={col3},
strand 4/.style={col4},
strand 5/.style={col5},
height=-0.8cm,]
{braid={|s_1-s_4| s_2  |s_1-s_4| s_3^{-1}}};
\end{tikzpicture}
\caption{Une représentation d'une tresse à cinq brins.}
\label{fig:tresse}
\end{center}
\end{figure}
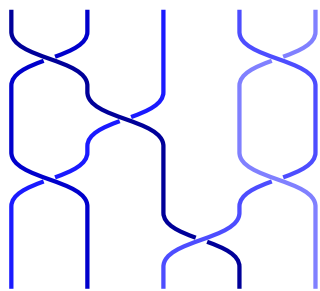

Chaque tresse peut se représenter comme un lacet dans l'espace de configuration de $n$ points du plan, c'est-à-dire comme $n$ chemins dans le plan ne s'intersectant pas à chaque instant. En représentant le temps par exemple de haut en bas, on peut ainsi obtenir un dessin classique de tresse comme dans la figure~\ref{fig:tresse}.

\mk

Si l'on considère un groupe de Coxeter~$W$ quelconque, il possède une action par réflexions sur un cône d'un espace vectoriel appelé cône de Tits \parencite{bourbaki_coxeter} généralisant l'action linéaire de~$\frak{S}_n$ sur~$\R^n$. On peut lui associer un espace topologique~$Y_W$ correspondant au complémentaire de l'arrangement des hyperplans de réflexion. L'exemple le plus simple est donné par le cas d'un groupe de Coxeter affine~$W$, agissant par réflexions sur $\R^n$. Dans ce cas, l'espace~$Y_W$ est le complémentaire dans~$\C^n$ de la réunion des complexifiés des hyperplans de réflexions, quotienté par~$W$. L'espace~$Y_W$ est également appelé l'espace de configuration de type~$W$, par analogie avec le cas des groupes de tresses.

\mk

Le groupe fondamental $G_W$ de $Y_W$ est appelé groupe d'Artin, ou groupe d'Artin--Tits, car il a été défini par 
\textcite{tits_normalisateurs}. Il possède une présentation très simple, analogue à la présentation d'Artin des groupes de tresses (voir dans la partie~\ref{sec:interval_garside_groups}). Cependant, hormis dans des cas très particuliers, la topologie de l'espace de configuration $Y_W$ et les propriétés algébriques du groupe d'Artin $G_W$ restent largement mystérieuses \parencite{godelle_paris,charney_problems,mccammond_mysterious}. Côté algébrique, on ne sait toujours pas si le groupe d'Artin $G_W$ est sans torsion, on ne connaît pas son centre, et on ne sait pas résoudre le problème du mot. Côté topologique, la grande question ouverte concernant cet espace est la suivante :

\begin{conjs}[Conjecture du $K(\pi,1)$]
L'espace de configuration $Y_W$ est un espace classifiant pour le groupe d'Artin $G_W$.
\end{conjs}

Rappelons que $Y_W$ est un espace classifiant pour $G_W$ si $\pi_1(Y_W)=G_W$ et si $Y_W$ est asphérique, c'est-à-dire que pour tout $i \geq 2$, nous avons $\pi_i(Y_W)=0$. Ceci est donc équivalent à la contractibilité du revêtement universel de $Y_W$.
Cette question a été étudiée, essentiellement dans le cas où $W$ est fini, par Arnol'd, Brieskorn, Deligne, Pham et Thom dans les années 1970 \parencite{brieskorn_kpi1,deligne_immeubles}. L'énoncé général de la conjecture du $K(\pi,1)$, pour un groupe de Coxeter quelconque, peut probablement être attribuée à Loojienga et à son élève en thèse van der Lek, qui a montré que le groupe fondamental de $Y_W$ est bien le groupe d'Artin $G_W$ \parencite{vanderlek}. Nous renvoyons le lecteur à \textcite{paris_kpi1} pour une présentation moderne de cette conjecture. Remarquons que si la conjecture du $K(\pi,1)$ est vraie, comme $Y_W$ est asphérique et de dimension finie, cela implique que $G_W$ est sans torsion, et cela implique également que le centre de $G_W$ est connu \parencite{jankiewicz_schreve}.

\mk

La conjecture du $K(\pi,1)$ a été résolue pour les groupes d'Artin de type sphérique, c'est-à-dire ceux dont le groupe de Coxeter est fini, par 
\textcite{deligne_immeubles}. Elle a par la suite été prouvée pour les groupes d'Artin affines de type $\tilde{A}_n$, $\tilde{C}_n$ \parencite{okonek}, $\tilde{B}_n$ \parencite{callegari_moroni_salvetti} et $\tilde{G_2}$ \parencite{charney_davis}. Les travaux remarquables de Paolini et Salvetti apportent une réponse positive unifiée pour tous les groupes d'Artin affines :

\begin{thms}[\cite{PS}]
La conjecture du $K(\pi,1)$ est vraie pour tous les groupes d'Artin affines.
\end{thms}

\mk

En dehors des cas sphériques et affines, la conjecture du $K(\pi,1)$ a été prouvée pour les groupes d'Artin de dimension au plus $2$ ou de type $FC$ \parencite{hendriks,charney_davis}, ainsi que pour les groupes d'Artin relativement extra-larges \parencite{juhasz}.

\mk

Remarquons que la stratégie de Deligne consiste à montrer que le revêtement universel de $Y_W$ a le type d'homotopie d'un complexe simplicial, appelé complexe de Deligne, dont il montre qu'il est contractile \parencite{deligne_immeubles}. Charney et Davis ont muni le complexe de Deligne d'une métrique CAT(0) dans certains cas, ce qui leur permet de montrer sa contractibilité \parencite{hendriks,charney_davis}. Cette stratégie de munir le complexe de Deligne d'une métrique à courbure négative est l'objet de plusieurs travaux récents \parencite{charney_deligne,godelle_cat0,boyd_charney_morris,haettel_kpi1,goldman_extralarge,morris_wright_parabolicFC}.

\mk

On peut également poser la question de l'asphéricité du complémentaire d'un arrangement d'hyperplans ne provenant pas d'un groupe de Coxeter, par exemple le complexifié d'un arrangement d'hyperplans réel quelconque, ou bien provenant d'un groupe de réflexions complexes. Il existe des arrangements d'hyperplans réels linéaires finis complexifiés dont le complémentaire n'est pas asphérique, voir par exemple \parencite{falk_kpi1}. Cependant, pour les groupes de réflexions complexes finis, le complémentaire de l'arrangement est toujours asphérique \parencite{bessis_complex}.

\mk

L'un des intérêts de la conjecture du $K(\pi,1)$ est qu'elle permet de calculer l'homologie et la cohomologie du groupe d'Artin $G_W$ à partir de l'espace de configuration $Y_W$, ce qui a déjà fait l'objet de nombreux travaux \parencite{cohen_braid,salvetti_homotopy,charney_davis_finitekpi1,deconcini_salvetti_artin,callegaro_milnor,callegaro_moroni_salvetti,paolini_salvetti_homology,paolini_homology}.

\mk

Une des raisons fondamentales pour lesquelles les groupes de tresses, et plus généralement les groupes d'Artin de type sphérique, sont bien compris est la notion de structure de Garside. Informellement, un groupe est dit de Garside si l'on peut trouver un élément particulier dont les diviseurs engendrent le groupe et forment un treillis \parencite{garside,dehornoy_paris,dehornoy_ENS,dehornoy_garside}. Un groupe de Garside a un problème du mot résoluble, et un espace classifiant combinatoire très simple.

\mk

Lorsque $W$~est un groupe de Coxeter fini, on peut considérer le relevé naturel dans~$G_W$ de l'élément le plus long du groupe de Coxeter~$W$ : cela définit une structure de Garside, appelée structure standard.

\mk

Cependant, lorsque le groupe de Coxeter $W$ est infini, il n'y a plus d'élément le plus long, et on ne connaît pas en général de structure de Garside pour $G_W$. En revanche, déjà pour les groupes d'Artin sphériques, il existe une autre structure de Garside appelée structure duale \parencite{birman_ko_lee,bessis}. L'idée est de considérer un élément de Coxeter à la place de l'élément le plus long de $W$. Cela revient, dans le cas du groupe de tresses à $n$ brins, à considérer un $n^{\text{ème}}$ de tour à la place d'un demi-tour. En termes de générateurs, cela revient à considérer toutes les transpositions, et non pas seulement les transpositions entre brins adjacents, voir la partie~\ref{sec:interval_garside_groups}.

\mk

Il s'avère qu'il existe également des structures de Garside duales pour des
groupes d'Artin non sphériques :  \textcite{digne_An,digne_Cn} a montré que
pour certains groupes d'Artin affines, le relevé d'un élément de Coxeter définit une
structure de Garside duale. Cependant, 
 \textcite{mccammond_dual_euclidian}
a montré que la propriété de treillis était fausse dans les autres cas.

\mk

Pour pallier ce manque, 
 \textcite{mccammond_sulway}
ont montré qu'il était possible d'agrandir tout groupe d'Artin affine en un groupe appelé groupe cristallographique tressé, qui possède lui une structure de Garside duale. Ceci leur a permis de construire, pour tout groupe d'Artin affine, un espace classifiant de dimension finie, mais avec cependant un nombre infini de cellules. L'une des conséquences des travaux de Paolini et Salvetti est l'amélioration suivante :

\begin{thms}[\cite{PS}]
Tout groupe d'Artin affine possède un espace classifiant fini.
\end{thms}

Nous allons donner un panorama de la preuve des résultats principaux de \textcite{PS}, que nous détaillerons dans la suite du texte.

\mk

Fixons un groupe de Coxeter~$W$ et un élément de Coxeter~$w$, c'est-à-dire un
produit des générateurs standard dans un ordre arbitraire. À partir de
l'intervalle $[1,w]$ des diviseurs de~$w$, on peut définir le groupe d'Artin
dual~$W_w$ associé à~$w$ (voir la partie~\ref{sec:dual_artin_groups}). Lorsque
$W$~est fini ou affine, le groupe d'Artin dual $W_w$ est isomorphe au groupe d'Artin~$G_w$ \parencite{bessis,mccammond_sulway}, mais c'est une question ouverte en général. Lorsque l'intervalle $[1,w]$ est un treillis, le groupe d'Artin dual~$W_w$ est un groupe de Garside, et possède un espace classifiant explicite~$K_W$, qui provient de la réalisation géométrique de cet ensemble ordonné~\mbox{$[1,w]$}.

\mk

Le premier argument important est de montrer que, même si l'intervalle $[1,w]$ n'est pas un treillis, le complexe d'intervalle $K_W$ reste un espace classifiant (infini) pour le groupe d'Artin dual $W_w$. Cette preuve repose sur les groupes cristallographiques tressés définis par 
\textcite{mccammond_sulway}, qui apparaissent comme produits amalgamés ayant pour facteur le groupe d'Artin $G_W$. Pour les définir, McCammond et Sulway ont étudié l'ensemble ordonné de toutes les isométries d'un espace euclidien, avec pour partie génératrice l'ensemble des réflexions. On renvoie aux parties~\ref{sec:dual_artin_groups}, \ref{sec:factorisation_euclidean_isometries} et \ref{sec:crystallographic_braid_groups} pour les détails.

\mk

La suite de la preuve de la conjecture du $K(\pi,1)$ consiste à définir un sous-complexe~$X'_W$ de~$K_W$ qui aura le même type d'homotopie que l'espace de configuration~$Y_W$. Pour cela, nous allons partir du complexe de Salvetti~$X_W$ \parencite{salvetti,paris_kpi1}, qui est un modèle combinatoire bien étudié de l'espace de configuration~$Y_W$, et qui peut se décrire par recollement de sous-complexes correspondant à des sous-groupes paraboliques sphériques. La version duale de cette construction fournit un complexe de Salvetti dual~$X'_W$, qui a le même type d'homotopie que le complexe de Salvetti standard~$X_W$, et qui se réalise naturellement comme sous-complexe de~$K_W$. Cette construction est détaillée dans la partie~\ref{sec:dual_model_configuration_space}.

\mk

La fin de la preuve consiste à montrer que le complexe d'intervalle~$K_W$ se rétracte par déformation forte sur son sous-complexe~$X'_W$. Tout d'abord, en étudiant l'action de l'élément de Coxeter~$w$ par conjugaison sur~$K_W$, on peut définir un sous-complexe fini~$K'_W$ de~$K_W$ contenant~$X'_W$, tel que $K_W$~se rétracte sur~$K'_W$. La preuve que $K'_W$~se rétracte sur~$X'_W$ est plus technique, et repose notamment sur la construction d'un étiquetage lexicographique sur l'intervalle $[1,w]$. Ceci a pour conséquence intéressante que l'intervalle $[1,w]$, aussi appelé ensemble des partitions non croisées affines, est décortiquable. La preuve de l'existence des deux rétractions de~$K_W$ sur~$K'_W$ puis~$X'_W$ repose sur la théorie de Morse discrète, via la construction de couplages acycliques. Les détails sont donnés dans les parties~\ref{sec:lexicographic_orderings}, \ref{sec:classifying_space_dual_artin}, \ref{sec:proof_Kpi1} et \ref{sec:discrete_morse_theory}.

\mk

Cette approche duale pour la conjecture du $K(\pi,1)$ est également présentée de manière synthétique dans \parencite{paolini_dualapproach}. En particulier, il est envisageable que cette stratégie puisse fournir une preuve de la conjecture du $K(\pi,1)$ pour d'autres groupes d'Artin : le cas des groupes d'Artin de rang $3$ est annoncé dans~\parencite[Theorem~6.1]{paolini_dualapproach}.

\mk
\mk

 {\bf Remerciements :} Nous souhaitons remercier chaleureusement Bérénice Delcroix-Oger, Clément Dupont, Hoel Queffelec, Giovanni Paolini et Luis Paris pour des discussions ayant aidé à la rédaction de ce texte. Nous remercions également Nicolas Bourbaki pour avoir contribué à améliorer ce texte.

\section{Groupes de Coxeter et groupes d'Artin}

Commençons par rappeler les définitions des groupes de Coxeter et d'Artin. Nous renvoyons le lecteur à \parencite{bourbaki_coxeter,humphreys,davis_coxeter,paris_kpi1,godelle_paris} pour plus de détails. 

\bdf[Matrice de Coxeter]
Considérons un ensemble fini $S$, une \emph{matrice de Coxeter} $(m_{st})_{s,t \in S}$ est une matrice symétrique telle que pour tout $s \in S$ nous ayons $m_{ss}=1$, et pour tous $s,t \in S$ distincts nous ayons $m_{st}=m_{ts} \in \{2,3,4,\dots,\infty\}$.
\edf

\bdf[Groupe de Coxeter]
Considérons un ensemble fini $S$, et une matrice de Coxeter $(m_{st})_{s,t \in S}$. Le \emph{groupe de Coxeter} $W$ associé à cette matrice a pour présentation
$$W = \<S \st \forall s,t \in S \mbox{ tels que $m_{st} \neq \infty$}, (st)^{m_{st}}=1\>.$$
\edf

\bdf[Groupe d'Artin]
Considérons un ensemble fini $S$ et une matrice de Coxeter $(m_{st})_{s,t \in S}$. Le \emph{groupe d'Artin} $G_W$ associé à cette matrice a pour présentation
$$G_W = \<S \st \forall s,t \in S \mbox{ tels que $m_{st} \neq \infty$}, [sts \cdots]_{m_{st}}=[tst \cdots]_{m_{st}}\>,$$
où $[sts\cdots]_m$ désigne le mot de longueur $m$ dont les lettres alternent entre $s$ et $t$.
\edf

Le cardinal de la partie génératrice $S$ est appelé le \emph{rang} du groupe de Coxeter ou du groupe d'Artin. La matrice de Coxeter est souvent représentée par un diagramme de Dynkin, qui est un graphe de sommets $S$, avec une arête entre les sommets $s$ et $t$ si $m_{s,t} \geq 3$, étiquetée par $m_{st}$ dès que $m_{st} \geq 4$. Lorsque le diagramme de Dynkin est connexe, c'est-à-dire lorsque le groupe de Coxeter ou d'Artin n'est pas un produit direct des sous-groupes engendrés par les composantes connexes, il est appelé \emph{irréductible}.

\bex
Fixons $n \geq 2$, et un ensemble $S=\{\sigma_1,\dots,\sigma_{n-1}\}$ de cardinal $n-1$. Considérons la matrice de Coxeter $(m_{st})_{s,t \in S}$ de type $A_{n-1}$, c'est-à-dire définie de la façon suivante : $m_{\sigma_i\sigma_{i+1}}=3$ et $m_{\sigma_i\sigma_j}=2$ dès que $|i-j| \geq 2$. Le diagramme de Dynkin associé est un chemin de longueur $n-1$ (voir la table~\ref{tab:classification}).

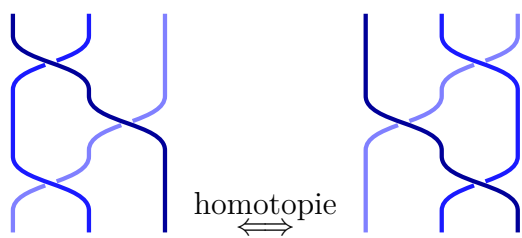
\begin{figure}
\begin{center}
\begin{tabular}{ccc}
\begin{tikzpicture}
\definecolor{col1}{rgb}{0,0,0.6}
\definecolor{col2}{rgb}{0,0,0.8}
\definecolor{col3}{rgb}{0.1,0.1,1}
\definecolor{col4}{rgb}{0.3,0.3,1}
\definecolor{col5}{rgb}{0.5,0.5,1}
\pic[
braid/.cd,
every strand/.style={ultra thick},
strand 1/.style={col1},
strand 2/.style={col3},
strand 3/.style={col5},
height=-0.8cm,]
{braid={s_1 s_2 s_1}};
\end{tikzpicture}
&
$\stackrel{\mbox{homotopie}}{\Longleftrightarrow}$
&
\begin{tikzpicture}
\definecolor{col1}{rgb}{0,0,0.6}
\definecolor{col2}{rgb}{0,0,0.8}
\definecolor{col3}{rgb}{0.1,0.1,1}
\definecolor{col4}{rgb}{0.3,0.3,1}
\definecolor{col5}{rgb}{0.5,0.5,1}
\pic[
braid/.cd,
every strand/.style={ultra thick},
strand 1/.style={col1},
strand 2/.style={col3},
strand 3/.style={col5},
height=-0.8cm,]
{braid={s_2 s_1 s_2}};
\end{tikzpicture}
\end{tabular}
\caption{La relation de tresse $\sigma_i \sigma_{i+1} \sigma_i = \sigma_{i+1} \sigma_i \sigma_{i+1}$.}
\label{fig:relation_tresse}
\end{center}
\end{figure}

Le groupe de Coxeter associé est le groupe symétrique $W=\frak{S}_n$, et le groupe d'Artin associé est le groupe de tresses $G_W=B_n$ à $n$ brins, avec la présentation d'Artin standard :
\begin{eqnarray*}
B_n &= \< \sigma_1,\dots,\sigma_{n-1} \st &\forall 1 \leq i \leq n-1, \sigma_i\sigma_{i+1}\sigma_i = \sigma_{i+1}\sigma_i\sigma_{i+1}, \\
&& \forall |i-j| \geq 2, \sigma_i\sigma_j = \sigma_j\sigma_i\>.
\end{eqnarray*}
Voir la figure~\ref{fig:relation_tresse} représentant la relation de tresse $\sigma_i\sigma_{i+1}\sigma_i = \sigma_{i+1}\sigma_i\sigma_{i+1}$.

\eex

\bex \label{ex:A2tilde}
Fixons $n \geq 2$, considérons un ensemble $S=\{\sigma_i, i \in \Z/n\Z\}$ de cardinal~$n$, et la matrice de Coxeter $(m_{st})_{s,t \in S}$ de type $\tilde{A}_{n-1}$, c'est-à-dire définie de la façon suivante : $m_{\sigma_i\sigma_{i+1}}=3$ et $m_{\sigma_i\sigma_j}=2$ dès que $|i-j| \geq 2$ modulo~$n$. Le diagramme de Dynkin associé est un $n$-cycle (voir la table~\ref{tab:classification}).

Le groupe de Coxeter associé est le groupe $W$ agissant sur l'espace euclidien $V=\{x \in \R^n \st x_1+x_2+\dots+x_n=0\}$, engendré par les réflexions orthogonales par rapport aux hyperplans d'équations $\{x_i-x_j=k\}$, pour $1 \leq i < j \leq n$ et $k \in \Z$. Le groupe d'Artin associé $G_W$ est parfois appelé groupe de tresses affine, et admet la présentation suivante :
\begin{eqnarray*}
G_W &= \< \sigma_i, i \in \Z/n\Z \st &\forall i \in \Z/n\Z, \sigma_i\sigma_{i+1}\sigma_i = \sigma_{i+1}\sigma_i\sigma_{i+1}, \\
&& \forall |i-j| \geq 2, \sigma_i\sigma_j = \sigma_j\sigma_i\>.
\end{eqnarray*}

\begin{figure}
\begin{center}
\begin{tikzpicture}
\def \p {0.05}
\def \op {0.1}
\def \gris {black}
\clip (-4.2,-4.2) rectangle (4.2,4.2);

\foreach \i in {-10,...,10}
\draw (\i,-10) -- (\i,10);

\foreach \i in {-10,...,10}
\draw (-17.32-\i/2,-10+\i*0.866) -- (17.32-\i/2,10+\i*0.866);

\foreach \i in {-10,...,10}
\draw (-17.32-\i/2,10-\i*0.866) -- (17.32-\i/2,-10-\i*0.866);

\draw[line width=2, color=blue] (0,-10) -- (0,10);
\draw[line width=2, color=blue] (-17.32,-10) -- (17.32,10);
\draw[line width=2, color=blue] (-17.32+1/2,10+1*0.866) -- (17.32+1/2,-10+1*0.866);

\node[line width=2, color=blue] at (0.2,-3) {$b$};
\node[line width=2, color=blue] at (3.6,-0.6) {$a$};
\node[line width=2, color=blue] at (3.6,1.8) {$c$};

\end{tikzpicture}
\caption{Le groupe de Coxeter de type $\tilde{A}_2$ est engendré par les réflexions orthogonales $a,b,c$.}
\label{fig:A2tilde_pavage}
\end{center}
\end{figure}
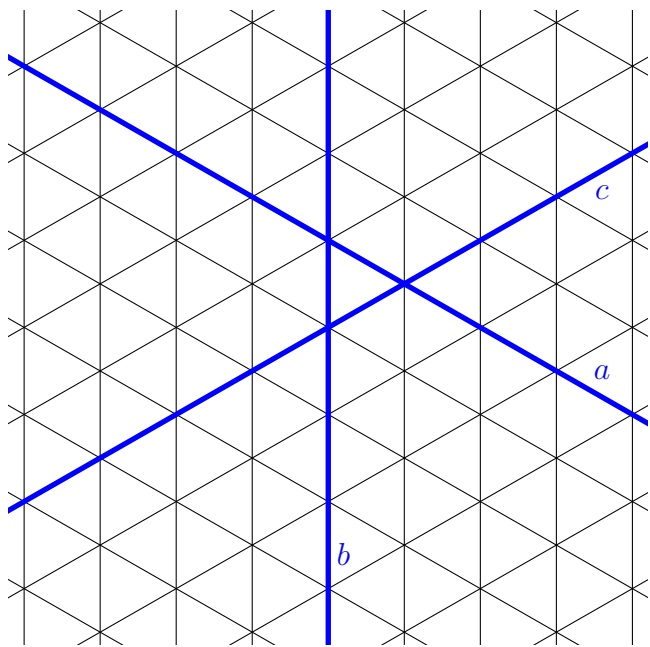

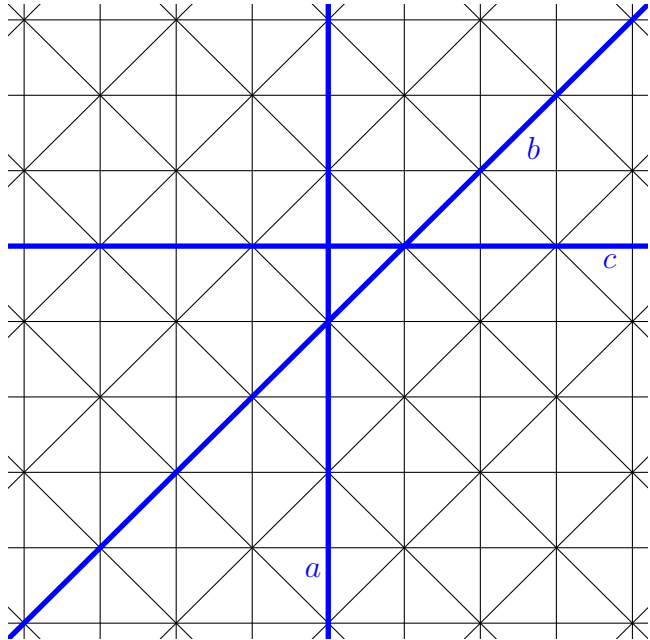
\begin{figure}
\begin{center}
\begin{tikzpicture}
\def \p {0.05}
\def \op {0.1}
\def \gris {black}
\clip (-4.2,-4.2) rectangle (4.2,4.2);

\foreach \i in {-10,...,10}
\draw (\i,-10) -- (\i,10);

\foreach \i in {-10,...,10}
\draw (-10,\i) -- (10,\i);

\foreach \i in {-10,...,10}
\draw (-10,-10+2*\i) -- (10,10+2*\i);

\foreach \i in {-10,...,10}
\draw (10,-10+2*\i) -- (-10,10+2*\i);

\draw[line width=2, color=blue] (0,-10) -- (0,10);
\draw[line width=2, color=blue] (-10,1) -- (10,1);
\draw[line width=2, color=blue] (-10,-10) -- (10,10);

\node[thick, color=blue] at (2.7,2.3) {$b$};
\node[thick, color=blue] at (-0.2,-3.3) {$a$};
\node[thick, color=blue] at (3.7,0.8) {$c$};

\end{tikzpicture}
\caption{Le groupe de Coxeter de type $\tilde{C}_2$ est engendré par les réflexions orthogonales $a,b,c$.}
\label{fig:C2tilde_pavage}
\end{center}
\end{figure}

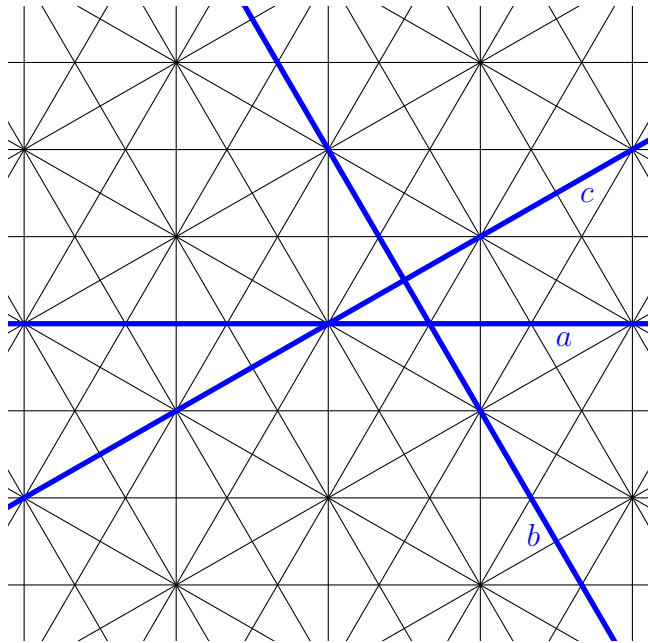
\begin{figure}
\begin{center}
\begin{tikzpicture}
\def \p {0.05}
\def \op {0.1}
\def \gris {black}
\clip (-4.2,-4.2) rectangle (4.2,4.2);

\foreach \i in {-10,...,10}
\draw (2*\i,-10) -- (2*\i,10);

\foreach \i in {-10,...,10}
\draw (-17.32-\i,-10+\i*2*0.866) -- (17.32-\i,10+\i*2*0.866);

\foreach \i in {-10,...,10}
\draw (-17.32-\i,10-\i*2*0.866) -- (17.32-\i,-10-\i*2*0.866);

\foreach \i in {-10,...,10}
\draw (-10+\i,-17.32-\i/1.732) -- (10+\i,17.32-\i/1.732);

\foreach \i in {-10,...,10}
\draw (10-\i,-17.32-\i/1.732) -- (-10-\i,17.32-\i/1.732);

\foreach \i in {-10,...,10}
\draw (-10,2*\i/1.732) -- (10,2*\i/1.732);

\draw[line width=2, color=blue] (-10,0) -- (10,0);
\draw[line width=2, color=blue] (-17.32,-10) -- (17.32,10);
\draw[line width=2, color=blue] (10+1,-17.32+1/1.732) -- (-10+1,17.32+1/1.732);

\node[thick, color=blue] at (2.7,-2.8) {$b$};
\node[thick, color=blue] at (3.1,-0.2) {$a$};
\node[thick, color=blue] at (3.4,1.7) {$c$};

\end{tikzpicture}
\caption{Le groupe de Coxeter de type $\tilde{G}_2$ est engendré par les réflexions orthogonales $a,b,c$.}
\label{fig:G2tilde_pavage}
\end{center}
\end{figure}

Lorsque $n=3$, nous noterons pour simplifier les générateurs standards $S=\{a,b,c\}$, et le pavage du plan euclidien $V$ est alors le 
pavage par des triangles équilatéraux, voir la figure~\ref{fig:A2tilde_pavage}. D'autres exemples de groupes de Coxeter affines de rang~$3$ sont ceux de type~$\tilde{C}_2$ et~$\tilde{G}_2$, voir les figures~\ref{fig:C2tilde_pavage} et \ref{fig:G2tilde_pavage}.
\eex

Tout groupe de Coxeter peut être réalisé comme groupe engendré par des réflexions linéaires dans un cône d'un espace vectoriel appelé cône de Tits \parencite{bourbaki_coxeter}. Par exemple, tout groupe de Coxeter fini de rang $n$ peut être réalisé comme groupe de réflexions orthogonales d'une sphère euclidienne de dimension $n-1$. Dans ce cas, le cône de Tits est l'espace vectoriel $\R^n$ entier. Ceci explique que lorsque le groupe de Coxeter est fini, le groupe d'Artin associé est dit \emph{(de type) sphérique}. Par exemple, c'est le cas des groupes de tresses.

Certains groupes de Coxeter peuvent être réalisés comme groupes de réflexions orthogonales affines d'un espace euclidien $\R^n$ : on les appelle alors groupes de réflexions affines, et le groupe d'Artin associé est dit \emph{(de type) affine}. Dans ce cas, le cône de Tits est un demi-espace de $\R^{n+1}$. Par exemple, le groupe de Coxeter de type $\tilde{A}_2$ présenté dans l'exemple~\ref{ex:A2tilde}, dans la représentation de Tits, est engendré par des réflexions linéaires, non orthogonales dans $\R^3$, mais qui stabilisent un hyperplan affine et agissent dessus comme réflexions affines orthogonales comme dans la figure~\ref{fig:A2tilde_pavage}.

\mk

Les diagrammes de Dynkin irréductibles de type sphérique ou affine sont classifiés \parencite{bourbaki_coxeter,davis_coxeter}. Nous avons rappelé la classification des diagrammes cristallographiques (ceux pour lesquels la matrice de Coxeter est à valeurs dans $\{1,2,3,4,6\}$) dans la table~\ref{tab:classification}.

\begin{table}
\begin{center}
\begin{tabular}{l|l|l|l}
& Type sphérique & & Type affine \\

\hline

$A_n$
&
\begin{tikzpicture}
\def \p {0.05}
\def \op {1}
\def \gris {black!10}

\draw[fill] (-3,0) circle (\p) node(s1) {};
\draw[fill] (-2,0) circle (\p) node(s2) {};
\draw[fill] (-1,0) circle (\p) node(s3) {};
\draw[fill] (1,0) circle (\p) node(s''n) {};
\draw[fill] (2,0) circle (\p) node(s'n) {};
\draw[fill] (3,0) circle (\p) node(sn) {};

\node (dots) at (0,0) {\bfseries $\dots$};

\draw [-] (s1) edge (s2) (s2) edge (s3) (s''n) edge (s'n) (s'n) edge (sn);

\end{tikzpicture}

& $\tilde{A_n}$
&
\begin{tikzpicture}
\def \p {0.05}
\def \op {1}
\def \gris {black!10}

\draw[fill] (-3,0) circle (\p) node(s1) {};
\draw[fill] (-2,0) circle (\p) node(s2) {};
\draw[fill] (-1,0) circle (\p) node(s3) {};
\draw[fill] (1,0) circle (\p) node(s''n) {};
\draw[fill] (2,0) circle (\p) node(s'n) {};
\draw[fill] (3,0) circle (\p) node(sn) {};
\draw[fill] (0,0.5) circle (\p) node(s0) {};

\node (dots) at (0,0) {\bfseries $\dots$};

\draw [-] (s1) edge (s2) (s2) edge (s3) (s''n) edge (s'n) (s'n) edge (sn) (sn) edge (s0) (s0) edge (s1);

\end{tikzpicture} \\

\hline

$B_n$
&
\begin{tikzpicture}
\def \p {0.05}
\def \op {1}
\def \gris {black!10}

\draw[fill] (-3,0) circle (\p) node(s1) {};
\draw[fill] (-2,0) circle (\p) node(s2) {};
\draw[fill] (-1,0) circle (\p) node(s3) {};
\draw[fill] (1,0) circle (\p) node(s''n) {};
\draw[fill] (2,0) circle (\p) node(s'n) {};
\draw[fill] (3,0) circle (\p) node(sn) {};

\node (dots) at (0,0) {\bfseries $\dots$};

\draw [-] (s1) edge (s2) (s2) edge (s3) (s''n) edge (s'n) (s'n) edge (sn);

\node (label) at (2.5,0.2) {$4$};

\end{tikzpicture}

& $\tilde{B_n}$
&
\begin{tikzpicture}
\def \p {0.05}
\def \op {1}
\def \gris {black!10}

\draw[fill] (-3,0.25) circle (\p) node(s1) {};
\draw[fill] (-3,-0.25) circle (\p) node(s'1) {};
\draw[fill] (-2,0) circle (\p) node(s2) {};
\draw[fill] (-1,0) circle (\p) node(s3) {};
\draw[fill] (1,0) circle (\p) node(s''n) {};
\draw[fill] (2,0) circle (\p) node(s'n) {};
\draw[fill] (3,0) circle (\p) node(sn) {};

\node (dots) at (0,0) {\bfseries $\dots$};

\draw [-] (s'1) edge (s2) (s1) edge (s2) (s2) edge (s3) (s''n) edge (s'n) (s'n) edge (sn);

\node (label) at (2.5,0.2) {$4$};

\end{tikzpicture}\\

\hline

$C_n$
&
\begin{tikzpicture}
\def \p {0.05}
\def \op {1}
\def \gris {black!10}

\draw[fill] (-3,0) circle (\p) node(s1) {};
\draw[fill] (-2,0) circle (\p) node(s2) {};
\draw[fill] (-1,0) circle (\p) node(s3) {};
\draw[fill] (1,0) circle (\p) node(s''n) {};
\draw[fill] (2,0) circle (\p) node(s'n) {};
\draw[fill] (3,0) circle (\p) node(sn) {};

\node (dots) at (0,0) {\bfseries $\dots$};

\draw [-] (s1) edge (s2) (s2) edge (s3) (s''n) edge (s'n) (s'n) edge (sn);

\node (label) at (2.5,0.2) {$4$};

\end{tikzpicture}

& $\tilde{C_n}$
&
\begin{tikzpicture}
\def \p {0.05}
\def \op {1}
\def \gris {black!10}

\draw[fill] (-3,0) circle (\p) node(s1) {};
\draw[fill] (-2,0) circle (\p) node(s2) {};
\draw[fill] (-1,0) circle (\p) node(s3) {};
\draw[fill] (1,0) circle (\p) node(s''n) {};
\draw[fill] (2,0) circle (\p) node(s'n) {};
\draw[fill] (3,0) circle (\p) node(sn) {};

\node (dots) at (0,0) {\bfseries $\dots$};

\draw [-] (s1) edge (s2) (s2) edge (s3) (s''n) edge (s'n) (s'n) edge (sn);

\node (label) at (-2.5,0.2) {$4$};
\node (label) at (2.5,0.2) {$4$};

\end{tikzpicture}\\

\hline
$D_n$
&
\begin{tikzpicture}
\def \p {0.05}
\def \op {1}
\def \gris {black!10}

\draw[fill] (-3,0) circle (\p) node(s1) {};
\draw[fill] (-2,0) circle (\p) node(s2) {};
\draw[fill] (-1,0) circle (\p) node(s3) {};
\draw[fill] (1,0) circle (\p) node(s''n) {};
\draw[fill] (2,0) circle (\p) node(s'n) {};
\draw[fill] (3,0.25) circle (\p) node(sn) {};
\draw[fill] (3,-0.25) circle (\p) node(sn') {};

\node (dots) at (0,0) {\bfseries $\dots$};

\draw [-] (s1) edge (s2) (s2) edge (s3) (s''n) edge (s'n) (s'n) edge (sn) (s'n) edge (sn');

\end{tikzpicture}

& $\tilde{D_n}$
&
\begin{tikzpicture}
\def \p {0.05}
\def \op {1}
\def \gris {black!10}

\draw[fill] (-3,0.25) circle (\p) node(s1) {};
\draw[fill] (-3,-0.25) circle (\p) node(s'1) {};
\draw[fill] (-2,0) circle (\p) node(s2) {};
\draw[fill] (-1,0) circle (\p) node(s3) {};
\draw[fill] (1,0) circle (\p) node(s''n) {};
\draw[fill] (2,0) circle (\p) node(s'n) {};
\draw[fill] (3,0.25) circle (\p) node(sn) {};
\draw[fill] (3,-0.25) circle (\p) node(sn') {};

\node (dots) at (0,0) {\bfseries $\dots$};

\draw [-] (s'1) edge (s2) (s1) edge (s2) (s2) edge (s3) (s''n) edge (s'n) (s'n) edge (sn) (s'n) edge (sn');

\end{tikzpicture}\\

\hline
$F_4$
&
\begin{tikzpicture}
\def \p {0.05}
\def \op {1}
\def \gris {black!10}

\draw[fill] (-3,0) circle (\p) node(s1) {};
\draw[fill] (-2,0) circle (\p) node(s2) {};
\draw[fill] (-1,0) circle (\p) node(s3) {};
\draw[fill] (0,0) circle (\p) node(s4) {};

\draw [-] (s1) edge (s2) (s2) edge (s3) (s3) edge (s4);

\node (label) at (-1.5,0.2) {$4$};

\end{tikzpicture}

& $\tilde{F_4}$
&
\begin{tikzpicture}
\def \p {0.05}
\def \op {1}
\def \gris {black!10}

\draw[fill] (-3,0) circle (\p) node(s1) {};
\draw[fill] (-2,0) circle (\p) node(s2) {};
\draw[fill] (-1,0) circle (\p) node(s3) {};
\draw[fill] (0,0) circle (\p) node(s4) {};
\draw[fill] (1,0) circle (\p) node(s5) {};

\draw [-] (s1) edge (s2) (s2) edge (s3) (s3) edge (s4) (s4) edge (s5);

\node (label) at (-1.5,0.2) {$4$};

\end{tikzpicture}\\

\hline

$G_2$
&
\begin{tikzpicture}
\def \p {0.05}
\def \op {1}
\def \gris {black!10}

\draw[fill] (-3,0) circle (\p) node(s1) {};
\draw[fill] (-2,0) circle (\p) node(s2) {};

\draw [-] (s1) edge (s2);

\node (label) at (-2.5,0.2) {$6$};

\end{tikzpicture}

& $\tilde{G_2}$
&
\begin{tikzpicture}
\def \p {0.05}
\def \op {1}
\def \gris {black!10}

\draw[fill] (-3,0) circle (\p) node(s1) {};
\draw[fill] (-2,0) circle (\p) node(s2) {};
\draw[fill] (-1,0) circle (\p) node(s3) {};

\draw [-] (s1) edge (s2) (s2) edge (s3);

\node (label) at (-2.5,0.2) {$6$};

\end{tikzpicture}\\

\hline

$E_6$
&
\begin{tikzpicture}
\def \p {0.05}
\def \op {1}
\def \gris {black!10}

\draw[fill] (-2,0) circle (\p) node(s1) {};
\draw[fill] (-1,0) circle (\p) node(s2) {};
\draw[fill] (0,0) circle (\p) node(s3) {};
\draw[fill] (1,0) circle (\p) node(s4) {};
\draw[fill] (2,0) circle (\p) node(s5) {};
\draw[fill] (0,-0.7) circle (\p) node(s'3) {};
\draw[opacity=0] (0,-1.4) circle (\p) node(s''3) {};

\draw [-] (s1) edge (s2) (s2) edge (s3) (s3) edge (s4) (s4) edge (s5) (s3) edge (s'3);

\end{tikzpicture}

& $\tilde{E_6}$
&
\begin{tikzpicture}
\def \p {0.05}
\def \op {1}
\def \gris {black!10}

\draw[fill] (-2,0) circle (\p) node(s1) {};
\draw[fill] (-1,0) circle (\p) node(s2) {};
\draw[fill] (0,0) circle (\p) node(s3) {};
\draw[fill] (1,0) circle (\p) node(s4) {};
\draw[fill] (2,0) circle (\p) node(s5) {};
\draw[fill] (0,-0.7) circle (\p) node(s'3) {};
\draw[fill] (0,-1.4) circle (\p) node(s''3) {};

\draw [-] (s1) edge (s2) (s2) edge (s3) (s3) edge (s4) (s4) edge (s5) (s3) edge (s'3) (s'3) edge (s''3);

\end{tikzpicture}\\

\hline

$E_7$
&
\begin{tikzpicture}
\def \p {0.05}
\def \op {1}
\def \gris {black!10}

\draw[fill] (-2,0) circle (\p) node(s1) {};
\draw[fill] (-1,0) circle (\p) node(s2) {};
\draw[fill] (0,0) circle (\p) node(s3) {};
\draw[fill] (1,0) circle (\p) node(s4) {};
\draw[fill] (2,0) circle (\p) node(s5) {};
\draw[fill] (3,0) circle (\p) node(s6) {};
\draw[fill] (0,-0.7) circle (\p) node(s'3) {};

\draw [-] (s1) edge (s2) (s2) edge (s3) (s3) edge (s4) (s4) edge (s5) (s5) edge (s6)  (s3) edge (s'3);

\end{tikzpicture}

& $\tilde{E_7}$
&
\begin{tikzpicture}
\def \p {0.05}
\def \op {1}
\def \gris {black!10}

\draw[fill] (-3,0) circle (\p) node(s0) {};
\draw[fill] (-2,0) circle (\p) node(s1) {};
\draw[fill] (-1,0) circle (\p) node(s2) {};
\draw[fill] (0,0) circle (\p) node(s3) {};
\draw[fill] (1,0) circle (\p) node(s4) {};
\draw[fill] (2,0) circle (\p) node(s5) {};
\draw[fill] (3,0) circle (\p) node(s6) {};
\draw[fill] (0,-0.7) circle (\p) node(s'3) {};

\draw [-] (s0) edge (s1) (s1) edge (s2) (s2) edge (s3) (s3) edge (s4) (s4) edge (s5) (s5) edge (s6) (s3) edge (s'3);

\end{tikzpicture}\\

\hline

$E_8$
&
\begin{tikzpicture}
\def \p {0.05}
\def \op {1}
\def \gris {black!10}

\draw[fill] (-2,0) circle (\p) node(s1) {};
\draw[fill] (-1,0) circle (\p) node(s2) {};
\draw[fill] (0,0) circle (\p) node(s3) {};
\draw[fill] (1,0) circle (\p) node(s4) {};
\draw[fill] (2,0) circle (\p) node(s5) {};
\draw[fill] (3,0) circle (\p) node(s6) {};
\draw[fill] (4,0) circle (\p) node(s7) {};
\draw[fill] (0,-0.7) circle (\p) node(s'3) {};f

\draw [-] (s1) edge (s2) (s2) edge (s3) (s3) edge (s4) (s4) edge (s5) (s5) edge (s6) (s6) edge (s7) (s3) edge (s'3);

\end{tikzpicture}

& $\tilde{E_8}$
&
\begin{tikzpicture}
\def \p {0.05}
\def \op {1}
\def \gris {black!10}

\draw[fill] (-2,0) circle (\p) node(s1) {};
\draw[fill] (-1,0) circle (\p) node(s2) {};
\draw[fill] (0,0) circle (\p) node(s3) {};
\draw[fill] (1,0) circle (\p) node(s4) {};
\draw[fill] (2,0) circle (\p) node(s5) {};
\draw[fill] (3,0) circle (\p) node(s6) {};
\draw[fill] (4,0) circle (\p) node(s7) {};
\draw[fill] (5,0) circle (\p) node(s8) {};
\draw[fill] (0,-0.7) circle (\p) node(s'3) {};

\draw [-] (s1) edge (s2) (s2) edge (s3) (s3) edge (s4) (s4) edge (s5) (s5) edge (s6) (s6) edge (s7) (s7) edge (s8) (s3) edge (s'3);

\end{tikzpicture}\\


\end{tabular}
\mk
\caption{Diagrammes de Dynkin cristallographiques sphériques et affines.}
\label{tab:classification}
\end{center}
\end{table}

Si $T \subset S$, on notera $W_T$ le sous-groupe du groupe de Coxeter $W$ engendré par $T$, appelé \emph{sous-groupe parabolique standard}. C'est aussi un groupe de Coxeter, dont la matrice de Coxeter est donnée par la restriction à $T$ de la matrice d'origine.

\mk

Considérons un groupe de Coxeter affine $W$ de rang $n$, agissant comme groupe de réflexions orthogonales par isométries affines sur l'espace euclidien $\R^n$. On appelle \emph{réflexion} de $W$ tout élément conjugué à un élément de $S$. L'ensemble des réflexions est noté $R$. Considérons l'ensemble ${\mathcal A}$ des hyperplans fixés par les réflexions de $R$, voir par exemple les figures~\ref{fig:A2tilde_pavage}, \ref{fig:C2tilde_pavage} et \ref{fig:G2tilde_pavage} représentant les arrangements d'hyperplans en types $\tilde{A}_2$, $\tilde{C}_2$ et $\tilde{G_2}$. Alors l'espace de configuration est le quotient
$$Y_W = \left(\C^n \smallsetminus \bigcup_{H \in {\mathcal A}}  H \underset{\R}{\otimes} \C\right) / W,$$
et son groupe fondamental est le groupe d'Artin de type affine $G_W$ \parencite{vanderlek}. La conjecture du $K(\pi,1)$ affirme dans ce cas que $Y_W$ est un espace classifiant pour le groupe d'Artin $G_W$, c'est-à-dire que le revêtement universel de $Y_W$ est contractile.

\section{Intervalles et groupes de Garside}

\label{sec:interval_garside_groups}

Nous allons maintenant donner une présentation rapide des groupes de Garside, et plus généralement des groupes construits à partir d'intervalles. Nous suivrons le point de vue de McCammond et Sulway \parencite{mccammond_dual_euclidian,mccammond_sulway}, et nous renvoyons le lecteur à \parencite{garside,dehornoy_paris,dehornoy_ENS,dehornoy_garside,digne_An,digne_Cn} pour plus de détails sur les structures de Garside.

\bdf[Groupe d'intervalle] \label{def:groupe_intervalle}
Soit $G$ un groupe, et $R \subset G$ un sous-ensemble (éventuellement infini) engendrant $G$ tel que $R=R^{-1}$. Supposons qu'à chaque élément $r \in R$ est associé un \emph{poids} $l(r)>0$, et supposons de plus que $l(R)$ soit un sous-ensemble discret de $\R$. Supposons également que, pour tout $r \in R$, nous ayons $l(r^{-1})=l(r)$.

Considérons le graphe de Cayley $\Cay(G,R)$ de $G$ par rapport à $R$, dont les arêtes orientées sont étiquetées par les éléments de $R$. Considérons $\Cay(G,R)$ comme un espace métrique, en déclarant qu'une arête étiquetée $r$ est de longueur $l(r)$. Pour tous $g,h \in G$, notons $d_R(g,h)$ la longueur du plus court chemin de $g$ à $h$ dans $\Cay(G,R)$.

Fixons $g \in G$, et considérons l'intervalle
$$[1,g]^G = \{h \in G \st d_R(1,h) + d_R(h,g)=d_R(1,g)\},$$
c'est aussi la réunion des géodésiques de $1$ à $g$ dans $\Cay(G,R)$. On peut également le voir comme un sous-graphe de $\Cay(G,R)$. Remarquons que $[1,g]^G$ est naturellement ordonné, en déclarant que $h \leq k$ si $d_R(1,h)+d_R(h,k)+d_R(k,g)=d_R(1,g)$.

Soit $R_g \subset R$ l'ensemble des étiquettes apparaissant parmi les géodésiques de $1$ à $g$. Le \emph{groupe d'intervalle} $G_g$ est le groupe engendré par $R_g$, avec les relations fournies par les mots apparaissant en lisant les étiquettes des cycles dans $[1,g]^G$. 
\edf

Remarquons qu'on peut choisir un poids constant égal à $1$. Cependant, pour certains groupes de Coxeter affines, nous verrons qu'il est important d'autoriser d'autres poids. Rappelons maintenant ce qu'est un treillis.

\bdf[Treillis]
Un ensemble ordonné $P$ est appelé \emph{treillis} si toute paire d'éléments $p,q \in P$ a une borne inférieure et une borne supérieure.
\edf

Nous pouvons maintenant définir les groupes de Garside.

\bdf[Groupe de Garside] \label{def:garside}
Un groupe $G$ est appelé \emph{groupe de Garside} s'il existe $R \subset G$ et
$l\colon R \ra \R_+^*$ comme ci-dessus, et un élément $\delta \in G$, tels qu'on ait les propriétés suivantes.
\bit
\item Le groupe $G$ est engendré par $R_\delta$.
\item L'ensemble ordonné $[1,\delta]^G$ est un treillis.
\item L'intervalle $[1,\delta]^G$ est \emph{équilibré}, c'est-à-dire que pour tout $h \in G$, nous avons que $h \in [1,\delta]^G$ si et seulement si $\delta h^{-1} \in [1,\delta]^G$.
\eit
\edf

Remarquons que, selon les auteurs, la notion de groupe de Garside peut éventuellement demander que l'ensemble $R_\delta$ soit fini. Dans ce texte, nous autoriserons cet ensemble à être infini, ce qui s'avère nécessaire pour les groupes d'Artin de type affine.

\bex
L'exemple le plus simple de groupe de Garside est, pour $n \geq 1$, le groupe abélien libre $\Z^n$, avec pour élément de Garside $\delta=(1,1,\dots,1)$. Nous verrons ci-dessous que les groupes de tresses donnent des exemples plus intéressants.
\eex

Remarquons que si $R$ est stable par conjugaison dans $G$, et que le poids de deux éléments conjugués de $R$ est identique, alors la dernière condition de la définition~\ref{def:garside} est toujours satisfaite. La condition la plus importante est alors que l'intervalle $[1,\delta]^G$ soit un treillis.

Tout groupe d'intervalle peut être réalisé comme groupe fondamental d'un complexe d'intervalle, que nous définissons maintenant.

\bdf[Complexe d'intervalle]
Soit $G$ un groupe, $R \subset G$ une partie genératrice et un poids $l$ comme dans la définition~\ref{def:groupe_intervalle}. Soit $g \in G$ avec un intervalle $[1,g]^G$ équilibré. Considérons comme modèle du simplexe $\Delta^d$ de dimension $d$ l'ensemble
$$\Delta^d = \{a \in \R^d \st 1 \geq a_1 \geq a_2 \geq \cdots \geq a_d \geq 0\},$$
aussi appelé orthosimplexe de dimension $d$.

Nous allons définir un $\Delta$-complexe (au sens de 
\textcite{hatcher}), appelé \emph{complexe d'intervalle} de $[1,g]^G$, qui est un recollement de simplexes par identifications de certaines faces. Il y a un $d$-simplexe noté $[x_1|x_2| \dots |x_d]$ pour chaque chaîne dans l'ensemble ordonné $[1,g]^G \smallsetminus \{1\}$, c'est-à-dire $x_1,x_2,\dots,x_d \in [1,g]^G \smallsetminus \{1\}$ tels que :
\bit
\item $x_1x_2 \cdots x_d \in [1,g]^G$ et
\item $l(x_1x_2 \cdots x_d)=l(x_1)+l(x_2)+ \cdots +l(x_d)$.
\eit
Décrivons les recollements des $d+1$ facettes (faces de codimension $1$) du $d$-simplexe $[x_1|x_2| \dots |x_d]$ :
\bit
\item La facette $\{1 = a_1 \geq a_2 \geq \cdots \geq a_d \geq 0\}$ du simplexe $[x_1|x_2| \dots |x_d]$ est identifiée avec le $(d-1)$-simplexe $[x_2| \dots |x_d]$ par $(1,a_2,\dots,a_d) \mapsto (a_2,\dots,a_d)$.
\item Pour $1 \leq i \leq d-1$, la facette $\{1 \geq a_1 \geq  \cdots \geq a_i = a_{i+1} \geq \cdots \geq a_d \geq 0\}$ du simplexe $[x_1|x_2| \dots |x_d]$ est identifiée avec le $(d-1)$-simplexe $[x_1| \dots |x_ix_{i+1}| \dots |x_d]$ par $(a_1,\dots,a_i,a_i,a_{i+2},\dots,a_d) \mapsto (a_1,\dots,a_i,a_{i+2},\dots,a_d)$.
\item La facette $\{1 \geq a_1 \geq a_2 \geq \cdots \geq a_d = 0\}$ du simplexe $[x_1|x_2| \dots |x_d]$ est identifiée avec le $(d-1)$-simplexe $[x_1| \dots |x_{d-1}]$ par $(a_1,\dots,a_{d-1},0) \mapsto (a_1,\dots,a_{d-1})$.
\eit
\edf

\begin{figure}
\begin{center}
\begin{tikzpicture}
\def \p {0.05}
\def \op {1}
\def \gris {black!10}
\def \lgth {0.5cm}
\def \wdth {0.2cm}

\draw (0,0) node(e) {};
\draw (3,0) node(a) {};
\draw (3,3) node(ab) {};

\draw[black,fill opacity=\op,fill=\gris] (e.center) -- (a.center) -- (ab.center) -- cycle;

\draw[-{Stealth[length=\lgth,width=\wdth]}] (e.center) -- (a.center);
\draw[-{Stealth[length=\lgth,width=\wdth]}] (e.center) -- (ab.center);
\draw[-{Stealth[length=\lgth,width=\wdth]}] (a.center) -- (ab.center);

\node at (1.5,-0.4) {$[x_1]$};
\node at (3.4,1.5) {$[x_2]$};
\node at (1.1,1.9) {$[x_1x_2]$};
\node at (1.8,0.7) {$[x_1|x_2]$};
\node at (-0.2,-0.2) {$[\,]$};
\node at (3.2,-0.2) {$[\,]$};
\node at (3.2,3.2) {$[\,]$};

\end{tikzpicture}
\caption{Le $2$-simplexe $[x_1|x_2]$, dont les $3$ sommets sont identifiés.}
\label{fig:2_simplex}
\end{center}
\end{figure}
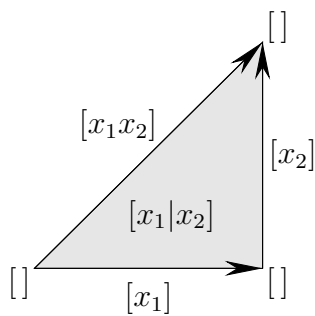

Voir la figure~\ref{fig:2_simplex} du $2$-simplexe $[x_1|x_2]$. Remarquons que ce complexe d'intervalle a un unique sommet noté $[\,]$, et que ses arêtes sont en bijection avec les éléments de $[1,g]^G \smallsetminus \{1\}$.

\bpro
Le groupe fondamental du complexe d'intervalle associé à $[1,g]^G$ est le groupe d'intervalle $G_g$.
\epro

\bp
Le groupe fondamental du complexe d'intervalle $K$ associé à $[1,g]^G$ admet une présentation déterminée par son $2$-squelette :
$$\pi_1(K,[\,]) = \<[1,g]^G \smallsetminus \{1\} \st (x_1)(x_2)=(x_1x_2) \mbox{ pour tout $2$-simplexe } [x_1|x_2]\>.$$
Le groupe $G_g$ a le même ensemble de générateurs $[1,g]^G \smallsetminus \{1\}$, et a un ensemble de relations plus grand, donné par l'ensemble des cycles dans $[1,g]^G$. Considérons un cycle $c=(1,y_1,y_2,\dots,y_k=1)$ basé en $1$ dans le graphe $[1,g]^G$. On peut le réécrire comme un produit de cycles de la forme $c_i=(1,y_i,y_{i+1},1)$, pour $1 \leq i \leq k-1$. Fixons $1 \leq i \leq k-1$, et supposons par exemple que $y_i < y_{i+1}$. Alors l'image du cycle $c_i$ dans le complexe $K$ borde le $2$-simplexe $[y_i|y_i^{-1}y_{i+1}]$, et est donc homotope à zéro. Ainsi chaque cycle $c$ dans $[1,g]^G$ a une image dans $K$ homotope à zéro, donc $G_g$ est bien le groupe fondamental de $K$.
\ep

Ce point de vue est une méthode efficace pour construire des groupes de Garside et leurs espaces classifiants :

\bthm[{\cite[Theorem~3.1]{charney_meier_whittlesey}; \cite[Theorem~0.1]{dehornoy_lafont}}] \label{thm:treillis_classifiant}
Si $[1,g]^G$ est un treillis équilibré, alors $G_g$ est un groupe de Garside, et le complexe d'intervalle associé à $[1,g]^G$ est un espace classifiant pour $G_g$.
\ethm

Remarquons qu'une autre preuve de ce résultat, à l'aide d'une métrique à courbure négative, est proposée dans \textcite[Theorem~E]{haettel_kpi1}.

\bex[Groupes de tresses] \label{ex:tresses_garside}
Considérons le groupe de tresses $B_n$ à $n$ brins, le groupe de Coxeter associé est le groupe symétrique $\mathfrak{S}_n$. Il peut être muni de deux structures de Garside différentes, la structure de Garside standard et la structure de Garside duale \parencite{birman_ko_lee,bessis}.

\mk

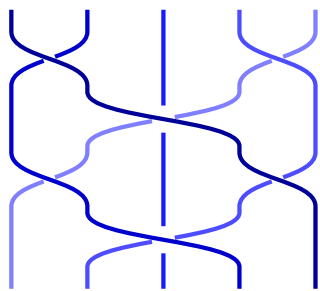
\begin{figure}
\begin{center}
\begin{tikzpicture}
\definecolor{col1}{rgb}{0,0,0.6}
\definecolor{col2}{rgb}{0,0,0.8}
\definecolor{col3}{rgb}{0.1,0.1,1}
\definecolor{col4}{rgb}{0.3,0.3,1}
\definecolor{col5}{rgb}{0.5,0.5,1}
\pic[
braid/.cd,
every strand/.style={ultra thick},
strand 1/.style={col1},
strand 2/.style={col2},
strand 3/.style={col3},
strand 4/.style={col4},
strand 5/.style={col5},
height=-0.8cm,]
{braid={|s_1-s_4| s_{2,4} |s_1-s_4| s_{2,4}}};
\end{tikzpicture}
\caption{Une représentation de l'élément de Garside standard $\Delta$ de $B_5$.}
\label{fig:garside_standard}
\end{center}
\end{figure}

Dans la structure de Garside standard, on considère l'élément de longueur maximale~$\Delta$ de~$\mathfrak{S}_n$ par rapport à la partie génératrice standard~$S$ constituée des transpositions $(i,i+1)$. L'élément $\Delta$ représente alors le \og demi-tour\fg{} qui à $i \in \{1,\dots,n\}$ associe $n+1-i$, il est de longueur $\f{n(n-1)}{2}$ (voir la figure~\ref{fig:garside_standard}). L'intervalle $[1,\Delta]^{\mathfrak{S}_n}$ est ainsi l'ensemble des écritures minimales de $\Delta$ comme produit de $\f{n(n-1)}{2}$ transpositions adjacentes.

\mk

Remarquons de plus que les ensembles $S$ et $[1,\Delta]^{\mathfrak{S}_n}$ se relèvent naturellement comme ensembles dans le groupe de tresses $B_n$, en associant à la transposition $(i,i+1)$ le générateur standard $\sigma_i$ du groupe de tresses.

\mk

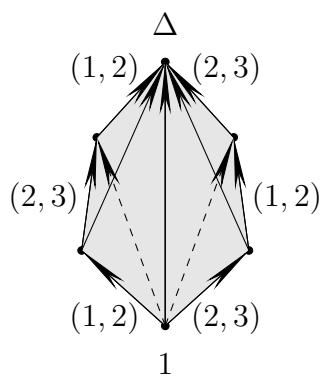
\begin{figure}
\begin{center}
\begin{tikzpicture}
\def \p {0.05}
\def \op {1}
\def \gris {black!10}
\def \lgth {0.7cm}
\def \wdth {0.15cm}

\draw[fill] (0,0) circle (\p) node(e) {};
\draw[fill] (-1.1,1) circle (\p) node(a) {};
\draw[fill] (-0.9,2.5) circle (\p) node(ab) {};
\draw[fill] (0,3.5) circle (\p) node(aba) {};
\draw[fill] (1.1,1) circle (\p) node(b) {};
\draw[fill] (0.9,2.5) circle (\p) node(ba) {};

\draw[black,fill opacity=\op,fill=\gris] (e.center) -- (a.center) -- (ab.center) -- (aba.center) -- cycle;
\draw[black,fill opacity=\op,fill=\gris] (e.center) -- (b.center) -- (ba.center) -- (aba.center) -- cycle;
\draw[-{Stealth[length=\lgth,width=\wdth]}] (e.center) -- (a.center);
\draw[-{Stealth[length=\lgth,width=\wdth]}] (e.center) -- (b.center);
\draw[-{Stealth[length=\lgth,width=\wdth]}] (e.center) -- (aba.center);
\draw[-{Stealth[length=\lgth,width=\wdth]}] (a.center) -- (aba.center);
\draw[-{Stealth[length=\lgth,width=\wdth]}] (b.center) -- (aba.center);
\draw[-{Stealth[length=\lgth,width=\wdth]}] (ab.center) -- (aba.center);
\draw[-{Stealth[length=\lgth,width=\wdth]}] (ba.center) -- (aba.center);
\draw[-{Stealth[length=\lgth,width=\wdth]}] (b.center) -- (ba.center);
\draw[-{Stealth[length=\lgth,width=\wdth]}] (a.center) -- (ab.center);
\draw[dashed,-{Stealth[length=\lgth,width=\wdth]}] (e.center) -- (ab.center);
\draw[dashed,-{Stealth[length=\lgth,width=\wdth]}] (e.center) -- (ba.center);

\node at (-0.8,0.1) {$(1,2)$};
\node at (-1.6,1.7) {$(2,3)$};
\node at (-0.8,3.4) {$(1,2)$};
\node at (0.8,0.1) {$(2,3)$};
\node at (1.6,1.7) {$(1,2)$};
\node at (0.8,3.4) {$(2,3)$};

\node at (0,-0.5) {$1$};
\node at (0,4) {$\Delta$};

\end{tikzpicture}
\caption{Le complexe d'intervalle de $\mathfrak{S}_3$ pour la structure de Garside standard, obtenu par identification des simplexes ayant le même étiquetage.}
\label{fig:complexe_b3_garside_standard}
\end{center}
\end{figure}

Dans le cas $n=3$, cet intervalle $[1,\Delta]^{\mathfrak{S}_3}$ est constitué de deux chaînes maximales de longueur $3$, car $\Delta = (1,2) (2,3) (1,2) = (2,3) (1,2) (2,3)$. Voir la figure~\ref{fig:complexe_b3_garside_standard} pour le complexe d'intervalle associé. Son revêtement universel est homéomorphe au produit d'une droite avec un complexe triangulaire ressemblant à un arbre $3$-régulier.

\mk

\begin{figure}
\begin{center}
\begin{tikzpicture}
\definecolor{col1}{rgb}{0,0,0.6}
\definecolor{col2}{rgb}{0,0,0.8}
\definecolor{col3}{rgb}{0.1,0.1,1}
\definecolor{col4}{rgb}{0.3,0.3,1}
\definecolor{col5}{rgb}{0.5,0.5,1}
\pic[
braid/.cd,
every strand/.style={ultra thick},
strand 1/.style={col1},
strand 2/.style={col2},
strand 3/.style={col3},
strand 4/.style={col4},
strand 5/.style={col5},
height=-0.8cm,]
{braid={s_1 s_2 s_3 s_4}};
\end{tikzpicture}
\caption{Une représentation de l'élément de Garside dual $\delta$ de $B_5$.}
\label{fig:garside_dual}
\end{center}
\end{figure}
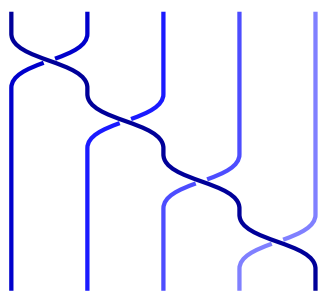

Dans la structure de Garside duale, on considère l'ensemble des conjugués de~$S$, c'est-à-dire l'ensemble~$T$ de toutes les transpositions de $\mathfrak{S}_n$. On considère un élément~$\delta$ obtenu comme produit des éléments de~$S$ dans un certain ordre, par exemple $\delta = (1,2)(2,3) \cdots (n-1,n)$ qui représente le $n$-cycle $(1,2,\dots,n)$ (voir la figure~\ref{fig:garside_dual}). Toutes les écritures minimales de~$\delta$ comme produit d'élements de~$T$ ont $n-1$ éléments, et l'ensemble de ces écritures forme l'intervalle $[1,\delta]^{\mathfrak{S}_n}$. Remarquons que $|\Delta|_S = \f{n(n-1)}{2} = |T|$, tandis que $|\delta|_T = n-1 = |S|$. Ceci est l'une des justifications pour la terminologie \og duale\fg{}.

\mk

L'intervalle $[1,\delta]^{\mathfrak{S}_n}$ possède également une belle interprétation comme ensemble de partitions non croisées. Considérons l'ensemble $U_n$ des $n$ sommets d'un $n$-gone régulier dans le plan, alors une partition de $U_n$ est dite \emph{non croisée} si les enveloppes convexes des éléments de la partition ne s'intersectent pas. Il y a un ordre naturel sur l'ensemble des partitions non croisées de $n$ points, donné par le raffinement de partitions. L'intervalle $[1,\delta]^{\mathfrak{S}_n}$ est alors isomorphe au treillis des partitions non croisées de $n$ points. Pour plus de détails sur les partitions non croisées, nous renvoyons le lecteur à \parencite{brady_partial_order,athanasiadis_brady_watt,brady_watt_NCP,brady_mccammond_orthoschemes,armstrong}.

\mk

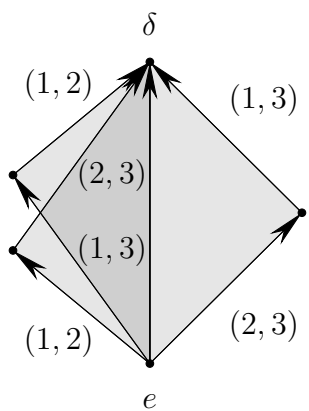
\begin{figure}
\begin{center}
\begin{tikzpicture}
\def \p {0.05}
\def \op {0.1}
\def \gris {black}

\def \lgth {0.5cm}
\def \wdth {0.2cm}

\draw[fill] (0,0) circle (\p) node(e) {};
\draw[fill] (-1.8,1.5) circle (\p) node(a) {};
\draw[fill] (2,2) circle (\p) node(b) {};
\draw[fill] (-1.8,2.5) circle (\p) node(c) {};
\draw[fill] (0,4) circle (\p) node(delta) {};

\draw[black,fill opacity=\op,fill=\gris] (e.center) -- (a.center) -- (delta.center) -- cycle;
\draw[black,fill opacity=\op,fill=\gris] (e.center) -- (c.center) -- (delta.center) -- cycle;
\draw[black,fill opacity=\op,fill=\gris] (e.center) -- (b.center) -- (delta.center) -- cycle;

\draw[-{Stealth[length=\lgth,width=\wdth]}] (e.center) -- (a.center);
\draw[-{Stealth[length=\lgth,width=\wdth]}] (e.center) -- (b.center);
\draw[-{Stealth[length=\lgth,width=\wdth]}] (e.center) -- (c.center);
\draw[-{Stealth[length=\lgth,width=\wdth]}] (e.center) -- (delta.center);
\draw[-{Stealth[length=\lgth,width=\wdth]}] (a.center) -- (delta.center);
\draw[-{Stealth[length=\lgth,width=\wdth]}] (b.center) -- (delta.center);
\draw[-{Stealth[length=\lgth,width=\wdth]}] (c.center) -- (delta.center);

\node at (-1.2,0.3) {$(1,2)$};
\node at (-1.2,3.7) {$(1,2)$};
\node at (1.5,0.5) {$(2,3)$};
\node at (1.5,3.5) {$(1,3)$};
\node at (-0.5,1.5) {$(1,3)$};
\node at (-0.5,2.5) {$(2,3)$};

\node at (0,-0.5) {$e$};
\node at (0,4.5) {$\delta$};

\end{tikzpicture}
\caption{Le complexe d'intervalle de $\mathfrak{S}_3$ pour la structure de Garside duale, obtenu par identification des simplexes ayant le même étiquetage.}
\label{fig:complexe_b3_garside_dual}
\end{center}
\end{figure}

Dans le cas $n=3$, cet intervalle $[1,\delta]^{\mathfrak{S}_3}$ est constitué de trois chaînes maximales de longueur $2$, car $\delta = (1,2) (2,3) = (2,3) (1,3) = (1,3) (1,2)$. Voir la figure~\ref{fig:complexe_b3_garside_dual} pour le complexe d'intervalle associé. Son revêtement universel est homéomorphe au produit d'une droite avec un arbre $3$-régulier.

\eex

\section{Groupes d'Artin duaux}

\label{sec:dual_artin_groups}

Soit $W$ un groupe de Coxeter, soit $R$ l'ensemble de ses réflexions, et soit
$S \subset R$ un ensemble simple de réflexions. Choisissons un poids $l\colon R \ra \R_+^*$ constant égal à $1$. Un \emph{élément de Coxeter} $w$ est le produit des éléments de $S$ dans un ordre quelconque. On appelle \emph{groupe d'Artin dual} associé à $w$ le groupe d'intervalle $W_w$ défini à partir de $R$.

\mk

\bex[Interprétation géométrique]
Soit $W$ un groupe de Coxeter irréductible fini ou affine, agissant comme groupe de réflexions sur la sphère $M=\SS^n$ ou l'espace euclidien $M=\R^n$. L'ensemble des hyperplans des réflexions de $R$ découpe $M$ en composantes connexes, appelées chambres (ouvertes) du complexe de Coxeter. Il existe une chambre $C_0$ de $M$ dont les réflexions par rapport aux faces sont précisément les éléments de $S$. Ainsi un élément de Coxeter $w$ peut être interprété comme le produit des réflexions par rapport aux faces de la chambre $C_0$, dans un ordre quelconque. 
\eex

\mk

On dit que deux éléments de Coxeter $w,w'$ sont \emph{géométriquement
  équivalents} s'il existe un automorphisme $\phi\colon W \ra W$ préservant l'ensemble $R$ des réflexions, envoyant $S$ sur un ensemble simple de réflexions, et tel que $\phi(w)=w'$. Pour la plupart des groupes de Coxeter sphériques ou affines, tous les éléments de Coxeter sont géométriquement équivalents.

\bthm[{\cite[Corollary~7.6, Corollary~7.8]{mccammond_dual_euclidian}}]
Soit $W$ un groupe de Coxeter irréductible de type sphérique ou affine, mais pas de type $\tilde{A}$. Alors tous les éléments de Coxeter sont géométriquement équivalents.

Soit $W$ un groupe de Coxeter affine de type $\tilde{A}_n$. Alors tout élément de Coxeter est géométriquement équivalent à un $(p,q)$-bigone de Coxeter (voir~\cite[Definition~7.7,Example~11.6]{mccammond_dual_euclidian}), où $p \geq q$ et $p+q=n+1$.
\ethm

\bex
Considérons le groupe de Coxeter $W$ affine de type $\tilde{A}_2$, engendré par les trois réflexions $a,b,c$ du plan euclidien bordant la chambre $C_0$ comme sur la figure~\ref{fig:A2tilde_element_coxeter}. Alors $w=abc$ est un exemple d'élément de Coxeter, qui consiste en la translation-réflexion le long de l'axe de Coxeter $\ell$ représenté sur la figure~\ref{fig:A2tilde_element_coxeter}. Remarquons que dans ce cas particulier, tous les éléments de Coxeter sont géométriquement équivalents.
Le cas similaire du type $\tilde{G_2}$ est représenté sur la figure~\ref{fig:G2tilde_element_coxeter}.
\eex

\begin{figure}
\begin{center}
\begin{tikzpicture}
\def \p {0.05}
\def \op {0.5}
\def \gris {black!50}
\clip (-4.2,-4.2) rectangle (4.2,4.2);

\foreach \i in {-10,...,10}
\draw (\i,-10) -- (\i,10);

\foreach \i in {-10,...,10}
\draw (-17.32-\i/2,-10+\i*0.866) -- (17.32-\i/2,10+\i*0.866);

\foreach \i in {-10,...,10}
\draw (-17.32-\i/2,10-\i*0.866) -- (17.32-\i/2,-10-\i*0.866);

\draw[line width=2, color=blue] (0,-10) -- (0,10);
\draw[line width=2, color=blue] (-17.32,-10) -- (17.32,10);
\draw[line width=2, color=blue] (-17.32+1/2,10+1*0.866) -- (17.32+1/2,-10+1*0.866);

\node[thick, color=blue] at (0.2,-3) {$b$};
\node[thick, color=blue] at (3.6,-0.6) {$a$};
\node[thick, color=blue] at (3.6,1.8) {$c$};

\draw[line width=2, color=red] (0.5,-10) -- (0.5,10);
\node[thick, color=red] at (0.7,3.4) {$\ell$};
\draw[line width=2, color=red,->] (0.65,0.5) -- (0.65,0.5+3/1.732);
\node[thick, color=red] at (0.85,1.1) {$w$};

\draw[black,fill opacity=\op,fill=\gris] (0,0) -- (1,1/1.732) -- (0,2/1.732) -- cycle;
\node[thick, color=\gris] at (-0.3,0.6) {$C_0$};

\end{tikzpicture}
\caption{Un élément de Coxeter pour le groupe de Coxeter affine de type $\tilde{A}_2$.}
\label{fig:A2tilde_element_coxeter}
\end{center}
\end{figure}
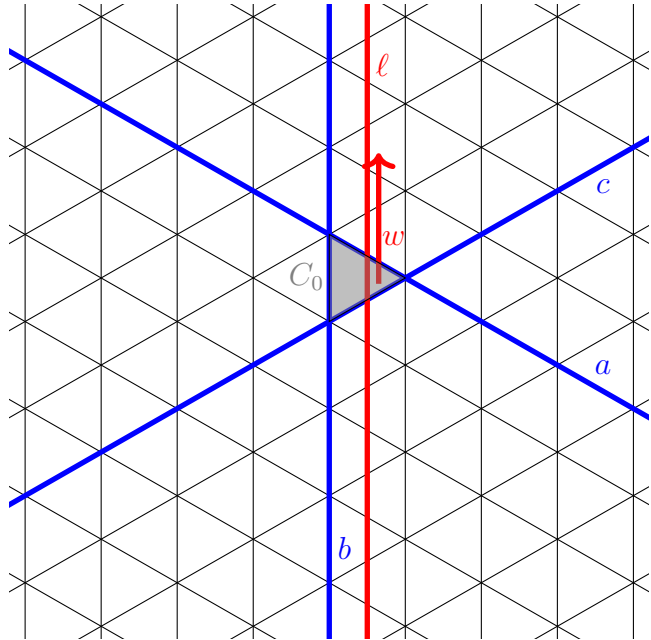

\begin{figure}
\begin{center}
\begin{tikzpicture}
\def \p {0.05}
\def \op {0.5}
\def \gris {black!50}
\clip (-4.2,-4.2) rectangle (4.2,4.2);

\foreach \i in {-10,...,10}
\draw (2*\i,-10) -- (2*\i,10);

\foreach \i in {-10,...,10}
\draw (-17.32-\i,-10+\i*2*0.866) -- (17.32-\i,10+\i*2*0.866);

\foreach \i in {-10,...,10}
\draw (-17.32-\i,10-\i*2*0.866) -- (17.32-\i,-10-\i*2*0.866);

\foreach \i in {-10,...,10}
\draw (-10+\i,-17.32-\i/1.732) -- (10+\i,17.32-\i/1.732);

\foreach \i in {-10,...,10}
\draw (10-\i,-17.32-\i/1.732) -- (-10-\i,17.32-\i/1.732);

\foreach \i in {-10,...,10}
\draw (-10,2*\i/1.732) -- (10,2*\i/1.732);

\draw[line width=2, color=blue] (-10,0) -- (10,0);
\draw[line width=2, color=blue] (-17.32,-10) -- (17.32,10);
\draw[line width=2, color=blue] (10+1,-17.32+1/1.732) -- (-10+1,17.32+1/1.732);

\node[thick, color=blue] at (2.7,-2.8) {$b$};
\node[thick, color=blue] at (3.1,-0.2) {$a$};
\node[thick, color=blue] at (3.4,1.7) {$c$};

\draw[line width=2, color=red] (1,-10) -- (1,10);
\node[thick, color=red] at (1.2,3.7) {$\ell$};
\draw[line width=2, color=red,->] (1.2,0.1) -- (1.2,1.3);
\node[thick, color=red] at (1.4,1) {$w$};

\draw[black,fill opacity=\op,fill=\gris] (0,0) -- (1,1/1.732) -- (4/3,0) -- cycle;
\node[thick, color=\gris] at (0.3,0.6) {$C_0$};

\end{tikzpicture}
\caption{Un élément de Coxeter pour le groupe de Coxeter affine de type $\tilde{G}_2$.}
\label{fig:G2tilde_element_coxeter}
\end{center}
\end{figure}
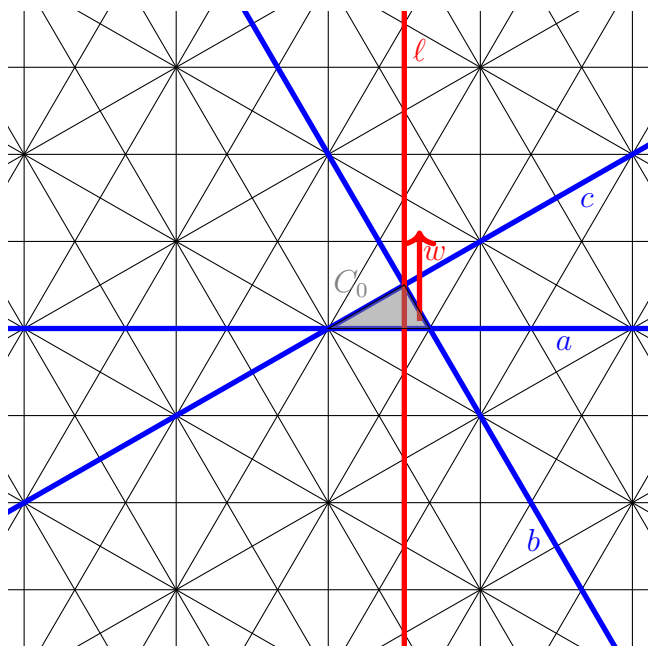

Remarquons que l'ensemble ordonné $[1,w]^W$ est borné, gradué, de rang $|S|$. De plus, les éléments de $S$ apparaissent tous comme étiquettes d'arêtes du graphe $[1,w]^W$ \parencite[Lemma~5.1]{PS}, ce qui permet de définir un morphisme naturel du groupe d'Artin usuel $G_W$ vers le groupe d'Artin dual $W_w$ associé à $w$.

\mk

Bessis a étudié le cas où $W$ est un groupe de Coxeter fini, et McCammond et Sulway ont étendu ce résultat à tous les groupes de Coxeter affines :

\bthm[\cite{bessis,brady_watt_Kpi1,mccammond_sulway}] \label{thm:dual_artin_isomorphic}
Si $W$ est un groupe de Coxeter fini ou affine, le morphisme naturel $G_W \ra W_w$ est un isomorphisme pour tout élément de Coxeter $w$.
\ethm

\bex
Dans le cas très simple du groupe de Coxeter $\mathfrak{S}_3$, pour l'élément de Coxeter $w=(1,2,3)$ correpondant au $3$-cycle, la présentation du groupe d'Artin dual est donnée par
$$\<\sigma_{1,2},\sigma_{1,3},\sigma_{2,3} \st \sigma_{1,2}\sigma_{2,3}=\sigma_{2,3}\sigma_{1,3}=\sigma_{1,3}\sigma_{1,2}\>,$$
dont il est facile de se convaincre qu'elle donne un groupe isomorphe au groupe de tresses~$B_3$, dont la présentation standard est
$$\<\sigma_1,\sigma_2 \st \sigma_1 \sigma_2 \sigma_1 = \sigma_2 \sigma_1 \sigma_2\>.$$
On peut en effet identifier $\sigma_{1,2}$ à $\sigma_1$, $\sigma_{2,3}$ à $\sigma_2$, et $\sigma_{1,3}$ à la tresse $\sigma_2^{-1} \sigma_1 \sigma_2$ permutant les brins $1$ et $3$.
\eex

On ne sait pas si ce résultat peut s'étendre à d'autres groupes de Coxeter que ceux finis ou affines.

\mk

Le théorème~\ref{thm:dual_artin_isomorphic} permet dans certains cas de munir le groupe d'Artin dual d'une structure de Garside. Pour cela, nous avons vu que la propriété clé est de savoir si l'intervalle $[1,w]^W$ est un treillis.

\bthm[\cite{brady_watt_NCP,bessis}]
Si $W$ est un groupe de Coxeter fini, l'intervalle $[1,w]^W$ est un treillis pour tout élément de Coxeter $w$. En particulier, le groupe d'Artin dual $W_w$ est un groupe de Garside.
\ethm

Si $W$ est un groupe de Coxeter quelconque et $w$ est un élément de Coxeter, on appelle ainsi $[1,w]^W$ l'ensemble ordonné des partitions non croisées associé à $W$. Il s'avère qu'à part dans le cas où $W$ est fini, il est rare que ce soit un treillis.

\bthm[\cite{digne_An,digne_Cn,mccammond_dual_euclidian}]
Soit $W$ un groupe de Coxeter irréductible affine, et soit $w$ un élément de Coxeter. Alors l'intervalle $[1,w]^W$ est un treillis si et seulement si $W$ est de type $\tilde{A}_n$ (lorsque $w$ est un $(n,1)$-bigone), $\tilde{C}_n$ ou $\tilde{G}_2$.
\ethm

Plus précisément, l'intervalle $[1,w]^W$ est un treillis si et seulement si le système de racines horizontal est irréductible, voir la table~\ref{tab:decomposition_horizontale}.

\mk

Cependant, McCammond et Sulway ont découvert comment compléter un groupe de Coxeter affine $W$ en un groupe d'isométries euclidiennes~$C$ plus grand, tel que l'intervalle $[1,w]^C \supset [1,w]^W$ soit un treillis équilibré. Le groupe d'intervalle associé à $[1,w]^C$, noté~$C_w$, est appelé \emph{groupe cristallographique tressé}, c'est un groupe de Garside. Ceci permet également de réaliser le groupe d'Artin dual~$W_w$ comme sous-groupe du groupe de~$C_w$.

\mk

Grâce aux groupes cristallographiques tressés, Paolini et Salvetti ont pu montrer que, même en l'absence de treillis, le complexe d'intervalle est un espace classifiant.

\bthm[{\cite[Theorem~6.6]{PS}}] \label{thm:artin_dual_classifiant}
Soit $W$ un groupe de Coxeter affine irréductible, et soit $w$ un élément de Coxeter. Le complexe d'intervalle $K_W$ associé à $[1,w]^W$ est un espace classifiant pour le groupe d'Artin dual $W_w$.
\ethm

Nous allons donner les idées de la preuve de ce théorème dans la partie~\ref{sec:crystallographic_braid_groups}, qui repose sur l'utilisation des groupes cristallographiques tressés.

\section{Un modèle dual pour l'espace de configuration}

\label{sec:dual_model_configuration_space}

Fixons un groupe de Coxeter quelconque $W$. Nous allons décrire un CW-complexe modèle de l'espace de configuration $Y_W$ associé à $W$.

\mk

Soit $S$ un ensemble simple de réflexions de $W$. Notons
$$\Delta_W = \{T \subset S \st W_T \mbox{ est fini}\}.$$

\mk

Rappelons la définition du complexe de Salvetti de $W$ \parencite{salvetti,paris_kpi1}. Pour tout $T \in \Delta_W$, notons $D_T$ le polytope de Coxeter associé à $W_T$: il peut être défini comme l'enveloppe convexe d'un point générique $x \in \R^T$ dans la représentation de $W_T$ agissant par réflexions orthogonales sur $\R^T$. Les faces de $D_T$ sont en correspondance avec les classes de $W_T / W_U$, où $U \subset T$, la face correspondant à $wW_U \in W_T/W_U$ s'identifiant avec $w \cdot D_U \subset D_T$ (où l'on voit $w \cdot D_U$ comme l'enveloppe convexe de $wW_U \cdot x$ dans~$D_T$).
\mk

Le complexe de Salvetti $X_W$ est le CW-complexe fini, quotient de la réunion disjointe $\bigsqcup_{T \in \Delta_W} D_T$ par les identifications suivantes : si $T \in \Delta_W$, $U \subset T$ et $wW_U \subset W_T / W_U$, alors la face $w D_U$ de $D_T$ est identifiée à $D_U$. Ce complexe a un seul sommet, correspondant à~$D_\emptyset$, et $S$ arêtes, correspondant à $D_{\{s\}}$, pour $s \in S$. On peut ainsi voir $X_W$ comme le recollement des complexes de Salvetti $X_T$, pour $T \in \Delta_W$.

\bex
Si $W$ est un groupe de Coxeter diédral d'ordre $2p$, avec $S=\{a,b\}$, alors le polytope de Coxeter $D_S$ est par exemple le $2p$-gone régulier. Il a un sommet privilégié correspondant à l'élément neutre $1$, et le sommet opposé correspondant à l'élément $[aba\cdots]_p=aba\cdots$ (le mot de longueur $p$). Les deux chemins reliant ces sommets sont étiquetés $[aba\cdots]_p$ et $[bab\cdots]_p$ respectivement.

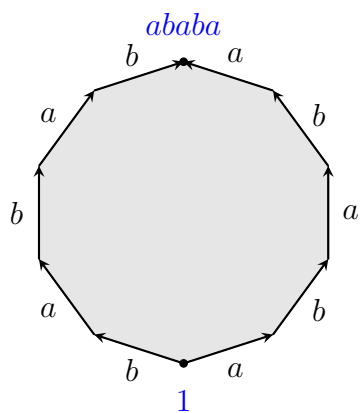
\begin{figure}
\begin{center}
\begin{tikzpicture}
\def \p {0.05}
\def \op {0.1}
\def \gris {black}
\def \r {2}

\draw[fill] (270+0:\r) circle (\p) node(e) {};
\draw[fill] (270+180:\r) circle (\p) node(P) {};
\draw[black,fill opacity=\op,fill=\gris,->] (270+0:\r) -- (270+36:\r) -- (270+72:\r) -- (270+108:\r) -- (270+144:\r) -- (270+180:\r) -- (270+216:\r) -- (270+252:\r) -- (270+288:\r) -- (270+324:\r) -- cycle;

\draw[thick,-stealth] (270+0:\r) -- (270+36:\r);
\draw[thick,-stealth] (270+36:\r) -- (270+72:\r);
\draw[thick,-stealth] (270+72:\r) -- (270+108:\r);
\draw[thick,-stealth] (270+108:\r) -- (270+144:\r);
\draw[thick,-stealth] (270+144:\r) -- (270+180:\r);
\draw[thick,stealth-] (270+180:\r) -- (270+216:\r);
\draw[thick,stealth-] (270+216:\r) -- (270+252:\r);
\draw[thick,stealth-] (270+252:\r) -- (270+288:\r);
\draw[thick,stealth-] (270+288:\r) -- (270+324:\r);
\draw[thick,stealth-] (270+324:\r) -- (270+0:\r);
\node[thick,blue] at (0,-\r-0.5) {$1$};
\node[thick,blue] at (0,\r+0.5) {$ababa$};

\node at (270+18:\r+0.2) {$a$};
\node at (270+36+18:\r+0.2) {$b$};
\node at (270+72+18:\r+0.2) {$a$};
\node at (270+108+18:\r+0.2) {$b$};
\node at (270+144+18:\r+0.2) {$a$};
\node at (270+180+18:\r+0.2) {$b$};
\node at (270+216+18:\r+0.2) {$a$};
\node at (270+252+18:\r+0.2) {$b$};
\node at (270+288+18:\r+0.2) {$a$};
\node at (270+324+18:\r+0.2) {$b$};

\end{tikzpicture}
\caption{Le complexe de Salvetti $X_W$ du groupe diédral $W$ d'ordre $2 \times 5$, obtenu en identifiant les arêtes selon leur étiquetage.}
\label{fig:complexe_salvetti_diedral}
\end{center}
\end{figure}

Le complexe de Salvetti associé à $W$ est le recollement des arêtes de $D_S$ selon leurs étiquettes, voir figure~\ref{fig:complexe_salvetti_diedral}.
\eex

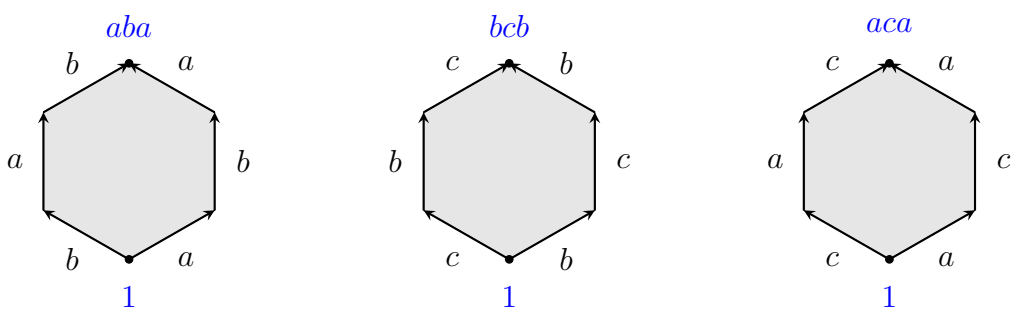
\begin{figure}
\begin{center}
\begin{tikzpicture}
\def \p {0.05}
\def \op {0.1}
\def \gris {black}
\def \r {1.3}
\def \t {4}

\draw[fill] (270+0:\r) circle (\p) node(e) {};
\draw[fill] (270+180:\r) circle (\p) node(P) {};
\draw[black,fill opacity=\op,fill=\gris,->] (270+0:\r) -- (270+60:\r) -- (270+120:\r) -- (270+180:\r) -- (270+240:\r) -- (270+300:\r) -- cycle;

\draw[thick,-stealth] (270+0:\r) -- (270+60:\r);
\draw[thick,-stealth] (270+60:\r) -- (270+120:\r);
\draw[thick,-stealth] (270+120:\r) -- (270+180:\r);
\draw[thick,stealth-] (270+180:\r) -- (270+240:\r);
\draw[thick,stealth-] (270+240:\r) -- (270+300:\r);
\draw[thick,stealth-] (270+300:\r) -- (270+0:\r);
\node[thick,blue] at (0,-\r-0.5) {$1$};
\node[thick,blue] at (0,\r+0.5) {$bcb$};

\node at (270+30:\r+0.2) {$b$};
\node at (270+90:\r+0.2) {$c$};
\node at (270+150:\r+0.2) {$b$};
\node at (270+210:\r+0.2) {$c$};
\node at (270+270:\r+0.2) {$b$};
\node at (270+330:\r+0.2) {$c$};

\begin{scope}[xshift=5cm]
\draw[fill] (270+0:\r) circle (\p) node(e) {};
\draw[fill] (270+180:\r) circle (\p) node(P) {};
\draw[black,fill opacity=\op,fill=\gris,->] (270+0:\r) -- (270+60:\r) -- (270+120:\r) -- (270+180:\r) -- (270+240:\r) -- (270+300:\r) -- cycle;

\draw[thick,-stealth] (270+0:\r) -- (270+60:\r);
\draw[thick,-stealth] (270+60:\r) -- (270+120:\r);
\draw[thick,-stealth] (270+120:\r) -- (270+180:\r);
\draw[thick,stealth-] (270+180:\r) -- (270+240:\r);
\draw[thick,stealth-] (270+240:\r) -- (270+300:\r);
\draw[thick,stealth-] (270+300:\r) -- (270+0:\r);
\node[thick,blue] at (0,-\r-0.5) {$1$};
\node[thick,blue] at (0,\r+0.5) {$aca$};

\node at (270+30:\r+0.2) {$a$};
\node at (270+90:\r+0.2) {$c$};
\node at (270+150:\r+0.2) {$a$};
\node at (270+210:\r+0.2) {$c$};
\node at (270+270:\r+0.2) {$a$};
\node at (270+330:\r+0.2) {$c$};
\end{scope}

\begin{scope}[xshift=-5cm]
\draw[fill] (270+0:\r) circle (\p) node(e) {};
\draw[fill] (270+180:\r) circle (\p) node(P) {};
\draw[black,fill opacity=\op,fill=\gris,->] (270+0:\r) -- (270+60:\r) -- (270+120:\r) -- (270+180:\r) -- (270+240:\r) -- (270+300:\r) -- cycle;

\draw[thick,-stealth] (270+0:\r) -- (270+60:\r);
\draw[thick,-stealth] (270+60:\r) -- (270+120:\r);
\draw[thick,-stealth] (270+120:\r) -- (270+180:\r);
\draw[thick,stealth-] (270+180:\r) -- (270+240:\r);
\draw[thick,stealth-] (270+240:\r) -- (270+300:\r);
\draw[thick,stealth-] (270+300:\r) -- (270+0:\r);
\node[thick,blue] at (0,-\r-0.5) {$1$};
\node[thick,blue] at (0,\r+0.5) {$aba$};

\node at (270+30:\r+0.2) {$a$};
\node at (270+90:\r+0.2) {$b$};
\node at (270+150:\r+0.2) {$a$};
\node at (270+210:\r+0.2) {$b$};
\node at (270+270:\r+0.2) {$a$};
\node at (270+330:\r+0.2) {$b$};
\end{scope}

\end{tikzpicture}
\caption{Le complexe de Salvetti $X_W$ du groupe diédral affine $W$ de type $\tilde{A}_2$, obtenu en identifiant les arêtes selon leur étiquetage.}
\label{fig:complexe_salvetti_affine_A2tilde}
\end{center}
\end{figure}

\bex
Si $W$ est un groupe de Coxeter affine de type $\tilde{A}_2$, avec $S=\{a,b,c\}$, alors $\Delta_W$ a trois sous-ensembles maximaux $\{a,b\}$, $\{b,c\}$ et $\{a,c\}$. Le complexe de Salvetti de $W$ est donc le recollement des trois hexagones $D_{a,b}$, $D_{b,c}$ et $D_{a,c}$ selon les étiquettes de leurs arêtes, voir figure~\ref{fig:complexe_salvetti_affine_A2tilde}.
\eex

L'un des intérêts du complexe de Salvetti est de fournir un modèle combinatoire fini pour l'espace de configuration.

\bthm[\cite{salvetti,salvetti_homotopy}]
Pour tout groupe de Coxeter $W$, le complexe de Salvetti $X_W$ a le même type d'homotopie que l'espace de configuration $Y_W$.
\ethm

Soit $R$ l'ensemble des réflexions de $W$, et $w=s_1s_2\cdots s_n$ un élément de Coxeter obtenu en écrivant le produit des éléments de $S$ dans un certain ordre. Notons $K_W$ le complexe d'intervalle asssocié à $[1,w]^W$. Pour tout $T \in \Delta_W$, remarquons qu'on peut considérer un élément de Coxeter privilégié $w_T$ de $W_T$, celui obtenu en prenant le produit des éléments de $T$ dans le même ordre que ceux apparaissant dans $w$.

\mk

Nous allons décrire un sous-complexe $X'_W$ de $K_W$ qui aura le type d'homotopie de~$X_W$. 

\bdf
Notons $X'_W$ le sous-complexe de $K_W$ constitué des simplexes $[x_1|x_2|\dots|x_d]$ de $K_W$ tels qu'il existe $T \in \Delta_W$ pour lequel $x_1x_2\cdots x_d \in [1,w_T]^W=[1,w_T]^{W_T}$. Nous proposons de l'appeler \emph{complexe de Salvetti dual} de $W$.
\edf

\bex
Soit $W$ un groupe de Coxeter affine de type $\tilde{A}_2$, avec $S=\{a,b,c\}$, et considérons comme élément de Coxeter $w=abc$. Alors le complexe de Salvetti dual $X'_W$ est le recollement des trois complexes d'intervalle associés aux intervalle $[1,ab]$, $[1,bc]$ et $[1,ac]$ selon les étiquettes des arêtes, voir figure~\ref{fig:complexe_salvetti_dual_artin_affine_A2_tilde}. On pourra noter la similitude avec la figure~\ref{fig:complexe_salvetti_affine_A2tilde}.
\eex

\mk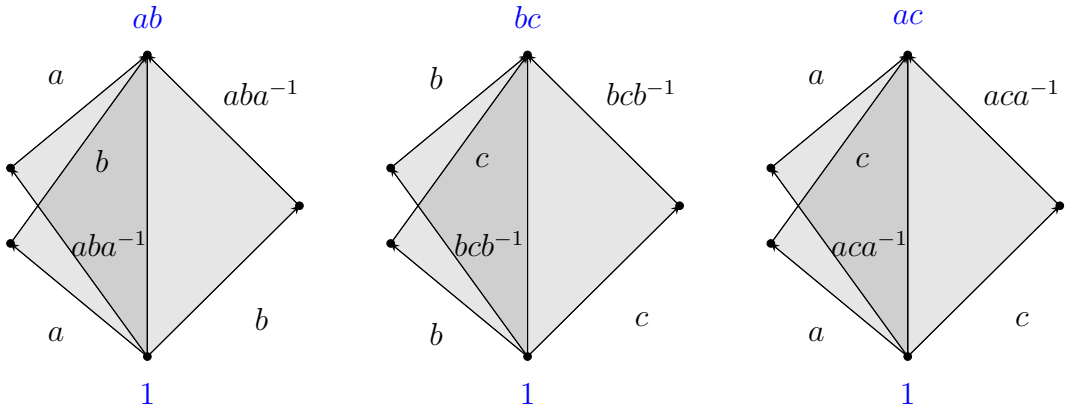
\begin{figure}
\begin{center}
\begin{tikzpicture}
\def \p {0.05}
\def \op {0.1}
\def \gris {black}

\draw[fill] (0,0) circle (\p) node(e) {};
\draw[fill] (-1.8,1.5) circle (\p) node(a) {};
\draw[fill] (2,2) circle (\p) node(b) {};
\draw[fill] (-1.8,2.5) circle (\p) node(c) {};
\draw[fill] (0,4) circle (\p) node(delta) {};

\draw[-stealth] (e.center) -- (a.center);
\draw[-stealth] (e.center) -- (b.center);
\draw[-stealth] (e.center) -- (c.center);
\draw[-stealth] (a.center) -- (delta.center);
\draw[-stealth] (b.center) -- (delta.center);
\draw[-stealth] (c.center) -- (delta.center);
\draw[black,fill opacity=\op,fill=\gris,->] (e.center) -- (a.center) -- (delta.center) -- cycle;
\draw[black,fill opacity=\op,fill=\gris] (e.center) -- (c.center) -- (delta.center) -- cycle;
\draw[black,fill opacity=\op,fill=\gris] (e.center) -- (b.center) -- (delta.center) -- cycle;

\node at (-1.2,0.3) {$b$};
\node at (-1.2,3.7) {$b$};
\node at (1.5,0.5) {$c$};
\node at (1.5,3.5) {$bcb^{-1}$};
\node at (-0.5,1.5) {$bcb^{-1}$};
\node at (-0.6,2.6) {$c$};

\node[blue] at (0,-0.5) {$1$};
\node[blue] at (0,4.5) {$bc$};

\begin{scope}[xshift=-5cm]
\draw[fill] (0,0) circle (\p) node(e) {};
\draw[fill] (-1.8,1.5) circle (\p) node(a) {};
\draw[fill] (2,2) circle (\p) node(b) {};
\draw[fill] (-1.8,2.5) circle (\p) node(c) {};
\draw[fill] (0,4) circle (\p) node(delta) {};

\draw[-stealth] (e.center) -- (a.center);
\draw[-stealth] (e.center) -- (b.center);
\draw[-stealth] (e.center) -- (c.center);
\draw[-stealth] (a.center) -- (delta.center);
\draw[-stealth] (b.center) -- (delta.center);
\draw[-stealth] (c.center) -- (delta.center);
\draw[black,fill opacity=\op,fill=\gris] (e.center) -- (a.center) -- (delta.center) -- cycle;
\draw[black,fill opacity=\op,fill=\gris] (e.center) -- (c.center) -- (delta.center) -- cycle;
\draw[black,fill opacity=\op,fill=\gris] (e.center) -- (b.center) -- (delta.center) -- cycle;

\node at (-1.2,0.3) {$a$};
\node at (-1.2,3.7) {$a$};
\node at (1.5,0.5) {$b$};
\node at (1.5,3.5) {$aba^{-1}$};
\node at (-0.5,1.5) {$aba^{-1}$};
\node at (-0.6,2.6) {$b$};

\node[blue] at (0,-0.5) {$1$};
\node[blue] at (0,4.5) {$ab$};
\end{scope}

\begin{scope}[xshift=5cm]
\draw[fill] (0,0) circle (\p) node(e) {};
\draw[fill] (-1.8,1.5) circle (\p) node(a) {};
\draw[fill] (2,2) circle (\p) node(b) {};
\draw[fill] (-1.8,2.5) circle (\p) node(c) {};
\draw[fill] (0,4) circle (\p) node(delta) {};

\draw[-stealth] (e.center) -- (a.center);
\draw[-stealth] (e.center) -- (b.center);
\draw[-stealth] (e.center) -- (c.center);
\draw[-stealth] (a.center) -- (delta.center);
\draw[-stealth] (b.center) -- (delta.center);
\draw[-stealth] (c.center) -- (delta.center);
\draw[black,fill opacity=\op,fill=\gris] (e.center) -- (a.center) -- (delta.center) -- cycle;
\draw[black,fill opacity=\op,fill=\gris] (e.center) -- (c.center) -- (delta.center) -- cycle;
\draw[black,fill opacity=\op,fill=\gris] (e.center) -- (b.center) -- (delta.center) -- cycle;

\node at (-1.2,0.3) {$a$};
\node at (-1.2,3.7) {$a$};
\node at (1.5,0.5) {$c$};
\node at (1.5,3.5) {$aca^{-1}$};
\node at (-0.5,1.5) {$aca^{-1}$};
\node at (-0.6,2.6) {$c$};

\node[blue] at (0,-0.5) {$1$};
\node[blue] at (0,4.5) {$ac$};
\end{scope}

\end{tikzpicture}
\caption{Le complexe de Salvetti dual $X'_W$ du groupe diédral affine $W$ de type $\tilde{A}_2$.}
\label{fig:complexe_salvetti_dual_artin_affine_A2_tilde}
\end{center}
\end{figure}

Remarquons que si $W$ est fini, alors $X'_W=K_W$. Ainsi, si $T \in \Delta_W$, on voit que $X'_{W_T}$ est un espace classifiant pour le groupe d'Artin $G_{W_T}$. De même, comme la conjecture du $K(\pi,1)$ est connue pour les groupes d'Artin sphériques \parencite{deligne_immeubles}, on sait que $X_{W_T}$ est aussi un espace classifiant pour le groupe d'Artin $G_{W_T}$.

\bthm[{\cite[Theorem~5.5]{PS}}] \label{thm:YW_XprimeW}
Pour tout groupe de Coxeter $W$, le complexe de Salvetti dual $X'_W \subset K_W$ a le même type d'homotopie que le complexe de Salvetti $X_W$ et que l'espace de configuration $Y_W$.
\ethm

\bp
Pour alléger les notations, dans cette preuve nous allons noter $X,X_T,X',X'_T$ à la place de $X_W,X_{W_T},X'_W,X'_{W_T}$, pour $T \in \Delta_W$. Nous avons déjà remarqué que, si $T \in \Delta_W$, les complexes~$X_T$ et~$X'_T$ sont tous deux des espaces classifiants pour le groupe d'Artin sphérique $G_{W_T}$. Nous allons construire, par induction sur~$|T|$, des équivalences d'homotopie $\phi_T \colon X_T \ra X'_T$ telles que, pour tout $U \subset T \in \Delta_W$, nous ayons le diagramme commutatif suivant :
\beq X_U & \stackrel{\phi_U}{\ra} & X'_U \\
\rotatebox[origin=c]{-90}{$\hookrightarrow$} & & \rotatebox[origin=c]{-90}{$\hookrightarrow$} \\
X_T & \stackrel{\phi_T}{\ra} & X'_T \eeq
Supposons que $T \in \Delta_W$ soit tel que nous ayons déjà défini de telles applications $\phi_U$, pour $U \subsetneq T$.
\bit
\item Si $T=\emptyset$, alors $X_T$ et $X'_T$ sont constitués d'un unique sommet, il n'y a donc qu'une seule application $\phi_T \colon X_T \ra X'_T$.
\item Si $T=\{s\}$, alors $X_T$ et $X'_T$ sont constitués d'une arête orientée
  étiquetée $s$ attachée à l'unique sommet. Considérons un homorphisme
  cellulaire $\phi_T\colon X_T \ra X'_T$ préservant l'orientation.
\item Si $T=\{s,s'\}$, alors d'après la preuve de la proposition~1B.9 de
  \textcite{hatcher}, l'application $\phi_{\{s\}} \cup \phi_{\{s'\}}\colon
  X_{\{s\}} \cup X_{\{s'\}} \ra X'_T$ peut être étendue en une application
  $\phi_T\colon X_T \ra X'_T$ telle que l'application induite $(\phi_T)_\star \colon \pi_1(X_T,X_\emptyset) \ra \pi_1(X'_T,X'_\emptyset)$ soit un isomorphisme. Comme $X_T$ et $X'_T$ sont des espaces classifiants, on en déduit que $\phi_T$ est une équivalence d'homotopie.
\item Si $|T| \geq 3$, en utilisant comme précédemment la preuve de la
  proposition~1B.9 de \textcite{hatcher}, on construit de même $\phi_T$.
\eit
Cette construction induit ainsi une application $\phi_W\colon X_W \ra X'_W$. En
appliquant successivement le théorème~7.5.7 de \textcite{brown}, on déduit que $\phi_W$ est une équivalence d'homotopie.
\ep

Remarquons qu'on peut également décrire, à homotopie près, le complexe de Salvetti comme un recollement de complexes d'intervalle analogues à $[1,w_T]^{W_T}$, mais pour la structure de Garside classique de $W_T$.

Ceci permet ainsi de justifier la dénomination du complexe $X'_W$ comme complexe de Salvetti dual.

\section{Factorisations d'isométries euclidiennes}

\label{sec:factorisation_euclidean_isometries}

Afin de comprendre l'ensemble des factorisations d'un élément de Coxeter d'un groupe de Coxeter affine, il est utile d'étudier plus généralement les factorisations d'isométries euclidiennes quelconques comme produits de réflexions.

\subsection{Un ordre sur les isométries euclidiennes}

\label{subsec:ordre_isometries}

Notons $V \simeq \R^n$ un espace vectoriel réel de dimension~$n$. Nous allons décrire un ordre sur le groupe~$L$ de toutes les isométries euclidiennes affines de~$V$.

\mk

Une isométrie $u \in L$ est appelée \emph{elliptique} si elle fixe au moins un point de $V$.  Dans ce cas, notons $\Fix(u) \subset V$ l'ensemble des points fixes de $u$.

Une isométrie $u \in L$ est appelée \emph{hyperbolique} si elle ne fixe aucun point de $V$. Dans ce cas, notons $\Dep(u)=\{u(a)-a \st a \in V\} \subset V$ l'ensemble des déplacements de $u$. C'est un sous-espace affine de $V$, ayant un unique vecteur $\mu$ de norme minimale. Notons $\Min(u)=\{a \in V \st u(a)=a+\mu\}$ l'ensemble des points de $E$ le moins déplacés par $u$.

\mk

Le groupe $L=\Isom(V)$ est engendré par l'ensemble $R$ de toutes les réflexions orthogonales de $V$ (avec un poids constant égal à $1$). On peut ainsi calculer le plus petit nombre de réflexions nécessaires $l(u)$ pour écrire une isométrie $u$ donnée.

\bpro[{\cite[Theorem~5.7]{brady_mccammond_factoring}}]
Si $u \in L$ est elliptique, alors $l(u) = \codim \Fix(u)$. Si $u$ est hyperbolique, alors $l(u)=\dim \Dep(u) + 2$.
\epro

On peut ainsi représenter l'ordre associé sur $L$ grâce au modèle simple suivant.

\bdf
Soit $(P,\leq)$ l'ensemble ordonné constitué d'un élément noté $e^F$ pour chaque sous-espace affine $F \subset L$, dont la direction est notée $\vec{F}$, ainsi qu'un élément noté $h^D$ pour chaque sous-espace affine $D \subset V \smallsetminus \{0\}$. L'ordre sur $P$ est défini comme suit :
\bit
\item $e^F \leq e^{F'}$ si et seulement si $F' \subset F$,
\item $h^D \leq h^{D'}$ si et seulement si $D \subset D'$ et
\item $e^F \leq h^D$ si et seulement si $D^\perp \subset \vec{F}$. 
\eit
Remarquons que $P$ a pour élément minimal $e^V$. Considérons l'application invariant, notée $inv$ et définie par
\beq L & \ra & P \\
v \in L \mbox{ elliptique } & \mapsto & e^{\Fix(v)} \\
v \in L \mbox{ hyperbolique } & \mapsto & h^{\Dep(v)}.\eeq
\edf

Ce modèle permet de décrire simplement l'intervalle en-dessous de toute isométrie.

\bthm[{\cite[Theorem~8.7]{brady_mccammond_factoring}}]
Pour toute isométrie $u \in L$, l'application invariant est un isomorphisme d'ensembles ordonnés entre l'intervalle $[1,u]^L$ et l'intervalle $[e^V,inv(u)]^P$.
\ethm

\subsection{Factorisations dans les groupes de Coxeter affines}

Soit $W$ un groupe de Coxeter affine irréductible, agissant comme groupe de réflexions sur l'espace euclidien $V=\R^n$, où $n$ est le rang de $W$. Soit $R$ l'ensemble de ses réflexions, et soit $S \subset R$ un ensemble simple de réflexions. Rappelons qu'un élément de Coxeter $w$ est le produit des éléments de $S$ dans un ordre quelconque.

\mk

Alors l'isométrie $w$ de $V$ est hyperbolique, de longueur $l(w)=n+1$, et $\Min(w)$ est une droite $\ell$ appelée \emph{axe de Coxeter}.

\mk

Considérons l'arrangement $\mathcal{A}$ de tous les hyperplans des réflexions de $R$ dans $V$. Les composantes connexes du complémentaire $V \smallsetminus \cup \mathcal{A}$ sont appelées les \emph{chambres} (ouvertes) de $\mathcal{A}$, qui sont les simplexes maximaux d'une structure simpliciale sur $V$. On pourra ainsi parler des sommets de cette structure.

\mk

Une chambre qui intersecte l'axe de Coxeter est appelée \emph{chambre axiale}, et ses sommets sont appelés \emph{sommets axiaux}. Remarquons que les hyperplans fixés par les réflexions de~$S$ bordent une chambre axiale particulière notée~$C_0$, voir les figures~\ref{fig:A2tilde_element_coxeter} et \ref{fig:G2tilde_element_coxeter} pour les types~$\tilde{A}_2$ et~$\tilde{G}_2$.

\mk

L'ordre sur le groupe $L$ des isométries de $V$ défini dans la partie~\ref{subsec:ordre_isometries} peut être comparé avec celui de l'intervalle $[1,w]^W$ du groupe de Coxeter $W$ :

\blem[{\cite[Lemma~2.15]{PS}}]
L'inclusion $[1,w]^W \ra [1,w]^L$ préserve l'ordre et le rang.
\elem

Cependant, l'ordre sur $W$ n'est pas nécessairement induit par celui de $L$.

\bdf[Vertical / Horizontal]
La direction de l'axe de Coxeter $\ell$ est appelée \emph{verticale}, et les directions orthogonales à $\ell$ sont appelées \emph{horizontales}. Une isométrie elliptique $u$ est appelée \emph{horizontale} si elle déplace chaque point dans une direction horizontale (i.e. $\ell \subset \Fix(u)$), et \emph{verticale} sinon.
\edf

Cette notion permet de décrire sommairement les éléments de l'intervalle $[1,w]^W$, et donne une ébauche de l'ordre que l'on va construire sur l'ensemble des réflexions.

\bpro[{\cite[Proposition~2.17]{PS}}] \label{pro:factorisations_elliptiques}
Soit $W$ un groupe de Coxeter affine irréductible, et soit $w$ un élément de Coxeter. Les éléments $u \in [1,w]^W$ se répartissent dans trois rangées, selon les cas suivants (où $v \in W$ tel que $uv=w$ est le complément à droite de $u$) :
\bit
\item (Rangée du bas) $u$ est elliptique horizontal et $v$ est hyperbolique,
\item (Rangée du milieu) $u$ et $v$ sont elliptiques verticaux,
\item (Rangée du haut) $u$ est hyperbolique et $v$ est elliptique horizontal.
\eit
De plus, les rangées du haut et du bas sont finies, tandis que la rangée du milieu est infinie.
\epro

On peut même décrire précisément quelles réflexions de $W$ apparaissent dans l'intervalle $[1,w]^W$.

\bthm[{\cite[Theorem~9.6]{mccammond_dual_euclidian}; \cite[Theorem~3.17]{PS}}] \label{thm:reflection_sommets_axiaux}
Soit $W$ un groupe de Coxeter affine irréductible, et soit $w$ un élément de Coxeter. Toute réflexion verticale $r \in W$ est dans $[1,w]^W$, et fixe au moins deux sommets axiaux. Une réflexion horizontale $r \in W$ appartient à $[1,w]^W$ si et seulement si $r$ fixe au moins un sommet axial.
\ethm

Nous noterons $R_0 = R \cap [1,w]^W$ l'ensemble des réflexions fixant au moins un sommet axial.

\mk

La dichotomie entre isométries hyperboliques et elliptiques permet même d'étudier de manière précise l'intervalle situé en-dessous d'une isométrie de $[1,w]^W$, à l'aide de la décomposition hyperbolique-horizontale suivante.

\mk

Si $u \in [1,w]^W$, notons $W_u$ le sous-groupe de $W$ engendré par $[1,u]^W$ : c'est un groupe de Coxeter, dont l'ensemble des réflexions est $R \cap W_u$ \parencite[Theorem~8.2]{humphreys}.

\bthm[{\cite[Lemma~3.20, Theorem~3.22]{PS}}] \label{thm:decomposition_horizontale_hyperbolique}
Soit $W$ un groupe de Coxeter affine irréductible, soit $w \in W$ un élément de Coxeter, et soit $u \in [1,w]^W$. Alors $u$ est un élément de Coxeter de $W_u$. De plus, si $u$ est hyperbolique, il existe une unique décomposition $u=u'h$ dans $W$ telle que :
\ben
\item $u'$ est hyperbolique, $h$ est elliptique horizontale, et $\ell(u)=\ell(u')+\ell(h)$,
\item $W_{u'}$ et $W_h$ commutent, et $W_u = W_{u'} \times W_h$,
\item $W_{u'}$ est un sous-groupe affine de Coxeter irréductible, et $u'$ en est un élément de Coxeter,
\item $W_h$ est un sous-groupe sphérique de Coxeter engendré par des réflexions horizontales, et $h$ en est un élément de Coxeter,
\item $[1,u]^W \simeq [1,u']^W \times [1,h]^W$ (où l'isomorphisme est donné par la multiplication).
\een
\ethm

Notons que la preuve de ce résultat repose sur un lemme technique, qui est
démontré pour les quatre familles infinies dans \textcite{PS}, et est démontré
par ordinateur pour les cas exceptionnels dans \textcite{paolini_computer}.

\subsection{Décomposition horizontale}

Nous allons maintenant décrire la structure des isométries horizontales d'un
groupe de Coxeter affine, qui ont essentiellement été étudiées dans \textcite{mccammond_sulway}.

Considérons un groupe de Coxeter affine irréductible $W$, agissant comme groupe de réflexions sur $\R^n$, où $n$ est le rang de $W$. Soit $w \in W$ un élément de Coxeter, $R$ l'ensemble des réflexions de $W$. 

Soit $\Phi$ le système de racines de $W$, et $\Phi_{hor}$ le système de
racines horizontal, c'est-à-dire l'ensemble des racines de $\Phi$ qui sont horizontales. D'après \textcite[Section~6]{mccammond_sulway}, ce système
se décompose en sous-systèmes de type $A$ : $\Phi_{hor}=\Phi_1 \sqcup \cdots
\sqcup \Phi_k$, où $\Phi_i$ est un système de racines de type $A_{n_i}$. Nous
renvoyons à la table~\ref{tab:decomposition_horizontale} pour la liste des
différentes décompositions possibles selon le type de $W$. Rappelons que tous
les éléments de Coxeter en type $\tilde{A}_n$ ne sont pas équivalents, et nous
renvoyons à \textcite[Definition~7.7,Example~11.6]{mccammond_dual_euclidian} pour la définition de $(p,q)$-bigone de Coxeter.

\begin{table} 
\begin{center}
\begin{tabular}{l|l}
Type de $W$ & Décomposition du système de racines horizontal \\
\hline
$\tilde{A}_n$, où $w$ est un $(p,q)$-bigone & $\Phi_{A_{p-1}} \sqcup \Phi_{A_{q-1}}$\\
$\tilde{B}_n$ & $\Phi_{A_1} \sqcup \Phi_{A_{n-2}}$\\
$\tilde{C}_n$ & $\Phi_{A_{n-1}}$\\
$\tilde{D}_n$ & $\Phi_{A_1} \sqcup \Phi_{A_1} \sqcup \Phi_{A_{n-3}}$\\
$\tilde{G}_2$ & $\Phi_{A_1}$\\
$\tilde{F}_4$ & $\Phi_{A_1} \sqcup \Phi_{A_2}$\\
$\tilde{E}_6$ & $\Phi_{A_1} \sqcup \Phi_{A_2} \sqcup \Phi_{A_2}$\\
$\tilde{E}_7$ & $\Phi_{A_1} \sqcup \Phi_{A_2} \sqcup \Phi_{A_3}$\\
$\tilde{E}_8$ & $\Phi_{A_1} \sqcup \Phi_{A_2} \sqcup \Phi_{A_4}$\\
\end{tabular}
\caption{Décomposition du système de racines horizontal $\Phi_{hor}$, d'après \textcite[Table~1]{mccammond_sulway}.}
\label{tab:decomposition_horizontale}
\end{center}
\end{table}

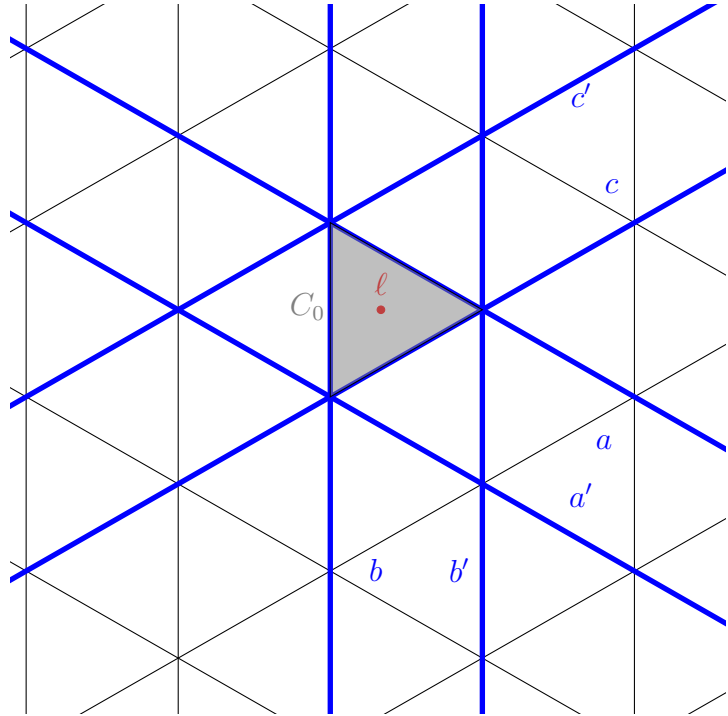
\begin{figure}
\begin{center}
\begin{tikzpicture}
\def \p {0.05}
\def \op {0.5}
\def \gris {black!50}
\clip (-4.2,-4.2) rectangle (5.2,5.2);

\foreach \i in {-10,...,10}
\draw (2*\i,-10) -- (2*\i,10);

\foreach \i in {-10,...,10}
\draw (-17.32-2*\i/2,-10+2*\i*0.866) -- (17.32-2*\i/2,10+2*\i*0.866);

\foreach \i in {-10,...,10}
\draw (-17.32-2*\i/2,10-2*\i*0.866) -- (17.32-2*\i/2,-10-2*\i*0.866);

\draw[fill, color=red] (0.666,2/1.732) circle (\p) node(l) {};
\node[thick, color=red] at (0.666,1.5) {$\ell$};

\draw[line width=2, color=blue] (0,-10) -- (0,10);
\draw[line width=2, color=blue] (-17.32,-10) -- (17.32,10);
\draw[line width=2, color=blue] (-2*17.32+1,20+2*0.866) -- (2*17.32+1,-20+2*0.866);

\draw[line width=2, color=blue] (2,-10) -- (2,10);
\draw[line width=2, color=blue] (-17.32,-10+4/1.732) -- (17.32,10+4/1.732);
\draw[line width=2, color=blue] (-2*17.32+1,20+2*0.866-4/1.732) -- (2*17.32+1,-20+2*0.866-4/1.732);

\node[thick, color=blue] at (0.6,-2.3) {$b$};
\node[thick, color=blue] at (1.7,-2.3) {$b'$};
\node[thick, color=blue] at (3.6,-0.6) {$a$};
\node[thick, color=blue] at (3.3,-1.3) {$a'$};
\node[thick, color=blue] at (3.7,2.8) {$c$};
\node[thick, color=blue] at (3.3,4) {$c'$};

\draw[black,fill opacity=\op,fill=\gris] (0,0) -- (2,2/1.732) -- (0,4/1.732) -- cycle;
\node[thick, color=\gris] at (-0.3,1.2) {$C_0$};

\end{tikzpicture}
\caption{Le système de racines horizontal en type $\tilde{C}_3$, et les réflexions horizontales de $[1,w]^W$.}
\label{fig:C3tilde_horizontal_roots}
\end{center}
\end{figure}

\bex
Considérons le groupe de Coxeter affine $W$ de type $\tilde{C_n}$, dont l'arrangement d'hyperplans $\mathcal{A}$ dans $\R^n$ est donné par les hyperplans $\{x_i = p\}_{1 \leq i \leq n, p \in \Z}$ et $\{x_i \pm x_j = p\}_{1 \leq i < j \leq n, p \in \Z}$. Un choix de chambre fondamentale $C_0$ est donné par le simplexe ouvert (aussi appelé orthosimplexe) :
$$C_0=\{x \in \R^n \st 0 < x_1 < x_2 < \cdots < x_n < 1\}.$$
Notons $S=\{s_1,s_{1,2},s_{2,3},\dots,s_{n-1,n},s_n\}$ l'ensemble simple de réflexions de~$W$, par rapport aux hyperplans
$\{x_1=0\},\{x_1=x_2\},\{x_2=x_3\},\dots,\{x_{n-1}=x_n\},\{x_n=1\}$ supportant les faces de~$C_0$.
Un élément de Coxeter est $w=s_1s_{1,2}s_{2,3} \cdots s_{n-1,n}s_n$, dont l'action sur~$\R^n$ est donnée par
$$\forall x \in \R^n, w \cdot x=(x_n-2,x_1,x_2,x_3,\dots,x_{n-1}),$$
ainsi l'axe de Coxeter est
$$\ell=\{a+\theta \mu \st \theta \in \R\}, \mbox{ où } a=\f{2}{n}(1,2,3,\dots,n) \mbox{ et } \mu=-\f{2}{n}(1,1,\dots,1).$$
Les réflexions horizontales sont celles dont l'hyperplan est orthogonal à $\mu$, c'est-à-dire $\{x_i - x_j = p\}_{1 \leq i < j \leq n, p \in \Z}$. Il s'agit bien d'un système de racines de type $\tilde{A}_{n-1}$. Voir la figure~\ref{fig:C3tilde_horizontal_roots} représentant pour le type $\tilde{C}_3$, dans un hyperplan orthogonal à $\ell$, le système de racines horizontal, ainsi que la trace de $\ell$ à l'intérieur de la chambre $C_0$.

Les réflexions horizontales appartenant à l'intervalle $[1,w]^W$ sont celles fixant au moins un sommet axial, d'après le théorème~\ref{thm:reflection_sommets_axiaux}. En type $\tilde{C}_3$, il y en a donc $6$, notées $a,a',b,b',c,c'$ sur la figure~\ref{fig:C3tilde_horizontal_roots}.
\eex

\bex \label{ex:B3tilde}
Considérons le groupe de Coxeter affine $W$ de type $\tilde{B_n}$, dont l'arrangement d'hyperplans $\mathcal{A}$ dans $\R^n$ est donné par les hyperplans $\{x_i = 2p\}_{1 \leq i \leq n, p \in \Z}$ et $\{x_i \pm x_j = p\}_{1 \leq i < j \leq n, p \in \Z}$. Un choix de chambre fondamentale $C_0$ est donné par le simplexe ouvert :
$$C_0=\{x \in \R^n \st 0 < x_1 < x_2 < \cdots < x_n, x_{n-1}+x_n<1\}.$$
Notons $S=\{s_1,s_{1,2},s_{2,3},\dots,s_{n-2,n-1},s_{n-2,n-1},s'_{n-2,n-1}\}$ l'ensemble simple de réflexions de $W$, par rapport aux hyperplans
$\{x_1=0\},\{x_1=x_2\},\{x_2=x_3\},\dots,\{x_{n-1}=x_n\},\{x_{n-1}+x_n=1\}$ supportant les faces de $C_0$.
Un élément de Coxeter est $w=s_1s_{1,2}s_{2,3} \cdots s_{n-1,n}s_{n-1,n}$, dont l'action sur $\R^n$ est donnée par
$$\forall x \in \R^n, w \cdot x=(x_{n-1}-1,x_1,x_2,\dots,x_{n-2},1-x_n),$$
ainsi l'axe de Coxeter est
\begin{equation*}
  \ell=\{a+\theta \mu \st \theta \in \R\},\, \mbox{où } a=\left(\f{1}{n-1},\f{2}{n-1},\dots,\f{n-1}{n-1},\f{1}{2}\right)\, \mbox{et } \mu=-\f{2}{n-1}(1,1,\dots,1,0).
\end{equation*}
Les réflexions horizontales sont celles dont l'hyperplan est orthogonal à $\mu$, c'est-à-dire $\{x_i - x_j = p\}_{1 \leq i < j \leq n-1, p \in \Z}$ et $\{x_n=2p\}_{p \in \Z}$. Il s'agit bien d'un système de racines réductible, de type $\tilde{A}_{n-2} \times \tilde{A}_1$, préservant la décomposition $\R^n = \R^{n-1} \times \R$.\eex

\section{Groupes cristallographiques tressés}

\label{sec:crystallographic_braid_groups}

Nous allons donner une idée de la preuve du théorème~\ref{thm:artin_dual_classifiant}, affirmant que le complexe d'intervalle $K_W$ est un espace classifiant pour le groupe d'Artin dual $W_w$. Nous allons pour cela présenter tout d'abord brièvement la construction des groupes cristallographiques tressés de 
\textcite{mccammond_sulway}.

\mk

Considérons un groupe de Coxeter affine irréductible $W$, agissant comme groupe de réflexions sur $\R^n$, où $n$ est le rang de $W$. Soit $w \in W$ un élément de Coxeter, $R$ l'ensemble des réflexions de $W$.

\mk

Notons $R_{hor}$ l'ensemble des réflexions horizontales de~$R$ (i.e. fixant l'axe de Coxeter~$\ell$), et $R_{ver}$ l'ensemble des réflexions verticales. Notons également $T$ l'ensemble fini des translations de $[1,w]^W$. Chaque élément de~$R$ a un poids de~$1$, et chaque élément de~$T$ a un poids de~$2$.

\mk

Notons $k \in \{1,2,3\}$ le nombre de composantes irréductibles du système de racines horizontales $\Phi_{hor}$, voir la partie~\ref{sec:factorisation_euclidean_isometries}. Considérons la décomposition orthogonale de $V$ en $\R\mu \oplus V_1 \oplus \dots \oplus V_k$, où $\mu$ est la direction de la droite de Coxeter $\ell$, et où $V_1,\dots,V_k$ correspondent à la décomposition $\Phi_{hor}=\Phi_1 \sqcup \dots \sqcup \Phi_k$. Pour chaque translation $t \in T$ et pour chaque $1 \leq i \leq k$, considérons la translation $t_i$ de $\R\mu \oplus V_i$ dont les projections vérifient $p_{\R\mu}(t_i) = \f{1}{k}p_{\R\mu}(t)$ et $p_{V_i}(t_i) = p_{V_i}(t)$. On a ainsi $t_1t_2 \dots t_k=t$. L'ensemble de ces translations est noté $T_F$, elles sont appelées \emph{translations de facteurs}, et elles ont un poids de $\f{2}{k}$.

\mk

Rappelons que $W$ est engendré par $R$, et qu'il contient $T$. Nous allons définir trois nouveaux groupes engendrés par certaines réflexions et translations :
\bit
\item Le \emph{groupe diagonal} $D$, engendré par $R_{hor}$ et $T$.
\item Le \emph{groupe factorisé} $F$, engendré par $R_{hor}$ et $T_F$.
\item Le \emph{groupe cristallographique} $C$, engendré par $R$ et $T_F$.
\eit
Remarquons que, lorsque $k=1$, nous avons $D=F$ et $W=C$.

\bex
Lorsque $W$ est de type $\tilde{B}_n$, nous avons vu dans l'exemple~\ref{ex:B3tilde} que le système de racines horizontales se décomposait en deux systèmes de types $\tilde{A}_{n-2}$ et $\tilde{A}_1$, correspondant à la décomposition $\R^n = \R^{n-1} \times \R$. Ainsi pour chaque translation $t \in W$, on ajoutera dans $T_F$ les deux translations correspondant aux composantes de $t$ dans la décomposition $\R^n = \R^{n-1} \times \R$.
\eex

L'intérêt de ces nouveaux groupes est double. Tout d'abord, l'introduction des translations de facteurs permet de rétablir la propriété de treillis qui manquait lorsque~\mbox{$k >1$:}

\bthm[{\cite[Theorem~A]{mccammond_sulway}}]
Les intervalles $[1,w]^F$ et $[1,w]^C$ sont des treillis. En particulier, les groupes d'intervalles associés $F_w$ et $C_w$ sont des groupes de Garside.
\ethm

Le groupe d'intervalle $C_w$ est appelé \emph{groupe cristallographique tressé}. Les quatre groupes sont étroitement reliés de la manière suivante.

\bthm[{\cite[Theorem~9.6]{mccammond_sulway}}]
Les intervalles entre $1$ et $w$ dans les quatre groupes $D,F,W,C$ sont reliés ainsi :
\beq [1,w]^C &=& [1,w]^W \cup [1,w]^F \\
\, [1,w]^D &=& [1,w]^W \cap [1,w]^F.\eeq
De plus, les groupes d'intervalles associés $D_w,F_w,W_w,C_w$ sont tels que $C_w$ est un produit amalgamé :
$$ C_w = W_w \underset{D_w}{\star} F_w.$$
\ethm

Notons $K_D,K_F,K_W,K_C$ les complexes associés respectivement aux intervalles \[[1,w]^D,[1,w]^F,[1,w]^F,[1,w]^C.\] Comme les intervalles $[1,w]^F$ et $[1,w]^C$ sont des treillis, nous savons d'après le théorème~\ref{thm:treillis_classifiant} que $K_F$ et $K_W$ sont des espaces classifiants pour $F_w$ et $C_w$. Nous souhaitons démontrer que $K_W$ est un espace classifiant pour $W_w$, et il s'avère plus simple de montrer d'abord que $K_D$ est un espace classifiant pour $D_w$, car ce groupe lui-même est plus simple.

\bthm[{\cite[Theorem~6.5]{PS}}]\label{thm:KD_classifiant}
Le groupe $D_w$ est une extension par $\Z$ d'un produit de groupes d'Artin de types~$\tilde{A_{n_i}}$, pour $1 \leq i \leq k$. Le complexe~$K_D$ est un espace classifiant pour~$D_w$. 
\ethm

\bp
Considérons le sous-groupe $H$ de $D$ engendré par $R_{hor}$. D'après les notations de la partie~\ref{sec:factorisation_euclidean_isometries}, le système de racines horizontal se décompose en systèmes irréductibles $\Phi_{hor} = \Phi_1 \sqcup \cdots \sqcup \Phi_k$, où chaque $\Phi_i$ est de type $A_{n_i}$. Ainsi $H$ est isomorphe au produit $H=W_1 \times \cdots \times W_k$, où chaque $W_i$ est un groupe de Coxeter affine de type~$\tilde{A}_{n_i}$. 

\mk

On a une décomposition de l'intervalle
$$[1,w]^W \cap H = ([1,w]^W \cap W_1) \times \cdots \times ([1,w]^W \cap W_k),$$
d'après \cite[Proposition~7.6]{mccammond_sulway}, et de plus le groupe $H_w$ associé à l'intervalle $[1,w]^W \cap H$ est un sous-groupe de $D_w$ d'après \cite[Lemma~9.3]{mccammond_sulway}. Et le groupe $H_w$ se décompose ainsi en produit de groupes d'Artin affines de type $\tilde{A}_{n_1},\dots,\tilde{A}_{n_k}$.

\mk 

Notons $K_H$ le sous-complexe de $K_D$ constitué des simplexes $\sigma=[x_1| \dots |x_d]$ tels que $\pi(\sigma)=x_1 \cdots x_d \in H$, on sait que le groupe fondamental de $K_H$ s'identifie à $H_w$.

\mk

Considérons l'automorphisme $\phi\colon g \in W \mapsto w^{-1}gw \in W$ de conjugaison par~$w$. Comme $\phi$ stabilise~$H$, $\phi$~agit par automorphisme sur~$K_H$, et on peut considérer la suspension~$Z$ de~$K_H$ par~$\phi$. Plus précisément, considérons le quotient $Z=K_H \times [0,1]/\sim$, où le simplexe $([x_1| \dots |x_d],1)$ est identifié avec le simplexe $([\phi(x_1)| \dots | \phi(x_d)],0)$. Nous allons voir que $Z$~est homéomorphe à~$K_D$.

\mk

Nous allons tout d'abord définir une structure simpliciale sur $Z$. Fixons un simplexe $\sigma=[x_1| \dots |x_d]$ de $K_H$. La cellule $\sigma \times [0,1]$ de $Z$ est découpée en $d+1$ simplexes $\tau_0,\dots,\tau_d$, où
$$\tau_i = \{(a_1,a_2,\dots,a_d,t) \in \sigma \times [0,1] \st  a_1 \geq \cdots \geq a_i \geq 1-t \geq a_{i+1} \geq \cdots \geq a_d \}.$$
Soit $y \in T$ le complément à droite de $\pi(\sigma)$, i.e. tel que $x_1x_2 \cdots x_dy = w$. Alors le simplexe~$\tau_i$ s'identifie au simplexe $[x_{i+1}| \dots |x_d|y|\phi(x_1)| \dots |\phi(x_i)]$ de $K_D$. Ainsi $Z$ s'identifie à un sous-complexe simplicial de~$K_D$.

\mk

Réciproquement, toute factorisation maximale de $w$ dans $D$ s'écrit \[w=x'_1 \cdots x'_i y' x'_{i+1} \cdots x'_d,\] où chaque $x'_j$ est une réflexion horizontale, et où $y' \in T$. Ainsi tout simplexe maximal de $K_D$ est dans $Z$, donc $Z=K_D$.

\mk

On conclut que $D_w$~est une extension par~$\Z$ de~$H_w$, et que le revêtement universel de~$K_D$, qui s'identifie à $\tilde{K_H} \times \R$, est contractile.
\ep

Comme les quatre complexes d'intervalle sont reliés par $K_C = K_W \cup K_F$ et $K_D = K_C \cap K_F$, nous pouvons maintenant apporter une preuve du théorème~\ref{thm:artin_dual_classifiant} affirmant que $K_W$ est un espace classifiant pour le groupe d'Artin dual $W_w$.

\mk

\bp[Preuve du théorème~\ref{thm:artin_dual_classifiant}]
Considérons le revêtement universel $\rho \colon \tilde{K_C} \ra K_C$. Nous savons donc que $\tilde{K_C} = \rho^{-1}(K_W) \cup \rho^{-1}(K_F)$ et que $\rho^{-1}(K_D) = \rho^{-1}(K_W) \cap \rho^{-1}(K_F)$. Considérons la suite exacte longue de Mayer-Vietoris (à coefficients entiers) :
$$ \cdots \ra H_i(\rho^{-1}(K_D)) \ra H_i(\rho^{-1}(K_W)) \oplus H_i(\rho^{-1}(K_F)) \ra H_i(\tilde{K_C}) \ra \cdots$$
D'après le théorème~\ref{thm:treillis_classifiant} et le théorème~\ref{thm:KD_classifiant}, nous savons que $\tilde{K_C}$ est contractile, et que chaque composante connexe de $\rho^{-1}(K_D)$ et de $\rho^{-1}(K_F)$ est contractile. Ainsi la composante connexe $\tilde{K_W}$ de $\rho^{-1}(K_W)$ est contractile.
\ep

\section{Ordres lexicographiques axiaux}

\label{sec:lexicographic_orderings}

Nous allons présenter une manière géométrique d'ordonner l'ensemble des réflexions d'un groupe de Coxeter inférieures à un élément de Coxeter donné. Ceci permettra dans la suite de montrer que le complexe d'intervalle $K_W$ se rétracte sur le complexe de Salvetti dual $X'_W$. Par ailleurs, cela implique également que l'ensemble ordonné des partitions non croisées de type affine est lexicographiquement décortiquable.

\subsection{Ordre lexicographique et décortiquabilité}

Nous allons rappeler ici la définition de la décortiquabilité lexicographique \parencite{bjorner_wachs_posets, bjorner_wachs_shellable_I,wachs}.

\mk

Soit $P$ un ensemble ordonné borné, c'est-à-dire ayant un élément minimal et un élément maximal. Si $p,q \in P$, notons $p \lessdot q$ si $p < q$ et il n'y a aucun élément $r \in P$ tel que $p < r < q$. Le diagramme de Hasse de~$P$ est le graphe de sommets~$P$, avec une arête entre~$p$ et~$q$ si $p \lessdot q$ ou $q \lessdot p$. Notons $\mathcal{E}(P)$ l'ensemble des arêtes du graphe de Hasse de~$P$.

\mk

Un \emph{étiquetage des arêtes} de $P$ est une application $\lambda \colon \mathcal{E}(P) \ra \Lambda$ à valeurs dans un ensemble totalement ordonné $\Lambda$. Toute chaîne maximale $c=(x \lessdot z_1 \lessdot z_2 \lessdot \cdots \lessdot z_n \lessdot y)$ entre deux éléments $x \leq y$ de $P$ est ainsi étiquetée par le mot
$$\lambda(c) = \lambda(x,z_1) \lambda(z_1,z_2) \dots \lambda(z_t,y).$$

On dit que la chaîne $c$ est \emph{croissante} si le mot associé $\lambda(c)$ est strictement croissant. De plus, si $x \leq y$ sont deux éléments de $P$, alors les chaînes maximales entre $x$ et $y$ peuvent être comparées lexicographiquement, ainsi que colexicographiquement (en comparant les lettres de droite à gauche).

\bdf
Un \emph{étiquetage lexicographique} de $P$ est un étiquetage tel que, pour tout intervalle fermé $[x,y] \subset P$, il existe une unique chaîne maximale croissante de $x$ à $y$, et cette chaîne précède lexicographiquement toutes les autres chaînes maximales. Un ensemble ordonné borné admettant un tel étiquetage lexicographique est dit \emph{lexicographiquement décortiquable}.
\edf

Si $P$ est un ensemble ordonné, son complexe d'ordre est le complexe simplicial de sommets $P$, et dont les simplexes sont donnés par les chaînes de $P$. L'un des intérêts de la notion d'étiquetage lexicographique réside dans le résultat suivant sur la topologie du complexe d'ordre de $P$ :

\bthm[{\cite[Theorem~3.2.2]{wachs}}]
Soit $P$ un ensemble ordonné borné lexicographiquement décortiquable. Alors le complexe d'ordre de $P$ a le type d'homotopie d'un bouquet de sphères.
\ethm

De plus, le produit de deux ensembles ordonnés décortiquables est lui-même décortiquable. Plus précisément :

\bthm[{\cite[Proposition~10.15]{bjorner_wachs_shellable_II}}] \label{thm:produit_EL_shellable}
Soient $P_1,P_2$ deux ensembles ordonnés admettant des étiquetages lexicographiques $\lambda_i\colon \mathcal{E}(P_i) \ra \Lambda_i$. Considérons un ordre total sur $\Lambda = \Lambda_1 \sqcup \Lambda_2$ se restreignant aux ordres de $\Lambda_1$ et de $\Lambda_2$. Alors $\lambda \colon \mathcal{E}(P_1 \times P_2) \ra \Lambda$ est un étiquetage lexicographique.
\ethm

\subsection{Ordres axiaux dans le cas fini}

Soit $W$ un groupe de Coxeter fini, agissant par isométries linéaires sur $V=\R^n$. Notons $\Phi \subset V$ le système de racines de $W$, $w$ un élément de Coxeter de $W$ et $\Phi^+ \subset \Phi$ le système positif associé. Notons $R$ l'ensemble des réflexions de $W$. Si $\alpha \in \Phi$, notons $r_\alpha \in R$ la réflexion orthogonale par rapport à $\alpha$.

\bdf
Un ordre total $\prec$ sur $R$ est appelé \emph{ordre de réflexion} si, pour toutes racines positives distinctes $\alpha_1,\alpha_2 \in \Phi^+$, et pour toute racine $\alpha \in \Phi^+$ qui est une combinaison linéaire positive de $\alpha_1$ et $\alpha_2$, on a
$$ r_{\alpha_1} \prec r_\alpha \prec r_{\alpha_2} \mbox{ ou } r_{\alpha_2} \prec r_\alpha \prec r_{\alpha_1}.$$
Cet ordre est dit \emph{compatible} avec $w$ si, dès que $\alpha,\beta \in \Phi^+$ sont les racines simples d'un sous-système de racines irréductible de rang $2$ et que $r_\alpha r_\beta \in [1,w]$, alors $r_\alpha \prec r_\beta$.
\edf

L'intérêt de cette notion concerne la décortiquabilité, comme l'ont montré Athanasiadis, Brady et Watt.

\bthm[{\cite[Theorem~3.5]{athanasiadis_brady_watt}}] \label{thm:fini_ordre_axial_EL}
Soit $W$ un groupe de Coxeter fini cristallographique, $w$ un élément de Coxeter, et $R$ l'ensemble des réflexions. Si $\prec$ est un ordre de réflexion sur $R$ compatible avec $w$, alors l'étiquetage associé de $\mathcal{E}([1,w])$ est un étiquetage lexicographique.
\ethm

Nous allons maintenant décrire une méthode géométrique simple permettant de construire de tels ordres de réflexions compatibles.

Considérons l'arrangement d'hyperplans $\mathcal{A}$ associé à $W$. Soit $C_0$ la chambre du complexe de Coxeter associé à $\Phi^+$. Considérons une droite affine $\ell' = \{a+\theta \mu \st \theta \in \R\} \subset V$ qui soit générique par rapport à $\mathcal{A}$, où $a \in C_0$ est un point base de $\ell'$, et $\mu \in V \smallsetminus \{0\}$ oriente la droite $\ell'$. On dit qu'un point $b \in \ell'$ est \emph{au-dessus} d'un point $b' \in \ell'$ si $b-b'$ est un multiple positif de $\mu$, et \emph{en-dessous} sinon.

Ceci permet de définir un ordre total sur $R$ :
\bit
\item en premier viennent les réflexions fixant un point de $\ell'$ situé au-dessus de $a$, et $r$ vient avant $r'$ si $\Fix(r) \cap \ell'$ est en-dessous de $\Fix(r') \cap \ell'$,
\item ensuite viennent les réflexions fixant un point de $\ell'$ situé en-dessous de $a$, et $r$ vient avant $r'$ si $\Fix(r) \cap \ell'$ est en-dessous de $\Fix(r') \cap \ell'$.
\eit

\bpro \label{pro:droite_generique_ordre_reflection}
Pour toute telle droite générique $\ell'$, l'ordre associé est un ordre de réflexion.
\epro

\bp
Comme $a \in C_0$, nous savons que pour toute racine positive $\alpha \in \Phi^+$ nous avons $\<a,\alpha\> >0$. Supposons pour simplifier que nous avons renormalisé les racines positives de sorte que $\<a,\alpha\> = 1$.

Pour toute racine positive $\alpha \in \Phi^+$, l'intersection entre l'hyperplan $\Fix(r_\alpha)$ et la droite~$\ell'$ est
$$\Fix(r_\alpha) \cap \ell' = \left\{a -\f{1}{\<\mu,\alpha\>}\mu\right\}.$$
Par définition de $\prec$, on a donc que $r_\alpha \prec r_\beta$ si et seulement si $\<\mu,\alpha\> < \<\mu,\beta\>$.

Ainsi, si $\alpha=c_1\alpha_1+c_2\alpha_2$ est une combinaison positive de racines positives, alors $\<\mu,\alpha\>$ est compris entre $\<\mu,\alpha_1\>$ et $\<\mu,\alpha_2\>$, donc  $r_{\alpha_1} \prec r_\alpha \prec r_{\alpha_2}$ ou $r_{\alpha_2} \prec r_\alpha \prec r_{\alpha_1}$.
\ep

\bex \label{ex:ordre_reflection}
Nous allons maintenant décrire un tel exemple d'ordre de réflexion présenté dans~\parencite[Example~3.3]{athanasiadis_brady_watt}. Considérons le groupe de Coxeter $W \simeq \frak{S}_{n}$ de type $A_{n-1}$. Choisissons comme élément de Coxeter $w \in W$ le $n$-cycle $w=(1,2, \dots ,n)=(1,2)(2,3) \cdots (n-1,n)$, où $(i,j)$ désigne la transposition permutant $i$ et $j$. L'intervalle $[1,w]$ s'identifie au treillis des partitions non croisées de $n$ points (voir l'exemple~\ref{ex:tresses_garside}). L'arrangement d'hyperplans associé dans $\R^n$ est $\{x_i-x_j=0, 1 \leq i < j \leq n\}$. Fixons $a \in \R^n$ tel que $a_1<a_2< \cdots <a_n$, et $\eps>0$ suffisamment petit. Considérons la droite
$$\ell'=\{a+\theta(1,\eps,\eps^2,\dots,\eps^{n-1}), \theta \in \R\}.$$
Alors $\ell'$ intersecte chaque hyperplan $\{x_i-x_j=0\}$ avec $\theta>0$, et l'ordre de réflexion associé est le suivant : la transposition $(i,j)$ vient avant la transposition $(i',j')$ si et seulement si $i < i'$, ou $i=i'$ et $j \leq j'$. Ceci définit un ordre de réflexion sur $W$.

De plus, cet ordre est compatible avec l'élément de Coxeter $w$: un sous-système $\Phi' \subset \Phi$ irréductible de rang $2$ correspond au choix d'indices $1 \leq i < j < k \leq n$, dont les racines simples sont $(i,j)$ et $(j,k)$. Alors $(i,j)(j,k) \in [1,w]^W$, et on a bien $(i,j) \prec (j,k)$.

Ainsi, d'après le théorème~\ref{thm:fini_ordre_axial_EL}, l'étiquetage associé de $\mathcal{E}([1,w])$ est un étiquetage lexicographique.
\eex

\subsection{Ordres sur les réflexions horizontales} \label{subsec:ordre_horizontal}

Nous allons maintenant décrire comment choisir une telle droite $\ell'$ afin de définir un ordre de réflexion sur l'ensemble des réflexions horizontales.

Notons $\Phi_{hor} = \Phi_1 \sqcup \cdots \sqcup \Phi_k$ la décomposition du système de racines horizontal en $k$ composantes irréductibles, où $\Phi_i$ est de type $A_{n_i}$ (voir la partie~\ref{sec:factorisation_euclidean_isometries}). Choisissons une factorisation horizontale $w=th_1 \cdots h_k$ de $w$, où $t$ est une translation, et $h_i$ est un élément de Coxeter pour le sous-groupe parabolique sphérique $W_{h_i} \subset W$ associé à l'intervalle $[1,h_i]^W$.


\blem[{\cite[Lemma~4.5]{PS}}] \label{lem:droite_elli_horizontal}
Pour tout point $a \in \ell$ sur l'axe de Coxeter, il existe une droite $\ell_i$ contenant $a$ de direction le sous-espace engendré par $\Phi_i$, telle que l'ordre de réflexion $\prec_{\ell_i}$ sur $W_{h_i}$ associé à $\ell_i$ soit compatible avec $h_i$.
\elem

\bp
Le groupe de Coxeter $W_{h_i}$ est de type $A_{n_i}$, on peut donc considérer une droite $\ell'$ comme dans l'exemple~\ref{ex:ordre_reflection} puis considérer sa projection sur le sous-espace affine de direction le sous-espace engendré par $\Phi_i$.
\ep

Notons $W_i \subset W$ le sous-groupe de Coxeter de type $\tilde{A}_{n_i}$ engendré par les réflexions par rapport aux racines de $\Phi_i$.

Nous étendons cet ordre de réflexion à un ordre $\prec_i$ sur $R_{hor} \cap W_i$ de la manière suivante : si $r_1 \prec_{\ell_i} r_2$, alors toute réflexion parallèle à $r_1$ vient avant toute réflexion parallèle à $r_2$. Remarquons que deux réflexions parallèles ne peuvent intervenir dans une factorisation minimale de $w$, ainsi l'ordre entre deux réflexions parallèles peut être choisi arbitrairement.

\mk

Considérons un ordre total $\prec_{hor}$ sur $R_{hor}$ tel que, pour tout $1 \leq i \leq k$, la restriction de $\prec_{hor}$ à $R_{hor} \cap W_i$ soit égale à $\prec_i$. Ceci nous permet de définir un étiquetage sur l'ensemble des isométries horizontales.

\blem[{\cite[Lemma~4.9]{PS}}]
Soit $W$ un groupe de Coxeter affine irréductible, et soit $w$ un élément de Coxeter. Pour tout élément horizontal $u \in [1,w]^W$, l'ordre $\prec_{hor}$ définit un étiquetage lexicographique $\lambda \colon \mathcal{E}([1,u]) \ra R_{hor}$.
\elem

\bp
Il suffit d'appliquer le théorème~\ref{thm:produit_EL_shellable} garantissant qu'un produit d'ensembles ordonnés décortiquables est décortiquable.
\ep

\subsection{Ordres axiaux dans le cas affine}

Soit $W$ un groupe de Coxeter affine irréductible, $w$ un élément de Coxeter, et $R_0$ l'ensemble des réflexions appartenant à $[1,w]^W$. Notons $\ell \subset V$ l'axe de Coxeter, et $C_0$ une chambre axiale du complexe de Coxeter.

\bdf
Un \emph{ordre axial} sur $R_0$ est un ordre total comme suit :
\bit
\item en premier viennent les réflexions verticales fixant un point de $\ell$ situé au-dessus de $C_0$, et $r$ vient avant $r'$ si $\Fix(r) \cap \ell$ est en-dessous de $\Fix(r') \cap \ell$ : ces réflexions verticales sont appelées \emph{positives};
\item ensuite viennent les réflexions horizontales de $R_{hor}$, dans l'un des ordres totaux $\prec_{hor}$ construits dans la partie~\ref{subsec:ordre_horizontal};
\item enfin viennent les réflexions verticales fixant un point de $\ell$ situé en-dessous de $C_0$, et $r$ vient avant $r'$ si $\Fix(r) \cap \ell$ est en-dessous de $\Fix(r') \cap \ell$ : ces réflexions verticales sont appelées \emph{négatives}.
\eit
Si deux réflexions verticales fixent le même point de $\ell$, leur ordre est choisi arbitrairement.
\edf

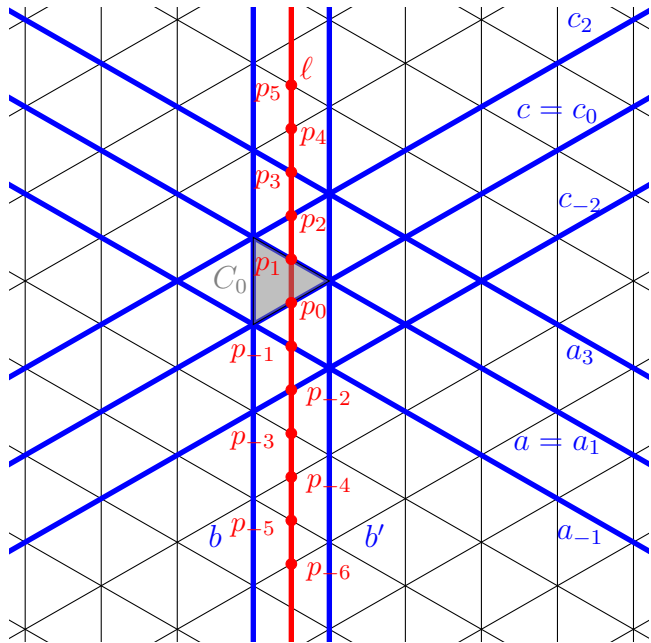
\begin{figure}
\begin{center}
\begin{tikzpicture}
\def \p {0.07}
\def \op {0.5}
\def \gris {black!50}
\clip (-3.2,-4.2) rectangle (5.2,4.2);

\foreach \i in {-10,...,10}
\draw (\i,-10) -- (\i,10);

\foreach \i in {-10,...,10}
\draw (-17.32-\i/2,-10+\i*0.866) -- (17.32-\i/2,10+\i*0.866);

\foreach \i in {-10,...,10}
\draw (-17.32-\i/2,10-\i*0.866) -- (17.32-\i/2,-10-\i*0.866);

\draw[line width=2, color=blue] (0,-10) -- (0,10);
\draw[line width=2, color=blue] (1,-10) -- (1,10);
\foreach \i in {-1,...,1}
\draw[line width=2, color=blue] (-17.32-\i/2,-10+\i*0.866) -- (17.32-\i/2,10+\i*0.866);
\foreach \i in {-2,...,0}
\draw[line width=2, color=blue] (-17.32-\i/2,10-\i*0.866) -- (17.32-\i/2,-10-\i*0.866);

\node[thick, color=blue] at (-0.5,-2.8) {$b$};
\node[thick, color=blue] at (1.6,-2.8) {$b'$};
\node[thick, color=blue] at (4.3,-0.4) {$a_3$};
\node[thick, color=blue] at (4,-1.6) {$a=a_1$};
\node[thick, color=blue] at (4.3,-2.8) {$a_{-1}$};
\node[thick, color=blue] at (4.3,1.6) {$c_{-2}$};
\node[thick, color=blue] at (4,2.8) {$c=c_0$};
\node[thick, color=blue] at (4.3,4.0) {$c_2$};

\draw[line width=2, color=red] (0.5,-10) -- (0.5,10);
\node[thick, color=red] at (0.7,3.4) {$\ell$};

\draw[black,fill opacity=\op,fill=\gris] (0,0) -- (1,1/1.732) -- (0,2/1.732) -- cycle;
\node[thick, color=\gris] at (-0.3,0.6) {$C_0$};

\foreach \i in {0,2,4}
{\draw[fill, color=red] (1/2,0.5/1.732+\i/1.732) circle (\p) node {};
\node[thick, color=red] at (1/2+0.3,0.5/1.732+\i/1.732-0.1) {\small $p_\i$};}

\foreach \i in {2,4,6}
{\draw[fill, color=red] (1/2,0.5/1.732-\i/1.732) circle (\p) node {};
\node[thick, color=red] at (1/2+0.5,0.5/1.732-\i/1.732-0.1) {\small $p_{-\i}$};}

\foreach \i in {1,3,5}
{\draw[fill, color=red] (1/2,0.5/1.732+\i/1.732) circle (\p) node {};
\node[thick, color=red] at (1/2-0.3,0.5/1.732+\i/1.732-0.1) {\small $p_\i$};}

\foreach \i in {1,3,5}
{\draw[fill, color=red] (1/2,0.5/1.732-\i/1.732) circle (\p) node {};
\node[thick, color=red] at (1/2-0.5,0.5/1.732-\i/1.732-0.1) {\small $p_{-\i}$};}

\end{tikzpicture}
\caption{Un ordre axial sur les réflexions en type $\tilde{A}_2$.}
\label{fig:A2tilde_ordre_axial}
\end{center}
\end{figure}

\bex \label{ex:ordre_axial_A2tilde}
Considérons le groupe de Coxeter $W$ de type $\tilde{A}_2$. Notons $a,b,c$ les générateurs standards de $W$, et $w=abc$ un élément de Coxeter. Notons $\ell \subset \R^2$ l'axe de Coxeter, comme sur la figure~\ref{fig:A2tilde_ordre_axial}. Les intersections de l'axe de Coxeter $\ell$ avec les hyperplans de réflexions forment une suite de points $(p_n)_{n \in \Z}$, avec $p_0$ et $p_1$ au bord de la chambre $C_0$.

D'après le théorème~\ref{thm:reflection_sommets_axiaux}, les réflexions horizontales de $[1,w]^W$ sont celles qui fixent un sommet axial, il s'agit donc de $b$ et $b'$. Pour tout $i \in \Z$ pair, notons $c_i \in R$ la réflexion parallèle à $c=c_0$ fixant $p_i$. Pour tout $i \in \Z$ impair, notons $a_i \in R$ la réflexion parallèle à $a=a_1$ fixant $p_i$. L'ensemble des réflexions verticales de $[1,w]^W$ est $\{c_i \st i \in \Z \mbox{ pair }\} \cup \{a_i \st i \in \Z \mbox{ impair }\}$.

Un ordre axial sur $R_0$ est donc donné par
$$ \underbrace{a_1 \prec c_2 \prec a_3 \prec c_4 \prec \cdots }_{\mbox{réflexions positives}} \prec \underbrace{b \prec b'}_{\mbox{réflexions horizontales}} \prec \underbrace{\cdots \prec a_{-3} \prec c_{-2} \prec a_{-1} \prec c_0}_{\mbox{réflexions négatives}}.$$
\eex

Cet ordre axial permet de montrer l'existence et l'unicité de chaînes lexicographiquement minimales dans les intervalles.

\blem[{\cite[Lemma~4.16]{PS}}] \label{lem:unique_chaine}
Considérons un ordre axial $\prec$ sur $R_0$. Tout intervalle $[u,v]$ contenu dans $[1,w]$ a une unique chaîne maximale lexicographiquement minimale, et celle-ci est croissante. De même, $[u,v]$ a une unique chaîne maximale colexicographiquement maximale, et celle-ci est croissante.
\elem

\bp
Nous allons le montrer par récurrence sur la longueur de l'intervalle $[u,v]$. Considérons les réflexions apparaissant dans cet intervalle. Ce sont les éléments de $R_0 \cap [1,u^{-1}v]$.

\mk

Si $u^{-1}v$ est elliptique, alors $R_0 \cap [1,u^{-1}v]$ est fini, et a donc une unique réflexion minimale $r$ pour $\prec$.

Si $u^{-1}v$ est hyperbolique, alors $R_0 \cap [1,u^{-1}v]$ contient au moins une réflexion verticale positive, et a donc une unique réflexion minimale $r$ pour $\prec$.

\mk

\'{E}crivons $u'=ur$, alors par hypothèse de récurrence sur l'intervalle $[u',v]$, il suffit de montrer que toutes les réflexions de $R_0 \cap [1,u'^{-1}v]$ sont supérieures à $r$. Si $r' \in R_0 \cap [1,u'^{-1}v]$, alors il existe une factorisation de $u^{-1}v$ débutant par $rr'$, donc $r \prec r'$.
\ep

Ainsi la preuve de la décortiquabilité se ramène à montrer qu'il y a au plus une chaîne maximale croissante dans un intervalle $[1,u]^W$. Nous allons distinguer les cas où $u$ est elliptique ou hyperbolique. Les preuves de ces deux lemmes étant un peu techniques, nous en donnerons seulement les grandes lignes.

\blem[{\cite[Lemma~4.17]{PS}}] \label{lem:unique_chaine_maximale_elliptique}
Considérons un ordre axial $\prec$ sur $R_0$, et soit $u \in [1,w]$ un élément elliptique. L'intervalle $[1,u]$ a au plus une chaîne maximale croissante.
\elem

\bp
Pour tout $1 \leq i \leq k$, on peut trouver un élément maximal $u_i$ de $[1,w] \cap W_i$ tel que $[1,u] \cap W_i \subset [1,u_i]$. Pour simplifier, on peut supposer que $u_i=h_i$.

Fixons un point $a \in C_0 \cap \ell$. Le lemme~\ref{lem:droite_elli_horizontal} donne une droite~$\ell_i$ contenant~$a$, dirigée par un vecteur~$\mu_i$ dans la direction de~$\Phi_i$. Notons~$\mu$ le vecteur orientant la droite de Coxeter~$\ell$. Pour $\eps>0$ suffisamment petit, considérons la droite~$\ell'$ passant par~$a$, dirigée par le vecteur
$$\mu' = \mu + \eps\mu_1 + \cdots + \eps \mu_k.$$
Les droites $\ell'$ et $\ell'_i$ intersectent les hyperplans des réflexions de $W_{h_i}$ dans le même ordre.

\mk

Perturbons légèrement la droite $\ell'$, de sorte qu'elle devienne générique par rapport aux hyperplans des réflexions de $W_u$. Alors les ordres $\prec_{\ell'}$ et $\prec$ ne diffèrent éventuellement que pour des paires de réflexions fixant le même point de $\ell$. Pour simplifier, supposons que ces ordres sont les mêmes.

\mk

Comme $\ell'$~est générique, d'après la proposition~\ref{pro:droite_generique_ordre_reflection}, l'ordre $\prec_{\ell'}$ est un ordre de réflexion. Nous admettons que cet ordre est compatible avec l'élément de Coxeter~$u$ de~$W_u$. Ainsi, d'après le théorème~\ref{thm:fini_ordre_axial_EL}, l'étiquetage associé de $\mathcal{E}([1,u]^W)$ est un étiquetage lexicographique.
\ep

\blem[{\cite[Lemma~4.18]{PS}}] \label{lem:unique_chaine_maximale_hyperbolique}
Considérons un ordre axial $\prec$ sur $R_0$, et soit $u \in [1,w]$ un élément hyperbolique tel que le sous-groupe $W_u$ soit irréductible. L'intervalle $[1,u]$ a au plus une chaîne maximale croissante.
\elem

\bp
Supposons que nous ayons une chaîne maximale croissante de $[1,u]^W$, correspondant à une factorisation $u=r_1r_2 \cdots r_m$. Comme $u$~est une isométrie verticale, $r_1$ ou~$r_m$ est une réflexion verticale.

\mk

Si $r_1$ est une réflexion verticale qui n'est pas minimale parmi $R_0 \cap [1,u]$, alors il existe une factorisation $u=r_1r'_2 \cdots r'_m$ telle que $r'_2 \prec r_1$. Notons $u=r_1u'$ : d'après le lemme~\ref{lem:unique_chaine_maximale_elliptique} on déduit aussi que l'intervalle $[1,u']$ a au plus une chaîne maximale croissante. Ainsi $u'=r_2 \cdots r_m$ est l'unique factorisation lexicographiquement minimale, donc $r_2 \preceq r'_2 \prec r_1$ : ceci contredit $r_1 \prec r_2$.

\mk

Donc si $r_1$ est une réflexion verticale, elle est minimale parmi $R_0 \cap [1,u]$, et $u'=r_2 \cdots r_m$ est elliptique. Ainsi, d'après le premier cas, on déduit que $u'=r_2 \cdots r_m$ est l'unique factorisation lexicographiquement minimale de $u'$.

\mk

De même, si $r_m$ est une réflexion verticale, alors elle est maximale parmi $R_0 \cap [1,u]$, et $r_1 \cdots r_{m-1}$ est l'unique factorisation lexicographiquement minimale.

\mk

Il y a donc au plus deux chaînes maximales croissantes. Nous admettons qu'il s'agit de la même chaîne.
\ep

Nous avons présenté les ingrédients nécessaires à la preuve de la décortiquabilité de l'ensemble ordonné $[1,w]^W$.

\bthm[{\cite[Theorem~4.19]{PS}}] \label{thm:ordre_axial_EL}
Soit $W$ un groupe de Coxeter affine irréductible, $w$ un élément de Coxeter, et $R_0$ l'ensemble des réflexions appartenant à $[1,w]^W$. Considérons un ordre axial sur $R_0$. Considérons l'étiquetage des arêtes $\lambda : \mathcal{E}([1,w]^W) \ra R_0$ naturel. Alors tout intervalle $[u,v]^W$ de $[1,w]^W$ possède une unique chaîne maximale strictement croissante, et cette chaîne est à la fois lexicographiquement minimale et colexigraphiquement maximale. En particulier, $\lambda$ est un étiquetage lexicographique.
\ethm

Une conséquence importante est que l'ensemble ordonné des partitions non croisées affines $[1,w]^W$ est lexicographiquement décortiquable.

\mk

\bp
Remarquons que la multiplication à gauche par $u^{-1}$ est un isomorphisme de $[u,v]^W$ sur $[1,u^{-1}v]^W$ : on peut donc se contenter d'étudier l'intervalle $[1,u]^W$. D'après le lemme~\ref{lem:unique_chaine}, il suffit de montrer que $[1,u]^W$ a au plus une chaîne maximale croissante. Si $u$ est elliptique, c'est l'objet du lemme~\ref{lem:unique_chaine_maximale_elliptique}. Si $u$ est hyperbolique et que $W_u$ est irréductible, c'est l'objet du lemme~\ref{lem:unique_chaine_maximale_hyperbolique}.

\mk

Supposons maintenant que $u$ est hyperbolique quelconque. D'après le théorème~\ref{thm:decomposition_horizontale_hyperbolique}, considérons la décomposition hyperbolique-horizontale $u=u'h$ de l'isométrie $u$. L'isométrie $u'$ est hyperbolique et le groupe de Coxeter $W_{u'}$ est irréductible, donc $\lambda$ est un étiquetage lexicographique de $[1,u']$. L'isométrie $h$ est elliptique horizontale, donc $\lambda$ est un étiquetage lexicographique de $[1,h]$.
Comme on a la décomposition $[1,u]=[1,u'] \times [1,h]$, d'après le théorème~\ref{thm:produit_EL_shellable}, $\lambda$ est un étiquetage lexicographique de $[1,u]$. En particulier, l'intervalle $[1,u]$ a une unique chaîne croissante maximale.
\ep

\section{Espaces classifiants finis pour les groupes d'Artin duaux}

\label{sec:classifying_space_dual_artin}

Nous avons maintenant les outils nécessaires à la description d'un sous-complexe fini~$K'_W$ de~$K_W$ sur lequel $K_W$~se rétracte. Cette construction va ainsi fournir un espace classifiant fini pour le groupe d'Artin dual~$W_w$.

\mk

L'ensemble des faces de $K_W$ possède une structure très particulière provenant de l'action de l'élément de Coxeter $w$ par conjugaison. Notons $\mathcal{F}(K_W)$ l'ensemble ordonné des faces de $K_W$.

\mk

Si $\sigma=[x_1|x_2| \dots |x_d] \in \mathcal{F}(K_W)$ est un $d$-simplexe de $K_W$, on notera $\pi(\sigma)=x_1x_2 \cdots x_d \in [1,w]^W$. Partant de $\sigma$, il y a une manière naturelle de se déplacer \og à droite\fg{} pour aller en $\rho(\sigma)$ (dans la direction de $w$) ou \og à gauche\fg{} pour aller en $\lambda(\sigma)$ (dans la direction de $w^{-1}$) parmi les simplexes de $K_W$ :
$$ \rho(\sigma) = \left\{ \begin{tabular}{l} $[x_2|\dots |x_d]$ si $\pi(\sigma)=w$ \\ $[x_1|\dots |x_d|y]$ sinon, où $x_1 \cdots x_dy = w$ \end{tabular}\right.$$
$$ \lambda(\sigma) = \left\{ \begin{tabular}{l} $[x_1|\dots |x_{d-1}]$ si $\pi(\sigma)=w$ \\ $[y|x_1|\dots |x_d]$ sinon, où $yx_1 \cdots x_d = w$ \end{tabular}\right.$$
Ceci est permis par le fait que le treillis $[1,w]^W$ est équilibré. Remarquons que $\rho$ est bien l'inverse de $\lambda$.

\mk

Nous allons appeler \emph{composantes fibrées} de $\mathcal{F}(K_W)$ les orbites sous l'action de $\lambda$ et $\rho$. Plus formellement, on peut considérer l'application entre ensembles ordonnés
\beq \eta \colon{} \mathcal{F}(K_W) & \ra & \N \\
\, [x_1|x_2| \dots |x_d] & \mapsto & \left\{ \begin{tabular}{l} $d$ si $\pi(\sigma)=x_1x_2 \cdots x_d=w$ \\ $d+1$ sinon \end{tabular}\right.
\eeq
Les composantes fibrées de $\mathcal{F}(K_W)$ sont alors les composantes connexes des images réciproques $\eta^{-1}(d)$, pour $d \geq 1$, dans le diagramme de Hasse de $\mathcal{F}(K_W)$. Le cas du type~$\tilde{A}_2$ est présenté dans l'exemple~\ref{ex:fiber_components_A2tilde}.

\mk

Si $\pi(\sigma)=w$, nous dirons que $\sigma$ est un simplexe \emph{supérieur}, et \emph{inférieur} sinon. Notons que chaque composante fibrée alterne entre simplexes supérieurs et inférieurs.

\mk

La proposition~\ref{pro:factorisations_elliptiques} permet de montrer la description simple suivante des simplexes de~$K_W$ :

\blem \label{lem:types de simplexes}
Soit $\sigma=[x_1|x_2| \dots |x_d]$ un $d$-simplexe supérieur de $K_W$, avec $d \geq 1$. Alors $\sigma$ est exactement de l'un de ces deux types :
\ben
\item soit chaque $x_i$ est elliptique, et au moins l'un d'entre eux est vertical,
\item soit chaque $x_i$ est elliptique horizontal ou hyperbolique.
\een
\elem

Nous allons étudier la topologie de $K_W$ par l'étude des composantes fibrées.

\bpro \label{pro:composantes_fibrees}
L'ensemble $\mathcal{F}(K_W)$ possède un nombre fini de composantes fibrées, qui sont de deux types :
\bit
\item Les composantes infinies, dont les simplexes sont du type 1. du lemme~\ref{lem:types de simplexes}.
\item Les composantes finies, dont les simplexes sont du type 2. du lemme~\ref{lem:types de simplexes}.
\eit
De plus, toute composante fibrée intersecte $\mathcal{F}(X'_W) \subset \mathcal{F}(K_W)$.
\epro

Nous allons donner une preuve de cette proposition en nous appuyant sur les deux lemmes suivants.

\blem
Soit $\mathcal{C} \subset \mathcal{F}(K_W)$ une composante fibrée finie. Alors il existe $\sigma \in \mathcal{C}$ tel que $\pi(\sigma)$ est une isométrie elliptique horizontale.
\elem

\bp
Soit $\sigma=[x_1|x_2| \dots |x_d] \in \mathcal{C}$ tel que $\pi(\sigma)=w$. Si l'un des $x_i$ est une isométrie elliptique verticale, la composante $\mathcal{C}$ est infinie. Ainsi $\sigma$ est du type 2 du lemme~\ref{lem:types de simplexes} : on peut donc supposer que $x_d$ est hyperbolique. Ainsi $\sigma'=[x_1|x_2| \dots |x_{d-1}] \in \mathcal{C}$, et $\pi(\sigma')=x_1x_2 \cdots x_{d-1}$ est elliptique horizontale.
\ep

\blem
Soit $\mathcal{C} \subset \mathcal{F}(K_W)$ une composante fibrée infinie. Alors il existe $\sigma \in \mathcal{C}$ tel que $\pi(\sigma)$ est une isométrie elliptique verticale.
\elem

\bp
Soit $\sigma=[x_1|x_2| \dots |x_d] \in \mathcal{C}$ tel que $\pi(\sigma)=w$. Comme $[1,w]^W$ n'a qu'un nombre fini d'isométries elliptiques horizontales ou hyperboliques, on déduit que $\sigma$ est du type 1 du lemme~\ref{lem:types de simplexes} : on peut donc supposer que $x_d$ est elliptique verticale. Ainsi $\sigma'=[x_1|x_2| \dots |x_{d-1}] \in \mathcal{C}$, et $\pi(\sigma')=x_1x_2 \cdots x_{d-1}$ est également elliptique verticale.
\ep

\bp[Démonstration de la proposition~\ref{pro:composantes_fibrees}]
Soit $\mathcal{C}$ une composante fibrée de $\mathcal{F}(K_W)$. D'après les lemmes précédents, il existe $\sigma \in \mathcal{C}$ tel que $x=\pi(\sigma)$ soit elliptique : ainsi $x$ fixe un sommet axial. Quitte à conjuguer par une puissance de $w$, on peut donc supposer que $x$ fixe un sommet de $C_0$. En particulier, $\sigma \in \mathcal{F}(X'_W)$ donc $\mathcal{C}$ intersecte $\mathcal{F}(X'_W)$. De plus, comme $X'_W$ est un complexe fini, il n'y a qu'un nombre fini de composantes fibrées.
\ep

Ceci nous permet de définir un sous-complexe intéressant de $K_W$.

\bdf
Le sous-complexe $K'_W$ de $K_W$ a pour simplexes la réunion des composantes fibrées finies de $\mathcal{F}(K_W)$, ainsi que les simplexes des composantes fibrées infinies de $\mathcal{F}(K_W)$ compris entre le premier et le dernier simplexe appartenant à $\mathcal{F}(X'_W)$. Le complexe $K'_W$ est appelé \emph{sous-complexe adapté} de $K_W$.
\edf

\bthm[{\cite[Theorem~7.9]{PS}}] \label{thm:KW_retracte_KprimeW}
Soit $W$ un groupe de Coxeter affine irréductible, et soit $w$ un élément de Coxeter. Le complexe d'intervalle~$K_W$ associé à $[1,w]^W$ possède un sous-complexe adapté fini~$K'_W$ contenant~$X'_W$, tel que $K_W$ se rétracte par déformation sur~$K'_W$.
\ethm

\bp
Nous allons utiliser la théorie de Morse discrète, et renvoyons à la partie~\ref{sec:discrete_morse_theory} pour plus de détails. Nous allons définir un couplage acyclique propre sur $\mathcal{F}(K_W)$ dont l'ensemble des simplexes critiques sera précisément $\mathcal{F}(K'_W)$.

\mk

Pour chaque composante fibrée infinie $\mathcal{C}$ de $\mathcal{F}(K_W)$, remarquons que $\mathcal{C}$ est une droite dont $\mathcal{C} \cap \mathcal{F}(K'_W)$ est un segment non vide. Considérons l'unique couplage acyclique propre $\mathcal{M}_\mathcal{C}$ dont les simplexes critiques sont $\mathcal{C} \cap \mathcal{F}(K'_W)$. D'après le théorème du patchwork (théorème~\ref{thm:patchwork}), la réunion de ces couplages est un couplage $\mathcal{M}$ sur $\mathcal{F}(K_W)$ dont les simplexes critiques sont $\mathcal{F}(K'_W)$. Il est clair que ce couplage est propre. D'après le théorème fondamental de la théorie de Morse discrète (théorème~\ref{thm:discrete_morse}), le complexe $K_W$ se rétracte par déformation sur son sous-complexe adapté $K'_W$.
\ep

Ce résultat s'étend également aux groupes cristallographiques tressés, et a la conséquence immédiate suivante.

\bthm[{\cite[Theorem~7.10]{PS}}]
Les groupes d'Artin affines, ainsi que les groupes cristallographiques tressés, ont un espace classifiant fini.
\ethm

\begin{figure}
	\newcommand{\height}{-1.5}
	\newcommand{\length}{0.95}
	\begin{center}
		\begin{tikzpicture}
		\node[inner sep=5pt] (0) at (0,0) {$[w]$};
		\node[inner sep=5pt] (1) at (0, \height) {$[\,]$};
		\draw (0.south) -> (1.north);
		\node[fill=white, draw, circle, inner sep=1.4pt] at (0.south) {};
		\node[fill=black, draw, circle, inner sep=1.4pt] at (1.north) {};
		\end{tikzpicture}\qquad\quad
		\begin{tikzpicture}[every node/.style={inner sep=5pt}]
		\node (0) at (0,0) {$[c_2c_0|b']$};
		\node (1) at (\length, \height) {$[b']$};
		\node (2) at (2*\length, 0) {$[b'|a_1a_{-1}]$};
		\node (3) at (3*\length, \height) {$[a_1a_{-1}]$};
		\node (4) at (4*\length, 0) {$[a_1a_{-1} | b]$};
		
		\node (-1) at (-\length, \height) {$[c_2c_0]$};
		\node (-2) at (-2*\length, 0) {$[b | c_2c_0]$};
		\node (-3) at (-3*\length, \height) {$[b]$};

		\draw[name path=a, white] (0.south) -> (1.north);
		\draw[name path=b, white] (2.south) -> (1.north);
		\draw[name path=c, white] (2.south) -> (3.north);
		
		\draw[name path=d, white] (0.south) -> (-1.north);
		\draw[name path=e, white] (-2.south) -> (-1.north);
	
		\draw[name path=z] (4.south) -> (-3.north);
		\begin{scope}[every node/.style={inner sep=2pt, fill=white, circle}]
			\path [name intersections={of=a and z,by=A}];
			\node at (A) {};
			\path [name intersections={of=b and z,by=B}];
			\node at (B) {};
			\path [name intersections={of=c and z,by=C}];
			\node at (C) {};
			\path [name intersections={of=d and z,by=D}];
			\node at (D) {};
			\path [name intersections={of=e and z,by=E}];
			\node at (E) {};
		\end{scope}

		\draw (0.south) -> (1.north);
		\draw (2.south) -> (1.north);
		\draw (2.south) -> (3.north);
		\draw (4.south) -> (3.north);
		
		\draw (0.south) -> (-1.north);
		\draw (-2.south) -> (-1.north);
		\draw (-2.south) -> (-3.north);
		
		\begin{scope}[every node/.style={fill=white, draw, circle, inner sep=1.4pt}]
			\node at (0.south) {};
			\node at (2.south) {};
			\node at (4.south) {};
			\node at (-2.south) {};
			\node at (-1.north) {};
			\node at (3.north) {};
		\end{scope}
		
		\begin{scope}[every node/.style={fill=black, draw, circle, inner sep=1.4pt}]
			\node at (-3.north) {};
			\node at (1.north) {};
		\end{scope}
		\end{tikzpicture}

		\vskip0.5cm

		\begin{tikzpicture}
		\clip (-1.6*\length, 1.2*\height) rectangle (1.6*\length, 0.2);
		
		\begin{scope}[every node/.style={inner sep=5pt}]
			\node (0) at (0,0) {$[a_1|bc_0]$};
			\node (1) at (\length, \height) {$[bc_0]$};
			\node (-1) at (-\length, \height) {$[a_1]$};
		\end{scope}
		
		\node (2) at (2*\length,0) {\phantom{$[\,]$}};
		\coordinate (1a) at ($(1)!0.5!(2)$) {};
		
		\node (-2) at (-2*\length, 0) {\phantom{$[\,]$}};
		\coordinate (-1a) at ($(-1)!0.5!(-2)$) {};
		
		\draw (0.south) -> (1.north);
		\draw[dashed] (1.north) -> (1a);
		
		\draw (0.south) -> (-1.north);
		\draw[dashed] (-1.north) -> (-1a);
		
		\begin{scope}[every node/.style={fill=white, draw, circle, inner sep=1.4pt}]
		\node at (0.south) {};
		\end{scope}
		
		\begin{scope}[every node/.style={fill=black, draw, circle, inner sep=1.4pt}]
		\node at (1.north) {};
		\node at (-1.north) {};
		\end{scope}
		\end{tikzpicture}\qquad
		\begin{tikzpicture}
		\clip (-1.6*\length, 1.2*\height) rectangle (1.6*\length, 0.2);

		\begin{scope}[every node/.style={inner sep=5pt}]
			\node (0) at (0,0) {$[a_1b|c_0]$};
			\node (1) at (\length, \height) {$[c_0]$};
			\node (-1) at (-\length, \height) {$[a_1b]$};
		\end{scope}
		
		\node (2) at (2*\length,0) {\phantom{$[\,]$}};
		\coordinate (1a) at ($(1)!0.5!(2)$) {};
		
		\node (-2) at (-2*\length, 0) {\phantom{$[\,]$}};
		\coordinate (-1a) at ($(-1)!0.5!(-2)$) {};
		
		\draw (0.south) -> (1.north);
		\draw[dashed] (1.north) -> (1a);
		
		\draw (0.south) -> (-1.north);
		\draw[dashed] (-1.north) -> (-1a);
		
		\begin{scope}[every node/.style={fill=white, draw, circle, inner sep=1.4pt}]
		\node at (0.south) {};
		\end{scope}
		
		\begin{scope}[every node/.style={fill=black, draw, circle, inner sep=1.4pt}]
		\node at (1.north) {};
		\node at (-1.north) {};
		\end{scope}
		\end{tikzpicture}\qquad
		\begin{tikzpicture}
		\clip (-3.6*\length, 1.2*\height) rectangle (1.6*\length, 0.2);
		
		\begin{scope}[every node/.style={inner sep=5pt}]
			\node (0) at (0,0) {$[a_1c_0|a_{-1}]$};
			\node (1) at (\length, \height) {$[a_{-1}]$};
			\node (-1) at (-\length, \height) {$[a_1c_0]$};
			\node (-2) at (-2*\length, 0) {$[c_2|a_1c_0]$};
			\node (-3) at (-3*\length, \height) {$[c_2]$};
		\end{scope}
		
		\node (2) at (2*\length,0) {\phantom{$[\,]$}};
		\coordinate (1a) at ($(1)!0.5!(2)$) {};

		\node (-4) at (-4*\length, 0) {\phantom{$[\,]$}};
		\coordinate (-3a) at ($(-3)!0.5!(-4)$) {};
		
		\draw (0.south) -> (1.north);
		\draw[dashed] (1.north) -> (1a);
		
		\draw (0.south) -> (-1.north);
		\draw (-2.south) -> (-1.north);
		\draw (-2.south) -> (-3.north);
		\draw[dashed] (-3.north) -> (-3a);
		
		\begin{scope}[every node/.style={fill=white, draw, circle, inner sep=1.4pt}]
		\node at (0.south) {};
		\node at (-2.south) {};
		\end{scope}
		
		\begin{scope}[every node/.style={fill=black, draw, circle, inner sep=1.4pt}]
		\node at (1.north) {};
		\node at (-1.north) {};
		\node at (-3.north) {};
		\end{scope}
		\end{tikzpicture}

		\vskip0.5cm
		
		\begin{tikzpicture}
		\clip (-1.6*\length, 1.2*\height) rectangle (1.6*\length, 0.2);
		
		\begin{scope}[every node/.style={inner sep=5pt}]
			\node (0) at (0,0) {$[a_1|b|c_0]$};
			\node (1) at (\length, \height) {$[b|c_0]$};
			\node (-1) at (-\length, \height) {$[a_1|b]$};
		\end{scope}
		
		\node (2) at (2*\length,0) {\phantom{$[\,]$}};
		\coordinate (1a) at ($(1)!0.5!(2)$) {};

		\node (-2) at (-2*\length, 0) {\phantom{$[\,]$}};
		\coordinate (-1a) at ($(-1)!0.5!(-2)$) {};
		
		\draw (0.south) -> (1.north);
		\draw[dashed] (1.north) -> (1a);
		
		\draw (0.south) -> (-1.north);
		\draw[dashed] (-1.north) -> (-1a);
		
		\begin{scope}[every node/.style={fill=white, draw, circle, inner sep=1.4pt}]
		\node at (0.south) {};
		\end{scope}
		
		\begin{scope}[every node/.style={fill=black, draw, circle, inner sep=1.4pt}]
		\node at (1.north) {};
		\node at (-1.north) {};
		\end{scope}
		\end{tikzpicture}\qquad
		\begin{tikzpicture}
		\clip (-3.6*\length, 1.2*\height) rectangle (1.6*\length, 0.2);
		
		\begin{scope}[every node/.style={inner sep=5pt}]
			\node (0) at (0,0) {$[a_1|c_0|a_{-1}]$};
			\node (1) at (\length, \height) {$[c_0 | a_{-1}]$};
			\node (-1) at (-\length, \height) {$[a_1|c_0]$};
			\node (-2) at (-2*\length, 0) {$[c_2|a_1|c_0]$};
			\node (-3) at (-3*\length, \height) {$[c_2|a_1]$};
		\end{scope}
		
		\node (2) at (2*\length,0) {\phantom{$[\,]$}};
		\coordinate (1a) at ($(1)!0.5!(2)$) {};
		
		\node (-4) at (-4*\length, 0) {\phantom{$[\,]$}};
		\coordinate (-3a) at ($(-3)!0.5!(-4)$) {};
		
		\draw (0.south) -> (1.north);
		\draw[dashed] (1.north) -> (1a);
		
		\draw (0.south) -> (-1.north);
		\draw (-2.south) -> (-1.north);
		\draw (-2.south) -> (-3.north);
		\draw[dashed] (-3.north) -> (-3a);
		
		\begin{scope}[every node/.style={fill=white, draw, circle, inner sep=1.4pt}]
		\node at (0.south) {};
		\node at (-2.south) {};
		\end{scope}
		
		\begin{scope}[every node/.style={fill=black, draw, circle, inner sep=1.4pt}]
		\node at (1.north) {};
		\node at (-1.north) {};
		\node at (-3.north) {};
		\end{scope}
		\end{tikzpicture}

		\vskip0.5cm
		
		\begin{tikzpicture}
		\clip (-3.6*\length, 1.2*\height) rectangle (1.6*\length, 0.2);
		
		\begin{scope}[every node/.style={inner sep=5pt}]
			\node (0) at (0,0) {$[c_2|c_0|b']$};
			\node (1) at (\length, \height) {$[c_0 | b']$};
			\node (-1) at (-\length, \height) {$[c_2|c_0]$};
			\node (-2) at (-2*\length, 0) {$[b|c_2|c_0]$};
			\node (-3) at (-3*\length, \height) {$[b|c_2]$};
		\end{scope}
		
		\node (2) at (2*\length,0) {\phantom{$[\,]$}};
		\coordinate (1a) at ($(1)!0.5!(2)$) {};
		
		\node (-4) at (-4*\length, 0) {\phantom{$[\,]$}};
		\coordinate (-3a) at ($(-3)!0.5!(-4)$) {};
		
		\draw (0.south) -> (1.north);
		\draw[dashed] (1.north) -> (1a);
		
		\draw (0.south) -> (-1.north);
		\draw (-2.south) -> (-1.north);
		\draw (-2.south) -> (-3.north);
		\draw[dashed] (-3.north) -> (-3a);
		
		\begin{scope}[every node/.style={fill=white, draw, circle, inner sep=1.4pt}]
		\node at (0.south) {};
		\node at (-1.north) {};
		\node at (-2.south) {};
		\end{scope}
		
		\begin{scope}[every node/.style={fill=black, draw, circle, inner sep=1.4pt}]
		\node at (1.north) {};
		\node at (-3.north) {};
		\end{scope}
		\end{tikzpicture}\qquad
		\begin{tikzpicture}
		\clip (-3.6*\length, 1.2*\height) rectangle (1.6*\length, 0.2);
		
		\begin{scope}[every node/.style={inner sep=5pt}]
			\node (0) at (0,0) {$[a_1|a_{-1}|b]$};
			\node (1) at (\length, \height) {$[a_{-1} | b]$};
			\node (-1) at (-\length, \height) {$[a_1|a_{-1}]$};
			\node (-2) at (-2*\length, 0) {$[b'|a_1|a_{-1}]$};
			\node (-3) at (-3*\length, \height) {$[b'|a_1]$};
		\end{scope}
		
		\node (2) at (2*\length,0) {\phantom{$[\,]$}};
		\coordinate (1a) at ($(1)!0.5!(2)$) {};
		
		\node (-4) at (-4*\length, 0) {\phantom{$[\,]$}};
		\coordinate (-3a) at ($(-3)!0.5!(-4)$) {};
		
		\draw (0.south) -> (1.north);
		\draw[dashed] (1.north) -> (1a);
		
		\draw (0.south) -> (-1.north);
		\draw (-2.south) -> (-1.north);
		\draw (-2.south) -> (-3.north);
		\draw[dashed] (-3.north) -> (-3a);
		
		\begin{scope}[every node/.style={fill=white, draw, circle, inner sep=1.4pt}]
		\node at (0.south) {};
		\node at (-1.north) {};
		\node at (-2.south) {};
		\end{scope}
		
		\begin{scope}[every node/.style={fill=black, draw, circle, inner sep=1.4pt}]
		\node at (1.north) {};
		\node at (-3.north) {};
		\end{scope}
		\end{tikzpicture}
	\end{center}
	\caption{Les composantes fibrées de $K_W$ en type $\tilde{A}_2$.}
	\label{fig:fiber_components_A2tilde}
\end{figure}
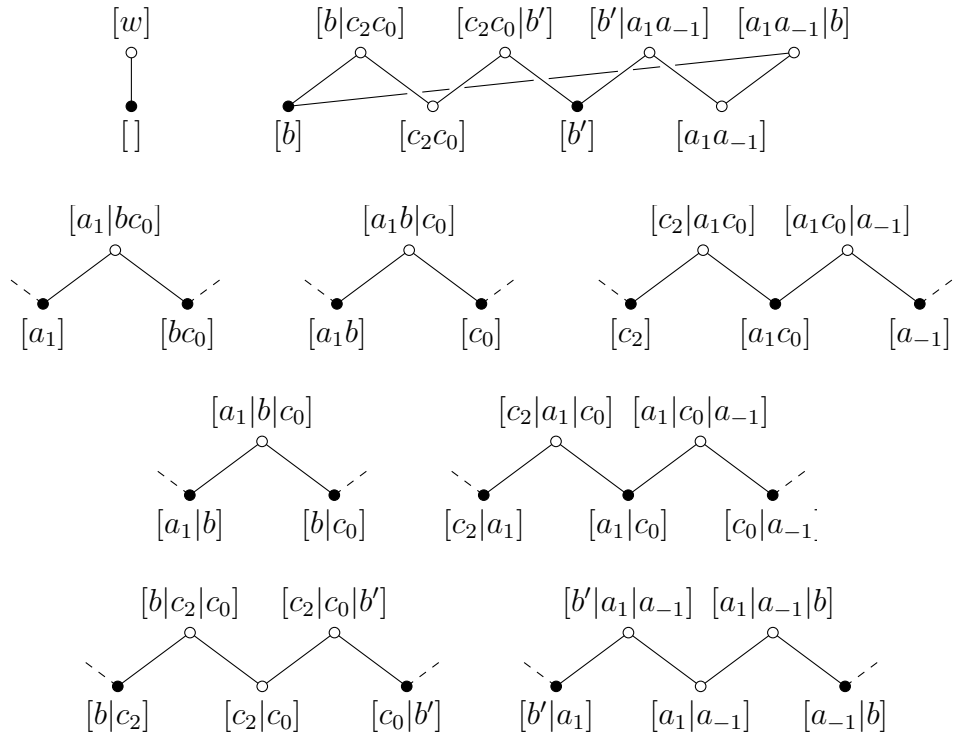

\bex \label{ex:fiber_components_A2tilde}
Nous allons présenter les composantes fibrées de $K_W$ dans le type $\tilde{A}_2$, en suivant \parencite[Example~7.12]{PS}. Nous utilisons les notations de l'exemple~\ref{ex:ordre_axial_A2tilde}, ainsi que la figure~\ref{fig:A2tilde_ordre_axial}. Le complexe $K_W$ possède 2 composantes finies, correspondant à la factorisation triviale de $w$, et aux factorisations
$$w=(b)(c_2c_0)=(c_2c_0)(b')=(b')(a_1a_{-1})=(a_1a_{-1})(b).$$
Ces deux composantes fibrées finies apparaissent en haut de la figure~\ref{fig:fiber_components_A2tilde}, qui provient de l'article~\parencite[Figure~8]{PS}. Les sept autres composantes sont infinies. Voici par exemple la suite infinie de factorisations de $w$ correspondant à la troisième composante de la figure~\ref{fig:fiber_components_A2tilde} :
$$ w = \cdots = (c_2b')(a_1) = (a_1)(bc_0) = (bc_0)(a_{-1}) =(a_{-1})(c_{-2}b') = (c_{-2}b')(a_{-3}) = \cdots.$$
Dans la figure~\ref{fig:fiber_components_A2tilde}, les sommets noirs correspondent aux simplexes du sous-complexe de Salvetti dual $X'_W$. Les sommets représentés correspondent à tous les simplexes du sous-complexe adapté fini $K'_W$. 
\eex

\section{La conjecture du $K(\pi,1)$}

\label{sec:proof_Kpi1}

Nous avons montré que $K_W$ est asphérique, que $K_W$ se rétracte sur $K'_W$, que $K'_W$ contient $X'_W$. Nous avons également montré que $X'_W$ a le type d'homotopie du complexe de Salvetti $X_W$, et donc également de l'espace de configuration $Y_W$.

Afin de montrer que $Y_W$ est asphérique, il suffit donc de montrer que $K'_W$ se rétracte sur $X'_W$. Pour cela, nous allons nous servir de l'ordre lexicographique construit dans la partie~\ref{sec:lexicographic_orderings} afin de construire un couplage sur les faces de $K'_W$ dont les faces critiques seront celles de $X'_W$.

\mk

Fixons un ordre axial $\prec$ sur l'ensemble des réflexions $R_0=R \cap [1,w]^W$, comme dans le théorème~\ref{thm:ordre_axial_EL}.

\mk

On peut aisément décrire les sous-complexes $K'_W$ et $X'_W$ de $K_W$ à l'aide des composantes fibrées et des applications $\lambda,\rho$. Soit $\sigma$ un simplexe de $K_W$.
\bit
\item $\sigma$ est un simplexe de $X'_W$ si et seulement si $\pi(\sigma)$ fixe un sommet de $C_0$.
\item $\sigma$ est un simplexe de $K'_W$ si et seulement s'il existe $k,k' \geq 0$ tels que $\lambda^k(\sigma) \in \mathcal{F}(X'_W)$ et $\rho^{k'}(\sigma) \in \mathcal{F}(X'_W)$.
\eit
En particulier, tous les simplexes $\sigma$ de $X'_W$ vérifient $\pi(\sigma) \neq w$, donc ils sont inférieurs.

\mk

Nous allons définir une notion supplémentaire, la profondeur d'un simplexe supérieur.

\bdf
Soit $\sigma[x_1| \dots |x_d]$ un simplexe de $K_W$ supérieur, i.e. tel que $\pi(\sigma)=w$. On appelle \emph{profondeur} de $\sigma$ le plus petit entier $\delta=\delta(\sigma) \in \{1,2,\dots,d\}$, tel que :
\bit
\item soit $\ell(x_\delta) \geq 2$,
\item soit $\ell(x_1) = \cdots = \ell(x_\delta)=1$ et, pour toute réflexion $r \in [1,w]^W$ telle que $r \leq x_{\delta+1}$, on a $x_\delta \prec r$. 
\eit
S'il n'existe aucun tel entier, on pose $\delta(\sigma)=\infty$.
\edf

Cette notion sera utile pour les simplexes supérieurs $\sigma$ dont le voisin de droite $\rho(\sigma)$ appartient à $X'_W$.

\blem
Soit $\sigma$ un simplexe de $K'_W$ supérieur tel que $\rho(\sigma)$ appartient à $X'_W$. Alors $\delta(\sigma) \neq \infty$.
\elem

Ceci va nous permettre de définir un couplage $\mu$ sur l'ensemble des faces de $K'_W$.

\bdf[Couplage]
Si $\sigma$ est un simplexe de $K'_W$ n'appartenant pas à $X'_W$, nous définissons un simplexe $\mu(\sigma)$ de $K_W$ comme suit.
\begin{enumerate}[(1)]
\item Si $\sigma$ est inférieur, on pose $\mu(\sigma)=\lambda(\sigma)$.
\item Si $\sigma$ est supérieur et $\rho(\sigma)$ n'appartient pas à $X'_W$, on pose $\mu(\sigma)=\rho(\sigma)$. 
\een
Supposons maintenant que $\sigma=[x_1| \dots |x_d]$ est supérieur et que $\rho(\sigma)$ appartient à $X'_W$.
\begin{enumerate}[(1)]
\setcounter{enumi}{2}
\item Si $\ell(x_\delta)\geq 2$, on pose $\mu(\sigma)=[x_1| \dots |x_{\delta-1}|y|z|x_{\delta+1}| \dots |x_d]$, où $x_\delta=yz$ et $y$ est la plus petite réflexion de $R_0 \cap [1,x_\delta]^W$ pour $\prec$.
\item Si $\ell(x_\delta) = 1$, on pose $\mu(\sigma)=[x_1| \dots |x_{\delta-1}|x_\delta x_{\delta+1}| x_{\delta+2}| \dots |x_d]$.
\een
\edf

Remarquons que, dans les cas (1) et (3), $\sigma$ est une facette de $\mu(\sigma)$, tandis que dans les cas (2) et (4), c'est $\mu(\sigma)$ qui est une facette de $\sigma$.

\bex
Nous reprenons l'exemple~\ref{ex:fiber_components_A2tilde} du type $\tilde{A}_2$, accompagné de la figure~\ref{fig:fiber_components_A2tilde}. Dans ce cas, le couplage $\mu$ est un couplage entre les sommets blancs de la figure~\ref{fig:fiber_components_A2tilde} correspondant aux simplexes de $K'_W \smallsetminus X'_W$. Le couplage $\mu$ est entièrement décrit comme suit \parencite[Figure~9]{PS} :
\beq
\begin{aligned}[c]
[c_2c_0] & \longleftrightarrow & [b|c_2 c_0] \\
\,[a_1 a_{-1}] & \longleftrightarrow & [b'|a_1 a_{-1}] \\
\,[c_2|c_0] & \longleftrightarrow & [b|c_2|c_0] \\
\,[a_1|a_{-1}] & \longleftrightarrow & [b'|a_1|a_{-1}] \\
\,[w] & \longleftrightarrow & [a_1|bc_0]
\end{aligned}
\qquad
\begin{aligned}[c]
[c_2c_0|b']  & \longleftrightarrow & [c_2|c_0|b'] \\
\, [a_1a_{-1}|b] & \longleftrightarrow & [a_1|a_{-1}|b] \\
\, [c_2|a_1c_0] & \longleftrightarrow & [c_2|a_1|c_0] \\
\, [a_1c_0|a_{-1}] & \longleftrightarrow & [a_1|c_0|a_{-1}] \\
\, [a_1b|c_0] & \longleftrightarrow & [a_1|b|c_0]
\end{aligned}
\eeq
Par exemple, si l'on considère le simplexe $\sigma=[w]$, alors $\ell(\sigma)=1$ et on est dans le cas (3) de la définition de $\mu(\sigma)$. Ainsi, comme on peut le voir dans l'exemple~\ref{ex:ordre_axial_A2tilde}, la plus petite réflexion de $R_0$ pour $\prec$ est $a_1$, donc $\mu([w])=[a_1|bc_0]$. Réciproquement, comme $a_1$ est inférieure pour $\prec$ à toutes les réflexions inférieures à $bc_0$, on a $\delta([a_1|bc_0])=1$ donc $\mu([a_1|bc_0])=[w]$.
\eex

Cette application $\mu$ préserve les faces de $K'_W$ :

\blem[{\cite[Lemma~8.8]{PS}}] \label{lem:composition_mu}
Si $\sigma \in \mathcal{F}(K'_W) \smallsetminus \mathcal{F}(X'_W)$, alors $\mu(\sigma) \in \mathcal{F}(K'_W) \smallsetminus \mathcal{F}(X'_W)$.
\elem

\bp
Nous allons simplement donner les idées de la preuve de ce lemme. Soit $\sigma \in \mathcal{F}(K'_W) \smallsetminus \mathcal{F}(X'_W)$, distinguons selon les cas de la définition de $\mu(\sigma)$ :

\begin{enumerate}[(1)]
\item Si $\sigma$ est inférieur, alors $\mu(\sigma)=\lambda(\sigma) \in \mathcal{F}(K'_W)$. De plus, comme $\mu(\sigma)$ est supérieur et que $X'_W$ est constitué de simplexes inférieurs, on déduit que $\mu(\sigma) \not\in \mathcal{F}(X'_W)$.
\item Si $\sigma$ est supérieur et $\rho(\sigma)$ n'appartient pas à $X'_W$, alors $\mu(\sigma)=\rho(\sigma) \in \mathcal{F}(K'_W) \smallsetminus \mathcal{F}(X'_W)$ par hypothèse.
\een
Supposons maintenant que $\sigma=[x_1| \dots |x_d]$ est supérieur et que $\rho(\sigma)$ appartient à $X'_W$. On peut montrer que $\rho(\mu(\sigma)) \in \mathcal{F}(X'_W)$. Il suffit donc de montrer qu'il y a un simplexe de $X'_W$ à gauche de $\mu(\sigma)$. Comme $\sigma$ appartient à $K'_W$, il existe un simplexe $\tau \in \mathcal{F}(X'_W)$ à gauche de $\sigma$, de la forme
$$\tau=[\phi^h(x_{i+1})| \dots |\phi^h(x_d)|\phi^{h+1}(x_1)| \dots |\phi^{h+1}(x_{i-1})],$$
où $h<0$ et $1 \leq i \leq d$, et où $\phi$ désigne la conjugaison par $w$.
\begin{enumerate}[(1)]
\setcounter{enumi}{2}
\item Si $\ell(x_\delta)\geq 2$, alors $\mu(\sigma)=[x_1| \dots |x_{\delta-1}|y|z|x_{\delta+1}| \dots |x_d]$. Si $i>\delta$, considérons le simplexe
$$\tau'=[\phi^h(x_{i+1})| \dots |\phi^h(x_d)|\phi^{h+1}(x_1)| \dots |\phi^{h+1}(y)|\phi^{h+1}(z)| \dots |\phi^{h+1}(x_{i-1})].$$
C'est un simplexe de $X'_W$ qui est à gauche de $\mu(\sigma)$, ainsi $\mu(\sigma) \in \mathcal{F}(K'_W)$. Le cas $i<\delta$ est similaire. Le cas $i=\delta$ est nettement plus technique, et nous l'admettons.
\item Si $\ell(x_\delta) = 1$, alors $\mu(\sigma)=[x_1| \dots |x_{\delta-1}|x_\delta x_{\delta+1}|x_{\delta+2}| \dots |x_d]$. Supposons par exemple que $i>\delta+1$, considérons alors le simplexe
$$\tau'=[\phi^h(x_{i+1})| \dots |\phi^h(x_d)|\phi^{h+1}(x_1)| \dots |\phi^{h+1}(x_\delta x_{\delta+1})| \dots |\phi^{h+1}(x_{i-1})],$$
c'est un simplexe de $X'_W$ qui est à gauche de $\mu(\sigma)$, ainsi $\mu(\sigma) \in \mathcal{F}(K'_W)$. Les autres valeurs de $i$ se traitent de manière similaire.
\een
\ep

Ceci nous permet de montrer que $\mu$ est une involution :

\bpro \label{pro:mu_involution_sans_point_fixe}
L'application $\mu$ est une involution de $\mathcal{F}(K'_W) \smallsetminus \mathcal{F}(X'_W)$ sans point fixe.
\epro

\bp
D'après le lemme~\ref{lem:composition_mu}, on a le droit de composer $\mu$ avec elle-même. Soit $\sigma \in \mathcal{F}(K'_W) \smallsetminus \mathcal{F}(X'_W)$, distinguons selon les cas de la définition de $\mu(\sigma)$ :
\begin{enumerate}[(1)]
\item Si $\sigma$ est inférieur, alors $\mu(\sigma)=\lambda(\sigma)$ est supérieur et $\rho(\mu(\sigma))=\sigma$ n'appartient pas à $X'_W$, donc $\mu(\mu(\sigma))=\rho(\lambda(\sigma))$.
\item Si $\sigma$ est supérieur et $\rho(\sigma)$ n'appartient pas à $X'_W$, comme $\mu(\sigma)=\rho(\sigma)$ est inférieur, on a $\mu(\mu(\sigma))=\lambda(\rho(\sigma))=\sigma$.
\een
Supposons maintenant que $\sigma=[x_1| \dots |x_d]$ est supérieur et que $\rho(\sigma)$ appartient à $X'_W$. Par définition de $\delta$ nous savons que $x_1 \succ x_2 \succ \cdots \succ x_{\delta-1} \succ y$, où $y$ désigne la réflexion minimale de $[1,x_\delta]$ pour l'ordre $\prec$.
\begin{enumerate}[(1)]
\setcounter{enumi}{2}
\item Si $\ell(x_\delta)\geq 2$, alors $x_\delta=yz$, et pour toute réflexion $r \leq z$ nous savons que $y \prec r$. Ainsi $\delta(\mu(\sigma))=\delta$, et $\mu(\mu(\sigma))=\sigma$.
\item Si $\ell(x_\delta) = 1$, alors $x_{\delta-1} \succ x_\delta$, ainsi $\delta(\mu(\sigma))=\delta$, et $x_\delta$ est la plus petite réflexion de $[1,x_\delta x_{\delta+1}]$ pour $\prec$, ainsi $\mu(\mu(\sigma))=\sigma$.
\een
\ep

Nous allons utiliser la théorie de Morse discrète, et nous renvoyons à la partie~\ref{sec:discrete_morse_theory} pour plus de détails. L'involution $\mu$ permet de définir un couplage $\mathcal{M}$ sur l'ensemble $\mathcal{F}(K'_W)$ des faces de $K'_W$ :
$$ \mathcal{M} = \{(\mu(\sigma),\sigma) \st \sigma \in \mathcal{F}(K'_W) \smallsetminus \mathcal{F}(X'_W) \mbox{ et } \mu(\sigma) \lessdot \sigma\},$$
où on rappelle que l'on note $\tau \lessdot \sigma$ si $\tau$ est une facette (face de codimension $1$) de $\sigma$. Comme $\mu$ est sans point fixe, l'ensemble des faces critiques de $\mathcal{M}$ est $\mathcal{F}(X'_W)$.

\bpro[{\cite[Lemma~8.13]{PS}}]
Le couplage $\mathcal{M}$ est acyclique et propre.
\epro

\bp
Comme $\mathcal{F}(K'_W)$ est fini, le couplage $\mathcal{M}$ est propre. Nous allons donner les grandes lignes de la preuve de l'acyclicité de $\mathcal{M}$.

\mk

Pour cela, nous allons définir une application $\xi$ de $\mathcal{F}(K'_W) \smallsetminus \mathcal{F}(X'_W)$ à valeurs dans un ensemble totalement ordonné $(P,\trianglelefteq)$, qui décroît le long de chemin orientés pour $\mathcal{M}$. Ceci permet de montrer qu'il n'y a pas de cycle orienté pour $\mathcal{M}$.

\mk

Considérons l'ensemble $P \subset {R_0}^{n+1}$ des factorisations minimales de $w$ comme produit de $(n+1)$ réflexions. L'ordre total $\trianglelefteq$ sur $P$ est défini de la manière suivante. Soient $\alpha,\alpha' \in P$, notons $r,r' \in R_0$ les réflexions maximales pour $\prec$ apparaissant dans $\alpha,\alpha'$ respectivement, ainsi que $1 \leq k,k' \leq n+1$ leurs positions d'apparition.
\bit
\item Si $r \neq r'$, alors $\alpha \triangleleft \alpha'$ si et seulement si $r \succ r'$.
\item Si $r=r'$ et $k \neq k'$, alors $\alpha \triangleleft \alpha'$ si et seulement si $k > k'$. 
\item Si $r=r'$ et $k = k'$, alors $\alpha \triangleleft \alpha'$ si et seulement si $\alpha$ vient avant $\alpha'$ dans l'ordre lexicographique pour $\prec$.
\eit

\mk

Nous allons maintenant définir une application $\xi \colon \mathcal{F}(K'_W) \smallsetminus \mathcal{F}(X'_W) \ra P$ comme suit. Soit $\sigma \in \mathcal{F}(K'_W) \smallsetminus \mathcal{F}(X'_W)$.
\bit
\item Si $\sigma=[x_1| \dots |x_d]$ est supérieur, i.e. $\pi(\sigma)=w$, définissons $\xi(\sigma) \in P$ comme la concaténation des factorisations croissantes de $x_1,\dots,x_d$.
\item Si $\sigma$ est inférieur, considérons $\mu(\sigma)=\lambda(\sigma)=[x_1| \dots |x_d]$, et définissons $\xi(\sigma) \in P$ comme la concaténation des factorisations croissantes de $x_1,\dots,x_d$.
\eit

\mk

La preuve de l'acyclicité repose sur deux faits, précisant le comportement de l'application $\xi$ lorsqu'on passe d'un simplexe à l'une de ses facettes.

{\bf Fait 1} Si $\sigma \in \mathcal{F}(K'_W) \smallsetminus \mathcal{F}(X'_W)$, alors $\xi(\mu(\sigma))=\xi(\sigma)$.

Quitte à échanger les rôles de $\sigma$ et $\mu(\sigma)$, on peut supposer que $\mu(\sigma) \lessdot \sigma$. Distinguons selon que $\sigma$ est dans le cas (2) ou (4) de la définition de $\mu(\sigma)$ :
\ben
\item[(2)] Si $\sigma$ est supérieur et $\rho(\sigma)$ n'appartient pas à $X'_W$, alors $\mu(\sigma)$ est inférieur et $\mu(\mu(\sigma))=\sigma$, donc $\xi(\mu(\sigma))=\xi(\sigma)$.
\item[(4)] Si $\sigma=[x_1| \dots |x_d]$ est supérieur, que $\rho(\sigma)$ appartient à $X'_W$ et $\ell(x_\delta)=1$, alors $\mu(\sigma)=[x_1| \dots |x_{\delta-1}|x_\delta x_{\delta+1}| x_{\delta+2}| \dots |x_d]$. Comme $x_\delta$ est la réflexion minimale de $[1,x_\delta x_{\delta+1}]$ pour $\prec$, on déduit que la factorisation croissante de $x_\delta x_{\delta+1}$ est la concaténation de $x_\delta$ avec la factorisation croissante de $x_{\delta+1}$. Ainsi $\xi(\mu(\sigma))=\xi(\sigma)$.
\een

\mk

{\bf Fait 2} Soient $\sigma,\tau \in \mathcal{F}(K'_W) \smallsetminus \mathcal{F}(X'_W)$ deux simplexes tels que $\sigma$ est supérieur et $\tau$ est une facette de $\sigma$. Alors $\xi(\tau) \trianglelefteq \xi(\sigma)$. Si de plus $\tau=\lambda(\sigma)$, alors $\xi(\tau) \triangleleft \xi(\sigma)$.

\mk

Nous admettons ce deuxième fait, et nous allons montrer que ces deux faits impliquent l'acyclicité de $\mathcal{M}$. Supposons par l'absurde qu'il existe un cycle orienté de simplexes distincts
$$\sigma_1 \gtrdot \tau_1 \lessdot \sigma_2 \gtrdot \tau_2 \lessdot \cdots \gtrdot \tau_m \lessdot \sigma_{m+1}=\sigma_1$$
dans $\mathcal{F}(K'_W)$, où $\mu(\tau_j)=\sigma_{j+1}$ pour tout $1 \leq j \leq m$.

\mk

D'après les faits 1 et 2, on déduit que
$$\xi(\sigma_1) \trianglerighteq \xi(\tau_1) = \xi(\sigma_2) \trianglerighteq \cdots \trianglerighteq \xi(\tau_m) = \xi(\sigma_{m+1}) =\xi(\sigma_1),$$
ainsi toutes ces inégalités sont des égalités.

\mk

Pour tout $1 \leq i \leq d$, comme $\mu(\tau_i) \neq \sigma_i$, on déduit que $\tau_i \neq \rho(\sigma_i)$. D'après le fait 2, on sait que $\tau_i \neq \lambda(\sigma_i)$. Ainsi $\tau_i$ est une facette de $\sigma_i$ différente de $\rho(\sigma_i)$ et de $\lambda(\sigma_i)$, ainsi tous les $\sigma_i$ et les $\tau_i$ sont des simplexes supérieurs.

\mk

De plus, pour tout $1 \leq i \leq d$, comme $\xi(\sigma_i)=\xi(\tau_i)$, nous avons $\delta(\sigma_i) \geq \delta(\tau_i)$. D'autre part, d'après la preuve de la Proposition~\ref{pro:mu_involution_sans_point_fixe}, nous savons que $\delta(\tau_i)=\delta(\mu(\tau_i))=\delta(\sigma_{i+1})$. Ainsi $\delta(\sigma_1)=\delta(\tau_1)= \cdots =\delta(\sigma_m)=\delta(\tau_m)$.

\mk

\'{E}crivons donc $\sigma_1=[x_1|x_2| \dots |x_d]$ et $\tau_1=[x_1| \dots |x_{i-1}|x_ix_{i+1}|x_{i+2}| \dots |x_d]$, où $1 \leq i \leq d-1$, et notons $\delta=\delta(\sigma_1)=\delta(\tau_1)$. De plus, $\sigma_1$ est dans le cas (4) de la définition de $\mu(\sigma_1)$, donc $\ell(x_\delta)=1$. Comme $\ell(x_ix_{i+1}) \geq 2$, nous savons que $\delta \leq i$. Comme $\tau_1$ est dans le cas (3) de la définition de $\mu(\tau_1)$, nous savons que $\delta \geq i$. Ainsi $\delta=i$, d'où $\sigma_1=\mu(\tau_1)=\sigma_2$, ce qui est
une contradiction.
\ep

Ceci permet enfin de donner la preuve de la conjecture du $K(\pi,1)$ pour tous les groupes d'Artin affines.

\bthm[{La conjecture du $K(\pi,1)$; \cite[Theorem~8.15]{PS}}]
Soit $W$ un groupe de Coxeter affine. Alors la conjecture du $K(\pi,1)$ pour le groupe d'Artin $G_W$ est vraie : l'espace de configuration $Y_W$ est un espace classifiant pour $G_W$.
\ethm

\bp
Il suffit de consider le cas où $W$ est irréductible. Fixons un ensemble simple de réflexions $S$, et un élément de Coxeter $w$. D'après le théorème~\ref{thm:artin_dual_classifiant}, le complexe d'intervalle $K_W$ est un espace classifiant (pour le groupe d'Artin dual $W_w$). 

\mk

D'après le théorème~\ref{thm:KW_retracte_KprimeW}, le complexe $K_W$ se rétracte par déformation forte sur son sous-complexe adapté fini $K'_W$.

\mk

Comme le couplage $\mathcal{M}$ sur $\mathcal{F}(K'_W)$ est acyclique et que son ensemble de faces critiques est $\mathcal{F}(X'_W)$, on déduit du théorème~\ref{thm:discrete_morse} que $K'_W$ se rétracte par déformation forte sur~$X'_W$.

\mk

D'après le théorème~\ref{thm:YW_XprimeW}, le complexe $X'_W$ a le même type d'homotopie que le complexe de Salvetti $X_W$ et que l'espace de configuration $Y_W$. On conclut que $Y_W$ est un espace classifiant pour son groupe fondamental $G_W$.
\ep

De plus, ceci fournit une nouvelle preuve de l'isomorphisme entre le groupe d'Artin~$G_W$ et le groupe d'Artin dual~$W_w$, dû à McCammond et Sulway :

\bthm[{\cite[Theorem~C]{mccammond_sulway}; \cite[Theorem~8.16]{PS}}]
Soit $W$ un groupe de Coxeter affine irréductible, et soit $w$ un élément de Coxeter. Le morphisme naturel du groupe d'Artin $G_W$ vers le groupe d'Artin dual $W_w$ est un isomorphisme. 
\ethm

\bp
La preuve ci-dessus montre déjà que $G_W$ et $W_w$ sont isomorphes. Pour toute réflexion simple $s \in S$, remarquons que l'arête $[s]$ représentant $s \in \pi_1(K_W)=W_w$ appartient au sous-complexe $X'_W$. De plus, l'équivalence d'homotopie entre $X_W$ et $X'_W$ du théorème~\ref{thm:YW_XprimeW} identifie cette arête $[s]$ avec le sous-complexe $X_{\{s\}}$ du complexe de Salvetti $X_W$. Ainsi cette arête représente l'élément $s \in \pi_1(X_W)=G_W$. Ceci montre que l'isomorphisme entre $G_W$ et $W_w$ donné par l'équivalence d'homotopie entre $X_W$ et $K_W$ est bien le morphisme naturel.
\ep

\section{Théorie de Morse discrète}

\label{sec:discrete_morse_theory}

Nous présentons ici quelques éléments de théorie de Morse discrète, qui ont permis de montrer que le complexe $K_W$ se rétracte sur le sous-complexe $K'_W$, puis sur le sous-complexe $X'_W$. Cette théorie est dûe à 
\textcite{forman_guide,forman}, et nous utilisons le point de vue de 
\textcite{chari} présenté dans \textcite{PS}.

\mk

Soit $P$ un ensemble gradué. Si $p,q \in P$, rappelons que l'on note $p \lessdot q$ si $p < q$ et il n'y aucun élément $r \in P$ tel que $p < r < q$. Notons $H$ le diagramme de Hasse de $P$ : rappelons que c'est le graphe de sommets $P$, avec une arête entre $p$ et $q$ si $p \lessdot q$ ou $q \lessdot p$. Notons $\mathcal{E}(P)$ l'ensemble des arêtes de $H$.

\mk

Si $\mathcal{M}$ est un sous-ensemble de $\mathcal{E}(P)$, on peut orienter les arêtes de $H$ de la manière suivante : une arête entre $p$ et $q$ tels que $p \lessdot q$ est orientée de $p$ vers $q$ si elle appartient à $\mathcal{M}$, et de $q$ vers $p$ sinon. Notons $H_\mathcal{M}$ le graphe orienté ainsi obtenu.

\bdf\

Un \emph{couplage} sur $P$ est un sous-ensemble $\mathcal{M} \subset \mathcal{E}(P)$ tel que tout élément de $P$ appartient à au plus une arête de $\mathcal{M}$.

Le couplage est \emph{acyclique} si le graphe $H_\mathcal{M}$ n'a pas de cycle orienté.

Le couplage est \emph{propre} si, pour tout $p \in P$, l'ensemble des chemins orientés dans $H_\mathcal{M}$ issus de $p$ est fini.

Un élément $p \in P$ est appelé \emph{critique} s'il n'appartient à aucune arête de $\mathcal{M}$.
\edf

La notion de couplage est pertinente pour les ensembles ordonnés de faces d'un complexe. Soit~$X$ un CW-complexe. L'ensemble $\mathcal{F}(X)$ des faces (cellules ouvertes) de~$X$ est ordonné par la relation suivante : $\tau \leq \sigma$ si et seulement si $\ov{\tau} \subset \ov{\sigma}$.

Fixons une cellule $\sigma$ de $X$ de dimension $n \geq 1$, et considérons l'application $\phi \colon \D^n \ra X$ définissant $\sigma$. Si $\tau$~est de codimension $1$ dans $\sigma$, on dit que $\tau$~est une \emph{facette régulière} de~$\sigma$ si :
\bit
\item l'application $\phi$ se restreint en un homéomorphisme de $\phi^{-1}(\tau)$ sur $\tau$, et
\item l'adhérence $\ov{\phi^{-1}(\tau)} \subset \D^n$ est une boule fermée de dimension $n-1$.
\eit

\mk

Voici le théorème principal de la théorie de Morse discrète.

\bthm[\cite{batzies,chari,forman}] \label{thm:discrete_morse}
Soit $X$ un CW-complexe, et $Y \subset X$ un sous-complexe. Supposons qu'il existe un couplage acyclique propre $\mathcal{M}$ sur l'ensemble ordonné $\mathcal{F}(X)$ des cellules de $X$ tel que :
\bit
\item l'ensemble des cellules critiques de $\mathcal{M}$ est $\mathcal{F}(Y)$ et
\item pour toute paire $(\tau,\sigma) \in \mathcal{M}$, la cellule $\tau$ est une facette régulière de $\sigma$.
\eit
Alors $X$ se rétracte par déformation forte sur $Y$. En particulier, l'inclusion $Y \hookrightarrow X$ est une équivalence d'homotopie.
\ethm

Voici un outil classique pour construire des couplages acycliques, appelé théorème du patchwork :

\bthm[{\cite[Theorem~11.10]{kozlov}}] \label{thm:patchwork}
Soit $\eta \colon P \ra Q$ une application d'ensembles ordonnés. Supposons que, pour tout $q \in Q$, nous ayons un couplage acyclique $\mathcal{M}_q \subset \mathcal{E}(P)$ ne comportant que des éléments de la fibre $\eta^{-1}(q)$. Alors la réunion de ces couplages est un couplage acyclique sur $P$.
\ethm

\printbibliography

\end{document}
